\input amstex
\input amsppt.sty
\magnification=\magstep1
\hsize=30truecc
\vsize=22.2truecm
\baselineskip=16truept
\NoBlackBoxes
\TagsOnRight \pageno=1 \nologo
\def\Z{\Bbb Z}
\def\N{\Bbb N}

\def\Q{\Bbb Q}

\def\l{\left}
\def\r{\right}
\def\bg{\bigg}
\def\({\bg(}
\def\[{\bg\lfloor}
\def\){\bg)}
\def\]{\bg\rfloor}
\def\t{\text}
\def\f{\frac}

\def\p{\ (\roman{mod}\ p)}

\def\sm{\setminus}

\def\bi{\binom}
\def\eq{\equiv}

\def\ls{\leqslant}
\def\gs{\geqslant}
\def\mo{\roman{mod}}

\def\ve{\varepsilon}

\def\da{\delta}

\def\M#1#2{\thickfracwithdelims[]\thickness0{#1}{#2}_q}

\def\Remark{\noindent{\it  Remark}}

\hbox {{\tt arXiv:0911.5665}}
\bigskip
\topmatter
\title Open Conjectures on congruences \endtitle
\author Zhi-Wei Sun\endauthor
\leftheadtext{Zhi-Wei Sun} \rightheadtext{Open  Conjectures on congruences}
\affil Department of Mathematics, Nanjing University\\
 Nanjing 210093, People's Republic of China
  \\  zwsun\@nju.edu.cn
  \\ {\tt http://math.nju.edu.cn/$\sim$zwsun}
\endaffil
\abstract We collect here various conjectures on congruences made by
the author in a series of papers, some of which involve binary
quadratic forms and other advanced theories. Part A consists of 100
unsolved conjectures of the author while conjectures in Part B have
been recently confirmed. We hope that this material will interest
number theorists and stimulate further research. Number theorists
are welcome to work on those open conjectures; for some of them we offer
prizes for the first correct proofs.
\endabstract
\thanks 2010 {\it Mathematics Subject Classification}. Primary 11B65, 11A07;
Secondary 05A10, 05A15, 11A41, 11B39, 11B68, 11E25, 11S99, 33C20, 65B10.
\newline\indent {\it Keywords}. Binomial coefficients, Catalan numbers,
Bernoulli numbers, Euler numbers, series involving $\pi$, binary quadratic forms, congruences modulo prime powers.
\newline\indent \copyright\  Copyright is owned by the author Zhi-Wei Sun. The material
on the author's homepage has been linked to Number Theory Web since
Nov. 27, 2009.
\endthanks
\endtopmatter
\document

\heading{Introduction}\endheading

Congruences modulo primes have been widely investigated since the
time of Fermat. However, we find that there are still lots of new
 challenging congruences that cannot be easily solved. They appeal
 for new powerful tools or advanced theory.

 Here we collect various conjectures of the author on congruences, most of which can be found in the author's papers
 available from {\tt arxiv} or his homepage. We use two sections to state conjectures and related remarks.
 Part A contains 100 unsolved conjectures of the author while Part B consists of conjectures that have
been recently confirmed. Most of the congruences here
are {\it super} congruences in the sense that they happen to hold modulo some higher power of $p$.
The topic of super congruences is related to the $p$-adic $\Gamma$-function, Gauss and Jacobi sums,
hypergeometric series, modular forms,
Calabi-Yau manifolds, and some sophisticated combinatorial identities involving harmonic numbers (cf. K. Ono's book [O]).
The recent theory of super congruences also involves Bernoulli and Euler numbers (see [S11b, S11e])
and various series related to $\pi$ (cf. [vH], [S-7] and [T1]).
Many congruences collected here are about $\sum_{k=0}^{p-1}a_k/m^k$ modulo powers of a prime $p$, where
$m$ is an integer not divisible by $p$ and  the quantity $a_k$ is a sum or a product
of some binomial coefficients which usually arises from enumerative combinatorics.

 For clarity, we often state the prime version of a conjecture instead of the general version.

 We also include several challenging conjectures on series related to $\pi$ or the Dirichlet $L$-function which arose
 from the author's investigations of congruences, see Conjectures A3, A5, A40, A44.

 For some conjectures we announce prizes for the first correct proofs, see remarks after Conjectures A1, A39, A44, A46, A48.

 Now we introduce some basic notation in this paper.

 As usual, we set
 $$\N=\{0,1,2,\ldots\}\ \ \t{and}\ \ \Z^+=\{1,2,3,\ldots\}.$$
 The Kronecker symbol $\da_{m,n}$ takes 1 or 0 according as $m=n$ or not.
 The rising factorial $(x)_n$ is defined by
 $(x)_n=\prod_{k=0}^{n-1}(x+k),$
 and $(x)_0$ is regarded as 1.
  For an integer $m$ and a positive odd number
 $n$, the notation $(\f mn)$ stands for the Jacobi symbol.
 For an odd prime $p$, we use $q_p(2)$ to denote the Fermat quotient $(2^{p-1}-1)/p$.
 For a prime $p$ and a rational number $x$, the {\it $p$-adic
 valuation}
 of $x$ is given by
 $$\nu_p(x)=\sup\{a\in\N:\ x\eq0\ (\mo\ p^a)\}.$$
 For a polynomial or a power series $P(x)$, we write $[x^n]P(x)$ for the coefficient of $x^n$ in the expansion of $P(x)$.
 For $k_1,\ldots,k_n\in\N$, we define the multinomial coefficient
 $$\bi{k_1+\cdots+k_n}{k_1,\ldots,k_n}:=\f{(k_1+\cdots+k_n)!}{k_1!\cdots k_n!}.$$

 Harmonic numbers are given by
 $$H_0=0 \ \ \t{and}\ \ H_n=\sum_{k=1}^n\f1k\ \ (n=1,2,3,\ldots).$$
 For $n\in\N$, $C_n$ denotes the {\it Catalan number} $\f1{n+1}\bi{2n}n=\bi{2n}n-\bi{2n}{n+1}$
 and $C_n^{(2)}$ stands for the (first kind) {\it second-order Catalan number} $\f1{2n+1}\bi{3n}n=\bi{3n}n-2\bi{3n}{n-1}$.
 Note that if $p$ is an odd prime then
 $$\bi{2k}k=\f{(2k)!}{(k!)^2}\eq 0\ (\mo\ p)\quad\t{for every}\ k=\f{p+1}2,\ldots,p-1.$$

 Bernoulli numbers $B_0,B_1,B_2,\ldots$ are rational numbers given by
$$B_0=1\ \ \t{and}\ \ \sum^n_{k=0}\bi {n+1}k B_k=0\ \ \ \t{for}\ n\in\Z^+=\{1,2,3,\ldots\}.$$
It is well known that $B_{2n+1}=0$ for all $n\in\Z^+$ and
$$\f x{e^x-1}=\sum_{n=0}^\infty B_n\f{x^{n}}{n!}\ \ \l(|x|<2\pi\r).$$
Euler numbers $E_0,E_1,E_2,\ldots$ are integers defined by
$$E_0=1\ \ \t{and}\ \ \sum^n\Sb k=0\\2\mid k\endSb \bi nk E_{n-k}=0\ \ \ \t{for}\ n\in\Z^+=\{1,2,3,\ldots\}.$$
It is well known that $E_{2n+1}=0$ for all $n\in\N$ and
$$\sec x=\sum_{n=0}^\infty(-1)^n E_{2n}\f{x^{2n}}{(2n)!}\ \ \l(|x|<\f{\pi}2\r).$$
Bernoulli polynomials and Euler polynomials are given by
$$B_n(x)=\sum_{k=0}^n\bi nk B_kx^{n-k}\ \, \t{and}\ \, E_n(x)=\sum_{k=0}^n\bi{n}{k}\frac{E_k}{2^k}\left(x-\frac{1}{2}\right)^{n-k}\
(n\in\N).$$

For $A,B\in\Z$ we define the Lucas sequences $u_n=u_n(A,B)\ (n\in\N)$ and $v_n=v_n(A,B)\ (n\in\N)$ as follows:
$$u_0=0,\ u_1=1,\ \t{and}\ u_{n+1}=Au_n-Bu_{n-1}\ (n=1,2,3,\ldots);$$
$$v_0=0,\ v_1=1,\ \t{and}\ v_{n+1}=Av_n-Bv_{n-1}\ (n=1,2,3,\ldots).$$

\medskip

 \heading{Part A. Conjectures that remain unsolved}\endheading

\proclaim{Conjecture A1 {\rm ([S11b])}}
Let $p$ be an odd prime. Then
$$\align&\sum_{k=0}^{p-1}\bi{2k}k^3
\\\eq&\cases4x^2-2p\ (\mo\ p^2)&\t{if}\ (\f p7)=1\ \&\ p=x^2+7y^2\ \t{with}\ x,y\in\Z,
\\0\ (\mo\ p^2)&\t{if}\ (\f p7)=-1,\ \t{i.e.,}\ p\eq3,5,6\ (\mo\ 7).\endcases
\endalign$$
\endproclaim

\Remark. (a) It is well known
that  $p=x^2+7y^2$ for some $x,y\in\Z$ if $p$ is an odd prime with $(\f p7)=1$.
The congruence modulo $p$ can be easily deduced from (6) in Ahlgren [A, Theorem 5].
Recently the author's twin brother Z. H. Sun [Su2] confirmed the conjecture in the case $(\f p7)=-1$.
I'd like to offer \$70 (US dollars) for the first correct proof of Conj. A1 in the case $(\f p7)=1$.

(b) Let $p$ be an odd prime. By [S11e, Theorem 1.3] or [S-7, Theorem 1.5], Conjecture A1 implies that
$$\align&3\sum_{k=0}^{p-1}k\bi{2k}k^3
\\\eq& \cases\f 87(3p-4x^2)=32y^2-\f 87p\ (\mo\ p^2)&\t{if}\ (\f p7)=1\ \&\ p=x^2+7y^2\ (x,y\in\Z),
\\\f 87p\ (\mo\ p^2)&\t{if}\ (\f p7)=-1.\endcases\endalign$$
We ever wrote that
the author was unable to guess $\sum_{k=0}^{p-1}k\bi{2k}k^3$ mod $p$ in the case $(\f p7)=1$
though we formulated Conj. A1 on November 13, 2009.
After reading this remark, on Nov. 28, 2009 Bilgin Ali and Bruno
Mishutka guessed that for $p=x^2+7y^2$ with $x,y\in\Z$ we have
$$\sum_{k=0}^{p-1}k\bi{2k}k^3\eq\cases 11y^2/3-x^2\ (\mo\ p)&\t{if}\ 3\mid y,
\\4(y^2-x^2)/3\ (\mo\ p)&\t{if}\ 3\nmid y.\endcases$$
Since $x^2\eq-7y^2\ (\mo\ p)$, we can simplify the congruence as follows:
$$\sum_{k=0}^{p-1}k\bi{2k}k^3\eq-\f{32}{21}x^2\eq\f{32}{3}y^2\ (\mo\ p).$$

\proclaim{Conjecture A2 {\rm ([S11e])}}
{\rm (i)} If $p$ is a prime and $a$ is a positive integer with $p^a\eq1\ (\mo\ 3)$,
then
$$\sum_{k=0}^{\lfloor\f23p^a\rfloor}(21k+8)\bi{2k}k^3\eq8p^a\ (\mo\ p^{a+5+(-1)^p}).$$

{\rm (ii)} Any integer $n>1$ satisfying the congruence
$\sum_{k=0}^{n-1}(21k+8)\bi{2k}k^3\eq 8n\ (\mo\ n^4)$ must be a prime.
\endproclaim

\Remark. (a) The author [S11e] proved that for any prime $p$ and positive integer $a$
we have
$$\f1{p^a}\sum_{k=0}^{p^a-1}(21k+8)\binom{2k}k^3\eq 8+16p^3B_{p-3}\ (\mo\ p^4),$$
where $B_{-1}$ is regarded as zero, and $B_0,B_1,B_2,\ldots$ are Bernoulli numbers.
(When $p>5$, the congruence even holds mod $p^5$ if we replace $16p^3B_{p-3}$ by $-48pH_{p-1}$.)
See Conj. A1 for our guess on $\sum_{k=0}^{p-1}\bi{2k}k^3$ mod $p^2$.

(b) The author [S-7] proved that for any prime $p>3$ we have
$$\sum_{k=0}^{(p-1)/2}(21k+8)\bi{2k}k^3\eq 8p+(-1)^{(p-1)/2}32p^3E_{p-3}\ (\mo\ p^4),$$
which has the following equivalent form:
$$\sum_{k=1}^{(p-1)/2}\f{21k-8}{k^3\bi{2k}k^3}\eq(-1)^{(p+1)/2}4E_{p-3}\ (\mo\ p)\ \ \ (\t{for any prime}\ p>3).$$
Note that $\sum_{k=1}^\infty(21k-8)/(k^3\bi{2k}k^3)=\zeta(2)=\pi^2/6$ by [Z] (see also [HP1,\,(7)]).

(c) The author has verified part (ii) of Conj. A2 for $n\ls 10^4$.

\proclaim{Conjecture A3 {\rm ([S11e])}} {\rm (i)} Set
$$a_n:=\f1{4n(2n+1)\bi{2n}n}\sum_{k=0}^{n-1}(35k+8)\bi{4k}{k,k,k,k} 81^{n-1-k}\ \ \t{for}\ n\in\Z^+.$$
Then $a_n\in\Z$ unless $2n+1$ is a power of $3$ in which case $3a_n\in\Z\sm3\Z$.

{\rm (ii)} Let $p$ be a prime. If $p>3$, then
$$\f1{p^a}\sum_{k=0}^{p^a-1}\f{35k+8}{81^k}\bi{4k}{k,k,k,k}\eq 8+\f{416}{27}p^3B_{p-3}\ (\mo\ p^4)\quad\t{for all}\ a\in\Z^+,$$
and
$$\sum_{k=0}^{(p-1)/2}\f{35k+8}{81^k}\bi{4k}{k,k,k,k}\eq 8p\times 3^{p-1}\ (\mo\ p^3).$$
If $(\f p7)=1$, i.e.,
$p\eq1,2,4\ (\mo\ 7)$, then
$$\sum_{k=0}^{p-1}\bi{2k}k^3\eq\sum_{k=0}^{p-1}\f{\bi{4k}{k,k,k,k}}{81^k}\
(\mo\ p^3)$$
and
$$3\sum_{k=0}^{p-1}k\bi{2k}k^3\eq 5\sum_{k=0}^{p-1}\f{k\bi{4k}{k,k,k,k}}{81^k}\ (\mo\ p^3).$$
If $(\f p7)=-1$ and $p\not=3$, then
$$\sum_{k=0}^{p-1}\f{\bi{4k}{k,k,k,k}}{81^k}\eq0\ (\mo\ p^2).$$

{\rm (iii)} We have
$$\sum_{k=1}^\infty\f{(35k-8)81^k}{k^3\bi{4k}{k,k,k,k}}=12\pi^2.$$
\endproclaim

\proclaim{Conjecture A4 {\rm ([S11b, S11e])}} {\rm (i)} For every $n=2,3,4,\ldots$
we have
 $$\f1{n(2n+1)\bi{2n}n}\sum_{k=0}^{n-1}(11k+3)\bi{2k}k^2\bi{3k}k64^{n-1-k}\in\Z.$$

{\rm (ii)} Let $p$ be an odd prime.
 Then
 $$\align&\sum_{k=0}^{p-1}\f{\bi{2k}k^2\bi{3k}k}{64^k}
 \\\eq&\cases
 x^2-2p\ (\mo\ p^2)&\t{if}\ (\f p{11})=1\ \&\ 4p=x^2+11y^2\ (x,y\in\Z),
 \\0\ (\mo\ p^2)&\t{if}\ (\f p{11})=-1,\ \t{i.e.},\ p\eq 2,6,7,8,10\ (\mo\ 11).
 \endcases\endalign$$
Furthermore,
$$\f1{p^a}\sum_{k=0}^{p^a-1}\f{11k+3}{64^k}\bi{2k}k^2\bi{3k}k\eq3+\f 72p^3B_{p-3}\ (\mo\ p^4)\ \ \t{for all}\ a\in\Z^+.$$

{\rm (iii)} If $p>3$ is a prime then
$$p\sum_{k=1}^{(p-1)/2}\f{(11k-3)64^k}{k^3\bi{2k}k^2\bi{3k}k}\eq 32q_p(2)-\f{64}3 p^2B_{p-3}\ (\mo\ p^3),$$
where $q_p(2)=(2^{p-1}-1)/p$.
 \endproclaim
 \Remark. It is well-known that the quadratic field $\Q(\sqrt{-11})$ has class number one
 and hence for any odd prime $p$ with $(\f p{11})=1$ we can write $4p=x^2+11y^2$ with $x,y\in\Z$.
 Concerning the parameters in the representation $4p=x^2+11y^2$, Jacobi (see, e.g., [HW]) proved the following result:
 If $p=11f+1$ is a prime and $4p=x^2+11y^2$ with $x\eq2\ (\mo\ 11)$, then
 $x\eq\bi{6f}{3f}\bi{3f}f/\bi{4f}{2f}\ (\mo\ p)$.

 \proclaim{Conjecture A5 {\rm ([S11e])}} {\rm (i)} For $n\in\Z^+$ set
 $$a_n:=\f1{n(2n+1)\bi{2n}n}\sum_{k=0}^{n-1}(10k+3)\bi{2k}k^2\bi{3k}k8^{n-1-k}.$$
 Then $a_n\in\Z$ for all $n=2,3,4,\ldots$.

 {\rm (ii)} Let $p$ be an odd prime.
 Then
 $$\align&\sum_{k=0}^{p-1}\f{\bi{2k}k^2\bi{3k}k}{8^k}
 \\\eq&\cases
 4x^2-2p\ (\mo\ p^2)&\t{if}\ (\f {-2}p)=1\ \&\ p=x^2+2y^2\ (x,y\in\Z),
 \\0\ (\mo\ p^2)&\t{if}\ (\f {-2}p)=-1.
 \endcases\endalign$$
 Also, for any $a\in\Z^+$ we have
$$\f1{p^a}\sum_{k=0}^{p^a-1}\f{10k+3}{8^k}\bi{2k}k^2\bi{3k}k\eq3+\f{49}8p^3B_{p-3}\ (\mo\ p^4).$$

{\rm (iii)} We have
$$\sum_{k=1}^\infty\f{(10k-3)8^k}{k^3\bi{2k}k^2\bi{3k}k}=\f{\pi^2}2.$$
 \endproclaim

\proclaim{Conjecture A6 {\rm ([S11e])}} Let $p>3$ be a prime.
If $p\eq 1\ (\mo\ 6)$ and $p=x^2+3y^2$
with $x,y\in\Z$, then
$$\sum_{k=0}^{p-1}\f{\bi{2k}k^3}{16^k}\eq 4x^2-2p\ (\mo\ p^2)\ \t{and}\
\sum_{k=0}^{p-1}\f{k\bi{2k}k^3}{16^k}\eq p-\f {4x^2}3\ (\mo\ p^2).$$
If $p\eq 5\ (\mo\ 6)$, then
$$\sum_{k=0}^{p-1}\f{\bi{2k}k^3}{16^k}\eq0\ (\mo\ p^2)
\ \t{and}\ \sum_{k=0}^{p-1}\f{k\bi{2k}k^3}{16^k}\eq\f p3\ (\mo\ p^2).$$
Furthermore,
$$\f1{p^a}\sum_{k=0}^{p^a-1}\f{3k+1}{16^k}\bi{2k}k^3\eq 1+\f 76p^3B_{p-3}\ (\mo\ p^4)\ \ \t{for all}\ a\in\Z^+,$$
and
$$\sum_{k=0}^{(p-1)/2}(3k+1)\f{\bi{2k}k^3}{16^k}\eq p+2\l(\f{-1}p\r)p^3E_{p-3}\ (\mo\ p^4).$$
Also,
$$\f1{2n\bi{2n}n}\sum_{k=0}^{n-1}(3k+1)\bi{2k}k^316^{n-1-k}\in\Z\ \ \t{for all}\ n=2,3,4,\ldots.$$
\endproclaim
\Remark. Z. H. Sun [Su2] confirmed the conjecture that
$\sum_{k=0}^{p-1}\bi{2k}k^3/16^k\eq0\ (\mo\ p^2)$ if $p$ is a prime with $p\eq 5\ (\mo\ 6)$.
Also, the first identity of J. Guillera [G3] in the case $a=1/2$ gives
$\sum_{k=1}^\infty(3k-1)16^k/(k^3\bi{2k}k^3)=\pi^2/2$.

\proclaim{Conjecture A7 {\rm ([S11e])}} Let $p$ be an odd prime. Then
$$\sum_{k=0}^{p-1}(3k+1)\f{\bi{2k}k^3}{(-8)^k}\eq p\l(\f{-1}p\r)+p^3E_{p-3}\ (\mo\ p^4)$$
and furthermore
$$a_n:=\f1{2n\bi{2n}n}\sum_{k=0}^{n-1}(3k+1)\bi{2k}k^3(-8)^{n-1-k}\in\Z^+\ \ \ \t{for all}\ n=2,3,4,\ldots.$$
\endproclaim
\Remark. Note that $a_1=1/4$ and $(2n+1)a_{n+1}+4na_n=(3n+1)\bi{2n-1}n^2$ for $n=1,2,3,\ldots$.
Also, the third identity of Guillera [G3] with $a=1/2$ gives
$\sum_{k=1}^\infty(3k-1)(-8)^k/(k^3\bi{2k}k^3)=-2G,$
where $G$ is the Catalan constant defined by $G=\sum_{k=0}^\infty(-1)^k/(2k+1)^2=0.915965594\ldots$.

\proclaim{Conjecture A8 {\rm ([S11e])}} {\rm (i)} For $n\in\Z^+$ set
$$a_n:=\f1{n(2n+1)\bi{2n}n}\sum_{k=0}^{n-1}(5k+1)\bi{2k}{k}^2\bi{3k}{k}(-192)^{n-1-k}.$$
Then $a_n\in\Z$ for $n=2,3,4,\ldots$ unless $2n+1$ is a power of $3$ in which case $3a_n\in\Z\sm 3\Z$.

{\rm (ii)} Let $p>3$ be a prime.
Then
$$\f1{p^a}\sum_{k=0}^{p^a-1}\f{5k+1}{(-192)^k}\bi{2k}{k}^2\bi{3k}{k}
\eq \l(\f{p^a}3\r)+\l(\f{p^{a-1}}3\r)\f5{18}p^2B_{p-2}\l(\f13\r)\ (\mo\ p^3)$$
for any $a\in\Z^+$.
We also have
 $$\align&\sum_{k=0}^{p-1}\f{\bi{2k}{k}^2\bi{3k}{k}}{(-192)^k}
 \\\eq&\cases
 x^2-2p\ (\mo\ p^2)&\t{if}\ p\eq1\ (\mo\ 3)\ \&\ 4p=x^2+27y^2\ (x,y\in\Z),
 \\0\ (\mo\ p^2)&\t{if}\  p\eq2\ (\mo\ 3).
 \endcases\endalign$$
 \endproclaim
\Remark. It is well known that for any prime $p\eq1\ (\mo\ 3)$ there are unique $x,y\in\Z^+$ such that $4p=x^2+27y^2$
 (see, e.g., [C]). Also,  Ramanujan [R] found  that
$$\sum_{k=0}^\infty(5k+1)\l(-\f 9{16}\r)^k\f{(1/2)_k(1/3)_k(2/3)_k}{(1)_k^3}
=\sum_{k=0}^\infty\f{5k+1}{(-192)^{k}}\bi{2k}{k}^2\bi{3k}{k}=\f{4\sqrt3}{\pi}.$$

\proclaim{Conjecture A9 {\rm ([S-12])}}  Let $p>3$ be a prime. Then
$$\align&\sum_{k=0}^{p-1}\f{\bi{6k}{3k}\bi{3k}{k,k,k}}{(-96)^{3k}}
\\\eq&\cases(\f{-6}p)(x^2-2p)\ (\mo\ p^2)&\t{if}\ (\f p{19})=1\ \&\ 4p=x^2+19y^2\ (x,y\in\Z),
\\0\ (\mo\ p^2)&\t{if}\ (\f p{19})=-1.\endcases
\endalign$$
Also,
$$\sum_{k=0}^{p-1}\f{342k+25}{(-96)^{3k}}\bi{6k}{3k}\bi{3k}{k,k,k}\eq 25p\l(\f{-6}p\r)\ (\mo\ p^3).$$
Furthermore, for any $n=2,3,\ldots$ we have
$$a_n:=\f1{2n(2n+1)\bi{2n}n}\sum_{k=0}^{n-1}(342k+25)(-96^3)^{n-1-k}\binom{6k}{3k}\binom{3k}{k,k,k}\in\Z$$
unless $2n+1$ is a power of $3$ in which case $3a_n\in\Z\sm3\Z$.
\endproclaim
\Remark. It is well known that $\Q(\sqrt{-19})$ has class number one and hence for any odd prime $p$ with $(\f p{19})=1$
there are unique positive integers $x$ and $y$ such that $4p=x^2+19y^2$.
The conjectured congruences modulo $p$ have been confirmed by Z. H. Sun [Su3].
D. V. Chudnovsky and G. V. Chunovsky [CC]
obtained that
$$\sum_{k=0}^\infty\f{342k+25}{(-96)^{3k}}\bi{6k}{3k}\bi{3k}{k,k,k}=\f{32\sqrt6}{\pi}.$$
See also [BB] and [G5] for more Ramanujan-type series involving $1/\pi$.

\proclaim{Conjecture A10 {\rm ([S-12])}} If $p>5$ is a prime, then
$$\align&\sum_{k=0}^{p-1}\f{\bi{6k}{3k}\bi{3k}{k,k,k}}{(-960)^{3k}}
\\\eq&\cases(\f{p}{15})(x^2-2p)\ (\mo\ p^2)&\t{if}\ (\f p{43})=1\ \&\ 4p=x^2+43y^2\ (x,y\in\Z),
\\0\ (\mo\ p^2)&\t{if}\ (\f p{43})=-1.\endcases
\endalign$$
Also, for any $n=2,3,\ldots$ we have
$$a_n:=\f1{2n(2n+1)\bi{2n}n}\sum_{k=0}^{n-1}(5418k+263)(-960^3)^{n-1-k}\binom{6k}{3k}\binom{3k}{k,k,k}\in\Z$$
unless $2n+1$ is a power of $3$ in which case $3a_n\in\Z\sm3\Z$.
\endproclaim
\Remark. It is well known that $\Q(\sqrt{-43})$ has class number one and hence for any odd prime $p$ with $(\f p{43})=1$
there are unique positive integers $x$ and $y$ such that $4p=x^2+43y^2$.
The conjectured congruence modulo $p$ has been confirmed by Z. H. Sun [Su3].
D. V. Chudnovsky and G. V. Chunovsky [CC] showed that
$$\sum_{k=0}^\infty\f{5418k+263}{(-960)^{3k}}\bi{6k}{3k}\bi{3k}{k,k,k}=\f{640\sqrt{15}}{3\pi},$$
and Zudilin [Zu] suggested that for any prime $p>5$ we should have
$$\sum_{k=0}^{p-1}\f{5418k+263}{(-960)^{3k}}\bi{6k}{3k}\bi{3k}{k,k,k}\eq 263p\l(\f{-15}p\r)\ (\mo\ p^3).$$

\proclaim{Conjecture A11 {\rm ([S-12])}} Let $p>5$ be a prime with $p\not=11$. Then
$$\align&\sum_{k=0}^{p-1}\f{\bi{6k}{3k}\bi{3k}{k,k,k}}{(-5280)^{3k}}
\\\eq&\cases(\f{-330}{p})(x^2-2p)\ (\mo\ p^2)&\t{if}\ (\f p{67})=1\ \&\ 4p=x^2+67y^2\ (x,y\in\Z),
\\0\ (\mo\ p^2)&\t{if}\ (\f p{67})=-1.\endcases
\endalign$$
Also,
$$\sum_{k=0}^{p-1}\f{261702k+10177}{(-5280)^{3k}}\bi{6k}{3k}\bi{3k}{k,k,k}\eq 10177p\l(\f{-330}p\r)\ (\mo\ p^3).$$
Furthermore, for any $n=2,3,\ldots$ we have
$$a_n:=\f1{2n(2n+1)\bi{2n}n}\sum_{k=0}^{n-1}(261702k+10177)(-5280^3)^{n-1-k}\binom{6k}{3k}\binom{3k}{k,k,k}\in\Z$$
unless $2n+1$ is a power of $3$ in which case $3a_n\in\Z\sm3\Z$.
\endproclaim
\Remark. It is well known that $\Q(\sqrt{-67})$ has class number one and hence for any odd prime $p$ with $(\f p{67})=1$
there are unique positive integers $x$ and $y$ such that $4p=x^2+67y^2$.
The conjectured congruences modulo $p$ have been confirmed by Z. H. Sun [Su3].
It is known that (cf. [CC] and [G5])
$$\align&\sum_{k=0}^\infty\f{(261702k+10177)(-1)^k(1/2)_k(1/6)_k(5/6)_k}{440^{3k}}
\\&\quad=\sum_{k=0}^\infty\f{261702k+10177}{(-5280)^{3k}}\bi{6k}{3k}\bi{3k}{k,k,k}=\f{3\times440^2}{\pi\sqrt{330}}.
\endalign$$

\proclaim{Conjecture A12 {\rm ([S-12])}} Let $p>5$ be a prime with $p\not=23,29$. Then
$$\align&\sum_{k=0}^{p-1}\f{\bi{6k}{3k}\bi{3k}{k,k,k}}{(-640320)^{3k}}
\\\eq&\cases(\f{-10005}{p})(x^2-2p)\ (\mo\ p^2)&\t{if}\ (\f p{163})=1\ \&\ 4p=x^2+163y^2\ (x,y\in\Z),
\\0\ (\mo\ p^2)&\t{if}\ (\f p{163})=-1.\endcases
\endalign$$
Also,
$$\sum_{k=0}^{p-1}\f{545140134k+13591409}{(-640320)^{3k}}\bi{6k}{3k}\bi{3k}{k,k,k}\eq 13591409p\l(\f{-10005}p\r)\ (\mo\ p^3).$$
Furthermore, for $n=2,3,\ldots$, if we denote by $a_n$ the number
$$\f1{2n(2n+1)\bi{2n}n}\sum_{k=0}^{n-1}(545140134k+13591409)(-640320^3)^{n-1-k}\binom{6k}{3k}\binom{3k}{k,k,k},$$
then $a_n\in\Z$ unless $2n+1$ is a power of $3$ in which case $3a_n\in\Z\sm3\Z$.
\endproclaim
\Remark. It is well known that the only imaginary quadratic fields with class number one are those
$\Q(\sqrt{-d})$ with $d=1,2,3,7,11,19,43,67,163$.
For any odd prime $p$ with $(\f p{163})=1$,
there are unique positive integers $x$ and $y$ such that $4p=x^2+163y^2$.
The conjectured congruences modulo $p$ have been confirmed by Z. H. Sun [Su3].
D. V. Chudnovsky and G. V. Chudnovsky [CC]
got the formula
$$\sum_{k=0}^\infty\f{545140134k+13591409}{(-640320)^{3k}}\bi{6k}{3k}\bi{3k}{k,k,k}=\f{3\times53360^2}{2\pi\sqrt{10005}},$$
which enabled them to hold the record for the calculation of $\pi$ during 1989-1994.

 \proclaim{Conjecture A13 {\rm ([S11e])}} {\rm (i)} If $p>3$ is a prime, then
 $$\align&\sum_{k=0}^{p-1}\f{\bi{2k}k^2\bi{3k}k}{(-27)^k}
 \\\eq&\cases
 4x^2-2p\ (\mo\ p^2)&\t{if}\ p\eq1,4\ (\mo\ 15)\ \&\ p=x^2+15y^2\ (x,y\in\Z),
 \\20x^2-2p\ (\mo\ p^2)&\t{if}\ p\eq2,8\ (\mo\ 15)\ \&\ p=5x^2+3y^2\ (x,y\in\Z),
 \\0\ (\mo\ p^2)&\t{if}\ (\f p{15})=-1.
 \endcases\endalign$$
 For $n\in\Z^+$ set
 $$\f1{2n(2n+1)\bi{2n}n}\sum_{k=0}^{n-1}(15k+4)\bi{2k}k^2\bi{3k}k(-27)^{n-1-k}.$$
Then $a_n\in\Z$ unless $2n+1$ is a power of $3$ in which case $3a_n\in\Z\sm 3\Z$.

 {\rm (ii)} For any prime $p>3$ and $a\in\Z^+$, we have
$$\f1{p^a}\sum_{k=0}^{p^a-1}\f{15k+4}{(-27)^k}\bi{2k}k^2\bi{3k}k\eq 4\l(\f {p^a}3\r)+\l(\f{p^{a-1}}3\r)\f43p^2B_{p-2}\l(\f13\r)\ (\mo\ p^3)$$
and
$$\f1{p^a}\sum_{k=0}^{p^a-1}\f{5k+1}{(-144)^k}\bi{4k}{k,k,k,k}
\eq\l(\f{p^a}3\r)+\l(\f{p^{a-1}}3\r)\f 5{12}p^2B_{p-2}\l(\f13\r)\ (\mo\ p^3).$$
We also have
$$ \sum_{k=1}^\infty\f{(5k-1)(-144)^k}{k^3\bi{2k}k^2\bi{4k}{2k}}=-\f{45}2K,$$
where
$$K:=L\l(2,\l(\f{-3}{\cdot}\r)\r)=\sum_{k=1}^\infty\f{(\f k3)}{k^2}
=0.781302412896486296867187429624\ldots.$$
 \endproclaim
\Remark. (a) Let $p>5$ be a prime. By the theory of binary quadratic forms (cf. [C]),
if $p\eq1,4\ (\mo\ 15)$ then $p=x^2+15y^2$ for some $x,y\in\Z$; if $p\eq 2,8\ (\mo\ 15)$ then
$p=5x^2+3y^2$ for some $x,y\in\Z$.

(b) Concerning the first congruence in Conj. A13(ii),
K. Hessami Pilehrood and T. Hessami Pilehrood [HP2] proved it modulo $p$ for $a=1$.

\proclaim{Conjecture A14 {\rm ([S11e])}} {\rm (i)} For  $n\in\Z^+$ set
$$a_n:=\f1{n(2n+1)\bi{2n}n}\sum_{k=0}^{n-1}(6k+1)\bi{2k}k^2\bi{3k}k6^{3(n-1-k)}$$
and
$$b_n:=\f1{2n(2n+1)\bi{2n}n}\sum_{k=0}^{n-1}(8k+1)\bi{4k}{k,k,k,k}48^{2(n-1-k)}.$$
Given an integer $n>1$, we have $a_n,b_n\in\Z$ unless $2n+1$ is a power of $3$ in which case
$3a_n,3b_n\in\Z\sm3\Z$.

{\rm (ii)} Let $p>3$ be a prime.
 Then
 $$\align&\sum_{k=0}^{p-1}\f{\bi{2k}k^2\bi{3k}k}{6^{3k}}
 \\\eq&\cases
 4x^2-2p\ (\mo\ p^2)&\t{if}\ p\eq1,7\ (\mo\ 24)\ \&\ p=x^2+6y^2\ (x,y\in\Z),
 \\8x^2-2p\ (\mo\ p^2)&\t{if}\ p\eq5,11\ (\mo\ 24)\ \&\ p=2x^2+3y^2\ (x,y\in\Z),
 \\0\ (\mo\ p^2)&\t{if}\ (\f {-6}p)=-1\ i.e.,\ p\eq 13,17,19,23\ (\mo\ 24);
 \endcases\endalign$$
 and
$$\align&\sum_{k=0}^{p-1}\f{\bi{4k}{k,k,k,k}}{48^{2k}}
 \\\eq&\cases
 4x^2-2p\ (\mo\ p^2)&\t{if}\ p\eq1,7\ (\mo\ 24)\ \&\ p=x^2+6y^2\ (x,y\in\Z),
 \\2p-8x^2\ (\mo\ p^2)&\t{if}\ p\eq5,11\ (\mo\ 24)\ \&\ p=2x^2+3y^2\ (x,y\in\Z),
 \\0\ (\mo\ p^2)&\t{if}\ (\f {-6}p)=-1\ i.e.,\ p\eq 13,17,19,23\ (\mo\ 24).
 \endcases\endalign$$
Also, for any $a\in\Z^+$ we have
 $$\align\f1{p^a}\sum_{k=0}^{p^a-1}\f{6k+1}{6^{3k}}\bi{2k}k^2\bi{3k}k
 \eq& \l(\f {p^a}3\r)-\l(\f{p^{a-1}}3\r)\f5{12}p^2B_{p-2}\l(\f13\r)\ (\mo\ p^3),
\\\f1{p^a}\sum_{k=0}^{p^a-1}\f{8k+1}{48^{2k}}\bi{4k}{k,k,k,k}\eq& \l(\f {p^a}3\r)-\l(\f{p^{a-1}}3\r)\f5{24}p^2B_{p-2}\l(\f13\r)\ (\mo\ p^3),
\endalign$$
and
$$\f 1{p^a}\sum_{p^a/2<k<p^a}\f{8k+1}{48^{2k}}\bi{4k}{k,k,k,k}\eq0\pmod{p^2}.$$
 \endproclaim
\Remark. (a) Let $p>3$ be a prime. By the theory of binary quadratic forms (see, e.g., [C]),
if $p\eq1,7\ (\mo\ 24)$ then $p=x^2+6y^2$ for some $x,y\in\Z$; if $p\eq 5,11\ (\mo\ 24)$ then
$p=2x^2+3y^2$ for some $x,y\in\Z$.

 (b) Ramanujan [R] found that
 $$\sum_{k=0}^\infty(6k+1)\f{(1/2)_k(1/3)_k(2/3)_k}{2^k(1)_k^3}=\sum_{k=0}^\infty\f{6k+1}{6^{3k}}\bi{2k}{k}^2\bi{3k}k=\f{3\sqrt3}{\pi}$$
and
$$\sum_{k=0}^\infty(8k+1)\f{(1/2)_k(1/4)_k(3/4)_k}{9^k(1)_k^3}=\sum_{k=0}^\infty\f{8k+1}{48^{2k}}\bi{4k}{k,k,k,k}=\f{2\sqrt3}{\pi}.$$

\proclaim{Conjecture A15 {\rm ([S11e])}} {\rm (i)} For $n\in\Z^+$ set
$$a_n:=\f1{2n(2n+1)\bi{2n}n}\sum_{k=0}^{n-1}(20k+3)\bi{4k}{k,k,k,k}(-2^{10})^{n-1-k}.$$
Then $(-1)^{n-1}a_n\in\Z^+$ for all $n=2,3,4,\ldots$.

{\rm (ii)} Let $p$ be an odd prime. Then
 $$\align&\sum_{k=0}^{p-1}\f{\bi{4k}{k,k,k,k}}{(-2^{10})^k}
 \\\eq&\cases
 4x^2-2p\ (\mo\ p^2)&\t{if}\ p\eq1,9\ (\mo\ 20)\ \&\ p=x^2+5y^2\ (x,y\in\Z),
 \\2(p-x^2)\ (\mo\ p^2)&\t{if}\ p\eq3,7\ (\mo\ 20)\ \&\ 2p=x^2+5y^2\ (x,y\in\Z),
 \\0\ (\mo\ p^2)&\t{if}\ (\f {-5}p)=-1,\ i.e.,\ p\eq 11,13,17,19\ (\mo\ 20).\endcases
 \endalign$$
 \endproclaim

\Remark. Let $p\not=2,5$ be a prime. By the theory of binary quadratic
forms (see, e.g., [C]), if $p\eq1,9\ (\mo\ 20)$ then $p=x^2+5y^2$
for some $x,y\in\Z$; if $p\eq 3,7\ (\mo\ 20)$ then $2p=x^2+5y^2$ for
some $x,y\in\Z$. See also [S-2] for a $p$-adic congruence mod $p^4$ which is an
analogue of the Ramanujan series
$$\sum_{k=0}^\infty\f{20k+3}{(-2^{10})^k}\bi{4k}{k,k,k,k}=\f{8}{\pi}.$$

\proclaim{Conjecture A16 {\rm ([S11e])}} {\rm (i)} For $n\in\Z^+$ set
$$a_n:=\f1{2n(2n+1)\bi{2n}n}\sum_{k=0}^{n-1}(10k+1)\bi{4k}{k,k,k,k}12^{4(n-1-k)}.$$
Given an integer $n>1$, we have $a_n\in\Z$ unless $2n+1$ is a power of $3$ in which case
$3a_n\in\Z\sm3\Z$.

{\rm (ii)} Let $p>3$ be a prime.
Then
$$\f1{p^a}\sum_{k=0}^{p^a-1}\f{10k+1}{12^{4k}}\bi{4k}{k,k,k,k}
\eq \l(\f{-2}{p^a}\r)-\l(\f{-2}{p^{a-1}}\r)\f{p^2}{48}E_{p-3}\l(\f14\r)\ (\mo\ p^3)$$
for all $a=1,2,3,\ldots$.
We also have
 $$\align&\sum_{k=0}^{p-1}\f{\bi{4k}{k,k,k,k}}{12^{4k}}
 \\\eq&\cases
 4x^2-2p\ (\mo\ p^2)&\t{if}\ p\eq1,9,11,19\ (\mo\ 40)\ \&\ p=x^2+10y^2\ (x,y\in\Z),
 \\2p-8x^2\ (\mo\ p^2)&\t{if}\ p\eq7,13,23,37\ (\mo\ 40)\ \&\ p=2x^2+5y^2\ (x,y\in\Z),
 \\0\ (\mo\ p^2)&\t{if}\ (\f {-10}p)=-1,\ i.e.,\ p\eq 3,17,21,27,29,31,33,39\ (\mo\ 40).
 \endcases\endalign$$
 \endproclaim
\Remark. (a) Let $p>5$ be a prime. By the theory of binary quadratic
forms (see, e.g., [C]), if $(\f{-2}p)=(\f p5)=1$ then $p=x^2+10y^2$
for some $x,y\in\Z$; if $(\f{-2}p)=(\f p5)=-1$ then $p=2x^2+5y^2$ for
some $x,y\in\Z$.

(b) Ramanujan [R] obtained that
$$\sum_{k=0}^\infty(10k+1)\f{(1/2)_k(1/4)_k(3/4)_k}{81^k(1)_k^3}
=\sum_{k=0}^\infty\f{10k+1}{12^{4k}}\bi{4k}{k,k,k,k}=\f{9\sqrt2}{4\pi}.$$

\proclaim{Conjecture A17 {\rm ([S-12])}} Let $p>3$ be a prime. Then
$$\align&\sum_{k=0}^{p-1}\f{\bi{4k}{k,k,k,k}}{(-2^{10} 3^4)^k}
\\\eq&\cases 4x^2-2p\ (\mo\ p^2)&\t{if}\ (\f{13}p)=(\f{-1}p)=1\ \&\ p=x^2+13y^2,
\\2p-2x^2\ (\mo\ p^2)&\t{if}\ (\f{13}p)=(\f{-1}p)=-1\ \&\ 2p=x^2+13y^2,
\\0\ (\mo\ p^2)&\t{if}\ (\f{13}p)=-(\f{-1}p).\endcases
\endalign$$
We also have
$$\sum_{k=0}^{p-1}\f{260k+23}{(-82944)^k}\bi{4k}{k,k,k,k}\eq 23p\l(\f{-1}p\r)+\f 53 p^3E_{p-3}\ (\mo\ p^4).$$
Furthermore, for $n=2,3,4,\ldots$ we have
$$a_n:=\f1{2n(2n+1)\bi{2n}n}\sum_{k=0}^{n-1}(260k+23)\bi{4k}{k,k,k,k}(-82944)^{n-1-k}\in\Z$$
unless $2n+1$ is a power of $3$ in which case $3a_n\in\Z\sm3\Z$.
\endproclaim
\Remark. Ramanujan (cf. [Be, p.\,353]) found that
$$\sum_{k=0}^\infty\f{(260k+23)(1/2)_k(1/4)_k(3/4)_k}{k!^318^{2k}}
=\sum_{k=0}^\infty\f{260k+23}{(-82944)^k}\bi{4k}{k,k,k,k}=\f{72}{\pi}.$$

\proclaim{Conjecture A18 {\rm ([S-12])}} Let $p>3$ be a prime with $p\not=11$. Then
$$\align&\sum_{k=0}^{p-1}\f{\bi{4k}{k,k,k,k}}{1584^{2k}}
\\\eq&\cases 4x^2-2p\ (\mo\ p^2)&\t{if}\ (\f{-11}p)=(\f{2}p)=1\ \&\ p=x^2+22y^2,
\\2p-8x^2\ (\mo\ p^2)&\t{if}\ (\f{-11}p)=(\f{2}p)=-1\ \&\ p=2x^2+11y^2,
\\0\ (\mo\ p^2)&\t{if}\ (\f{-11}p)=-(\f{2}p).\endcases
\endalign$$
We also have
$$\sum_{k=0}^{p-1}\f{280k+19}{1584^{2k}}\bi{4k}{k,k,k,k}\eq 19p\l(\f{p}{11}\r)\ (\mo\ p^3).$$
Furthermore, for $n=2,3,4,\ldots$ we have
$$a_n:=\f1{2n(2n+1)\bi{2n}n}\sum_{k=0}^{n-1}(280k+19)\bi{4k}{k,k,k,k}1584^{2(n-1-k)}\in\Z$$
unless $2n+1$ is a power of $3$ in which case $3a_n\in\Z\sm3\Z$.
\endproclaim
\Remark. Ramanujan (cf. [Be, p.\,354]) found that
$$\sum_{k=0}^\infty\f{(280k+19)(1/2)_k(1/4)_k(3/4)_k}{k!^399^{2k}}
=\sum_{k=0}^\infty\f{280k+19}{1584^{2k}}\bi{4k}{k,k,k,k}=\f{2\times99^2}{\pi\sqrt{11}}.$$

\proclaim{Conjecture A19 {\rm ([S-12])}} Let $p>3$ be a prime with $p\not=7$. Then
$$\align&\sum_{k=0}^{p-1}\f{\bi{4k}{k,k,k,k}}{(-2^{10} 21^4)^k}
\\\eq&\cases 4x^2-2p\ (\mo\ p^2)&\t{if}\ (\f{37}p)=(\f{-1}p)=1\ \&\ p=x^2+37y^2,
\\2p-2x^2\ (\mo\ p^2)&\t{if}\ (\f{37}p)=(\f{-1}p)=-1\ \&\ 2p=x^2+37y^2,
\\0\ (\mo\ p^2)&\t{if}\ (\f{-37}p)=-1.\endcases
\endalign$$
Furthermore, for $n=2,3,4,\ldots$ we have
$$a_n:=\f1{2n(2n+1)\bi{2n}n}\sum_{k=0}^{n-1}(21460k+1123)\bi{4k}{k,k,k,k}(-2^{10}21^4)^{n-1-k}\in\Z$$
unless $2n+1$ is a power of $3$ in which case $3a_n\in\Z\sm3\Z$.
\endproclaim
\Remark. Ramanujan (cf. [Be, p.\,353]) found that
$$\align&\sum_{k=0}^\infty\f{(21460k+1123)(-1)^k(1/2)_k(1/4)_k(3/4)_k}{k!^3 882^{2k}}
\\&\ \ =\sum_{k=0}^\infty\f{21460k+1123}{(-2^{10}21^4)^k}\bi{4k}{k,k,k,k}=\f{2^321^2}{\pi}.
\endalign$$

\proclaim{Conjecture A20 {\rm ([S-12])}} For $n\in\Z^+$ set
$$a_n:=\f1{n(2n+1)\bi{2n}n}\sum_{k=0}^{n-1}(51k+7)\bi{2k}{k}^2\bi{3k}{k}(-12^3)^{n-1-k}.$$
Given an integer $n>1$, we have $(-1)^{n-1}a_n\in\Z^+$ unless $2n+1$ is a power of $3$
in which case $3a_n\in\Z\sm 3\Z$.

{\rm (ii)} Let $p>3$ be a prime.
Then
$$\f1{p^a}\sum_{k=0}^{p^a-1}\f{51k+7}{(-12^3)^k}\bi{2k}{k}^2\bi{3k}{k}
\eq 7\l(\f{p^a}3\r)+\l(\f{p^{a-1}}3\r)\f56p^2B_{p-2}\l(\f13\r)\ (\mo\ p^3)$$
for every $a\in\Z^+$.
We also have
 $$\align&\sum_{k=0}^{p-1}\f{\bi{2k}{k}^2\bi{3k}{k}}{(-12^3)^k}
 \\\eq&\cases
 x^2-2p\ (\mo\ p^2)&\t{if}\ (\f p3)=(\f p{17})=1\ \&\ 4p=x^2+51y^2\ (x,y\in\Z),
 \\2p-3x^2\ (\mo\ p^2)&\t{if}\ (\f p3)=(\f p{17})=-1\ \&\ 4p=3x^2+17y^2\ (x,y\in\Z),
 \\0\ (\mo\ p^2)&\t{if}\  (\f p3)=-(\f p{17}).
 \endcases\endalign$$
 \endproclaim
\Remark. (a) Let $p>3$ be a prime. By the theory of binary quadratic forms (see, e.g., [C]),
if $(\f p3)=(\f p{17})=1$ then $4p=x^2+51y^2$ for some $x,y\in\Z$; if $(\f p3)=(\f p{17})=-1$
then $4p=3x^2+17y^2$ for some $x,y\in\Z$.
(b) Ramanujan [R] obtained that
$$\sum_{k=0}^\infty(51k+7)\f{(1/2)_k(1/3)_k(2/3)_k}{(-16)^k(1)_k^3}
=\sum_{k=0}^\infty\f{51k+7}{(-12^3)^{k}}\bi{2k}{k}^2\bi{3k}{k}=\f{12\sqrt3}{\pi}.$$

\proclaim{Conjecture A21 {\rm ([S-12])}} Let $p>3$ be a prime with $p\not=11$. Then
$$\align&\sum_{k=0}^{p-1}\f{\bi{4k}{k,k,k,k}}{396^{4k}}
\\\eq&\cases 4x^2-2p\ (\mo\ p^2)&\t{if}\ (\f{29}p)=(\f{-2}p)=1\ \&\ p=x^2+58y^2,
\\2p-8x^2\ (\mo\ p^2)&\t{if}\ (\f{29}p)=(\f{-2}p)=-1\ \&\ p=2x^2+29y^2,
\\0\ (\mo\ p^2)&\t{if}\ (\f{-58}p)=-1.\endcases
\endalign$$
Furthermore, for $n=2,3,4,\ldots$ we have
$$a_n:=\f1{2n(2n+1)\bi{2n}n}\sum_{k=0}^{n-1}(26390k+1103)\bi{4k}{k,k,k,k}396^{4(n-1-k)}\in\Z$$
unless $2n+1$ is a power of $3$ in which case $3a_n\in\Z\sm3\Z$.
\endproclaim
\Remark. Ramanujan (cf. [Be, p.\,354]) found that
$$\align&\sum_{k=0}^\infty\f{(26390k+1103)(1/2)_k(1/4)_k(3/4)_k}{k!^3 99^{4k}}
\\&\ \ =\sum_{k=0}^\infty\f{26390k+1103}{396^{4k}}\bi{4k}{k,k,k,k}=\f{99^2}{2\pi\sqrt2}.
\endalign$$

\proclaim{Conjecture A22 {\rm ([S-12])}} Let $p>3$ be a prime. Then
$$\align&\sum_{k=0}^{p-1}\f{\bi{2k}k^2\bi{3k}k}{(-48)^{3k}}
\\\eq&\cases x^2-2p\ (\mo\ p^2)&\t{if}\ (\f{p}3)=(\f{p}{41})=1\ \&\ 4p=x^2+123y^2,
\\2p-3x^2\ (\mo\ p^2)&\t{if}\ (\f{p}3)=(\f{p}{41})=-1\ \&\ 4p=3x^2+41y^2,
\\0\ (\mo\ p^2)&\t{if}\ (\f{p}{123})=-1.\endcases
\endalign$$
Also,
$$\f1{p^a}\sum_{k=0}^{p^a-1}\f{615k+53}{(-48)^{3k}}\bi{2k}k^2\bi{3k}k\eq 53\l(\f {p^a}3\r)
+\l(\f{p^{a-1}}3\r)\f5{12}p^2B_{p-2}\l(\f13\r)\ (\mo\ p^3)$$
for any positive integer $a$. Furthermore, for $n=2,3,4,\ldots$ we have
$$a_n:=\f1{n(2n+1)\bi{2n}n}\sum_{k=0}^{n-1}(615k+53)\bi{2k}{k}^2\bi{3k}{k}(-48)^{3(n-1-k)}\in\Z$$
unless $2n+1$ is a power of $3$ in which case $3a_n\in\Z\sm 3\Z$.
\endproclaim
\Remark. It is known (cf. [G5]) that
$$\align&\sum_{k=0}^\infty\f{(615k+53)(-1)^k(1/2)_k(1/3)_k(2/3)_k}{k!^3 2^{10k}}
\\&\ \ =\sum_{k=0}^\infty\f{615k+53}{(-48)^{3k}}\bi{2k}{k}^2\bi{3k}{k}=\f{96\sqrt3}{\pi}.
\endalign$$

\proclaim{Conjecture A23 {\rm ([S-12])}} Let $p>5$ be a prime. Then
$$\align&\sum_{k=0}^{p-1}\f{\bi{2k}k^2\bi{3k}k}{(-300)^{3k}}
\\\eq&\cases x^2-2p\ (\mo\ p^2)&\t{if}\ (\f{p}3)=(\f{p}{89})=1\ \&\ 4p=x^2+267y^2,
\\2p-3x^2\ (\mo\ p^2)&\t{if}\ (\f{p}3)=(\f{p}{89})=-1\ \&\ 4p=3x^2+89y^2,
\\0\ (\mo\ p^2)&\t{if}\ (\f{p}{267})=-1.\endcases
\endalign$$
Also,
$$\align &\f1{p^a}\sum_{k=0}^{p^a-1}\f{14151k+827}{(-300)^{3k}}\bi{2k}k^2\bi{3k}k
\\\eq& 827\l(\f {p^a}3\r)+\l(\f{p^{a-1}}3\r)\f{13}{150}p^2B_{p-2}\l(\f13\r)\ (\mo\ p^3)
\endalign$$
for any positive integer $a$.
Furthermore, for $n=2,3,4,\ldots$ we have
$$a_n:=\f1{2n\bi{2n}n}\sum_{k=0}^{n-1}(14151k+827)\bi{2k}{k}^2\bi{4k}{2k}(-300)^{3(n-1-k)}\in\Z$$
unless $n-1$ is a power of $2$ in which case $2a_n$ is an odd integer.
\endproclaim
\Remark. It is known (cf. [G5]) that
$$\align&\sum_{k=0}^\infty\f{(14151k+827)(-1)^k(1/2)_k(1/3)_k(2/3)_k}{k!^3 500^{2k}}
\\&\ \ =\sum_{k=0}^\infty\f{14151k+827}{(-300)^{3k}}\bi{2k}{k}^2\bi{4k}{2k}=\f{1500\sqrt3}{\pi}
\endalign$$
and W. Zudilin [Zu] suggested that for any prime $p>5$ we have
$$\sum_{k=0}^{p-1}\f{14151k+827}{(-300)^{3k}}\bi{2k}k^2\bi{3k}k\eq 827p\l(\f p3\r)\ (\mo\ p^3).$$

\proclaim{Conjecture A24 {\rm ([S11e])}} {\rm (i)} For $n\in\Z^+$ set
$$a_n:=\f1{2n(2n+1)\bi{2n}n}\sum_{k=0}^{n-1}(28k+3)\bi{4k}{k,k,k,k}(-3\times2^{12})^{n-1-k}.$$
Then we have $(-1)^{n-1}a_n\in\Z^+$ for all $n=2,3,4,\ldots$.

{\rm (ii)} Let $p>3$ be a prime.
Then
$$\f1{p^a}\sum_{k=0}^{p^a-1}\f{28k+3}{(-3\times2^{12})^k}\bi{4k}{k,k,k,k}
\eq 3\l(\f {p^a}3\r)+\l(\f{p^{a-1}}3\r)\f5{24}p^2B_{p-2}\l(\f13\r)\ (\mo\ p^3)$$
for every $a=1,2,3,\ldots$.
We also have
 $$\align&\sum_{k=0}^{p-1}\f{\bi{4k}{k,k,k,k}}{(-3\times 2^{12})^k}
 \\\eq&\cases
 (-1)^{\lfloor x/6\rfloor}(4x^2-2p)\ (\mo\ p^2)&\t{if}\ 12\mid p-1\ \&\ p=x^2+y^2\ (4\mid x-1\ \&\ 2\mid y),
 \\-4(\f{xy}3)xy\ (\mo\ p^2)&\t{if}\ 12\mid p-5\ \&\ p=x^2+y^2\  (4\mid x-1\ \&\ 2\mid y),
 \\0\ (\mo\ p^2)&\t{if}\ p\eq3\ (\mo\ 4).
 \endcases\endalign$$
 \endproclaim
\Remark.  Ramanujan [R] obtained that
$$\sum_{k=0}^\infty(28k+3)\f{(1/2)_k(1/4)_k(3/4)_k}{(-48)^k(1)_k^3}
=\sum_{k=0}^\infty\f{28k+3}{(-3\times2^{12})^{k}}\bi{4k}{k,k,k,k}=\f{16\sqrt3}{3\pi}.$$

\proclaim{Conjecture A25 {\rm ([S11e])}} Let $p>5$ be a prime. Then
$$\align &\sum_{k=0}^{p-1}\f{\bi{4k}{k,k,k,k}}{(-2^{14}3^4 5)^k}
\\\eq&\cases \ve(x)(4x^2-2p)&\t{if}\ 4\mid p-1,\ (\f{p}5)=1,\ p=x^2+y^2\ \&\ 2\nmid x,
\\4xy\ (\mo\ p^2)&\t{if}\ 4\mid p-1,\ (\f p5)=-1,\ p=x^2+y^2\ \&\ 5\mid x+y,
\\0\ (\mo\ p^2)&\t{if}\ p\eq3\ (\mo\ 4),\endcases
\endalign$$
where $\ve(x)$ takes $-1$ or $1$ according as $5\mid x$ or not. We also have
$$\sum_{k=0}^{p-1}\f{644k+41}{(-2^{14}3^4 5)^k}\bi{4k}{k,k,k,k}\eq 41p\l(\f{-5}p\r)\ (\mo\ p^3).$$
Furthermore, for $n=2,3,4,\ldots$ we have
$$a_n:=\f1{2n(2n+1)\bi{2n}n}\sum_{k=0}^{n-1}(644k+41)\bi{4k}{k,k,k,k}(-2^{14}3^4 5)^{n-1-k}\in\Z$$
unless $2n+1$ is a power of $3$ in which case $3a_n\in\Z\sm3\Z$.
\endproclaim
\Remark.  Ramanujan (cf. [Be, p.\,353]) found that
$$\sum_{k=0}^\infty(644k+41)\f{(1/2)_k(1/4)_k(3/4)_k}{k!^3(-5)^k72^{2k}}
=\sum_{k=0}^\infty\f{644k+41}{(-2^{14}3^4 5)^k}\bi{4k}{k,k,k,k}=\f{288}{\pi\sqrt5}.$$

\proclaim{Conjecture A26 {\rm ([S11e])}} {\rm (i)} For  $n\in\Z^+$ set
$$a_n:=\f1{10n(2n+1)\bi{2n}n}\sum_{k=0}^{n-1}(154k+15)\bi{6k}{3k}\bi{3k}{k,k,k}(-2^{15})^{n-1-k}.$$
Given an integer $n>1$, we have $(-1)^{n-1}a_n\in\Z^+$ unless $2n+1$ is a power of $5$ in which case
$5a_n\in\Z\sm5\Z$.

{\rm (ii)} Let $p$ be an odd prime.
Then
$$\align &\f1{p^a}\sum_{k=0}^{p^a-1}\f{154k+15}{(-2^{15})^k}\bi{6k}{3k}\bi{3k}{k,k,k}
\\\eq& 15\l(\f {-2}{p^a}\r)
+\l(\f{-2}{p^{a-1}}\r)\f{15}{16}p^2E_{p-3}\l(\f14\r)\ (\mo\ p^3)
\endalign$$
for any positive integer $a$.
We also have
 $$\align&\sum_{k=0}^{p-1}\f{\bi{6k}{3k}\bi{3k}{k,k,k}}{(-2^{15})^k}
 \\\eq&\cases
 (\f{-2}p)(x^2-2p)\ (\mo\ p^2)&\t{if}\ (\f p{11})=1\ \&\ 4p=x^2+11y^2\ (x,y\in\Z),
 \\0\ (\mo\ p^2)&\t{if}\ (\f p{11})=-1,\ i.e.,\ p\eq2,6,7,8,10\ (\mo\ 11).
 \endcases\endalign$$
 \endproclaim
\Remark.  The last congruence mod $p$ has been confirmed by Z. H. Sun [Su3].
Ramanujan [R] obtained that
$$\sum_{k=0}^\infty(28k+3)\f{(-27)^k}{2^{9k}}\cdot\f{(1/2)_k(1/6)_k(5/6)_k}{(1)_k^3}
=\sum_{k=0}^\infty\f{154k+15}{(-2^{15})^{k}}\bi{6k}{3k}\bi{3k}{k,k,k}=\f{32\sqrt2}{\pi}.$$

\proclaim{Conjecture A27 {\rm (i) ([S11b, S-3])}} If $p>5$ is a prime with
$p\eq1\ (\mo\ 4)$, then
$$\sum_{k=0}^{p^a-1}\f{k^3\bi{2k}k^3}{64^k}\eq0\ (\mo\ p^{2a})\quad\t{for all}\ a=1,2,3,\ldots.$$
Moreover, if $p$ is an odd prime, $m\in\Z^+$ and $p\eq m-1\ (\mo\ 2m)$, then
$$\sum_{k=0}^{p-1}(-1)^nk^n\bi{-1/m}k^n\eq0\pmod{p^2}$$
for all $n=3,5,\ldots,2(m-\da_{p,3m-1})-1$.

{\rm (ii) ([S11j])} Let $p>2$ be a prime and let $n\gs 2$ be an integer.
Assume that $x$ is a $p$-adic integer with $x\eq-2k\pmod{p}$ for some
$k\in\{1,\ldots,\lfloor(p+1)/(2n+1)\rfloor\}$. Then we have
$$\sum_{r=0}^{p-1}(-1)^r\bi xr^{2n+1}\eq0\pmod{p^2}.$$
\endproclaim
\Remark. Let $p$ be an odd prime. It is known that
$$\sum_{k=0}^{p-1}\f{\bi{2k}k^3}{64^k}\eq a(p)\ (\mo\ p^2),\ \
\ \t{where}\ \ \sum_{n=1}^\infty a(n)q^n=q\prod_{n=1}^\infty(1-q^{4n})^6.$$
This was proved by many authors, see, e.g., L. van Hammer [vH] and E. Mortenson [M2].
In 1892 F. Klein and R. Fricke proved that if $p=x^2+y^2$ with $x$ odd and $y$
even then $a(p)=4x^2-2p$ (see, e.g., Ishikawa [I]). The author
[S11b] determined $\sum_{k=0}^{p-1}k^3\bi{2k}k^3/64^k$ modulo $p$.
Recently Z. H. Sun [Su2] confirmed the first congruence in Conj. A27(i)
in the case $a=1$, and the author [S-3] proved the second congruence for $n=3$.
Concerning part (ii) the author [S11j] proved that if $p>3$ and $x$ is a $p$-adic integer with $x\eq-2k\pmod{p}$ for some
$k\in\{1,\ldots,\lfloor(p-1)/3\rfloor\}$ then
$$\sum_{r=0}^{p-1}(-1)^r\bi xr^3\eq0\pmod{p^2}.$$

\proclaim{Conjecture A28 {\rm ([S11e])}} Let $p$ be an odd prime.
If $(\f{-2}p)=1$ $($i.e., $p\eq 1,3\ (\mo\ 8))$ and $p=x^2+2y^2$ with $x,y\in\Z$, then
$$\l(\f{-1}p\r)\sum_{k=0}^{p-1}\f{\bi{2k}k^3}{(-64)^k}\eq\sum_{k=0}^{p-1}\f{\bi{4k}{k,k,k,k}}{28^{4k}}
\eq4x^2-2p\ (\mo\ p^2).$$
If $(\f{-2}p)=-1$ $($i.e., $p\eq 5,7\ (\mo\ 8))$, then
$$\sum_{k=0}^{p-1}\f{\bi{2k}k^3}{(-64)^k}
\eq0\ (\mo\ p^2),$$
and
$$\sum_{k=0}^{p-1}\f{\bi{4k}{k,k,k,k}}{28^{4k}}\eq 0\ (\mo\ p^2)\ \ \ \t{provided}\ p\not=7.$$
Also,
$$\f1{p^a}\sum_{k=0}^{p^a-1}\f{40k+3}{28^{4k}}\bi{4k}{k,k,k,k}
\eq 3\l(\f {p^a}3\r)-\l(\f{p^{a-1}}3\r)\f{5p^2}{392}B_{p-2}\l(\f13\r)\ (\mo\ p^3)$$
for all $a\in\Z^+$  provided that $p\not=3,7$.
Moreover,
$$a_n:=\f1{2n(2n+1)\bi{2n}n}\sum_{k=0}^{n-1}(40k+3)\bi{4k}{k,k,k,k} 28^{4(n-1-k)}$$
are integers for all $n=2,3,4,\ldots$.
\endproclaim
\Remark. (a) E. Mortenson [M4] proved the following conjecture of van Hamme [vH]:
$$\sum_{k=0}^{(p-1)/2}(4k+1)\bi{-1/2}k^3\eq (-1)^{(p-1)/2}p\ (\mo\ p^3)\quad\t{for any odd prime}\ p.$$
(Note that $\bi{-1/2}k^3=\bi{2k}k^3/(-64)^k$ for $k=0,1,2,\ldots$.)
See [S-2] for a further refinement. Recently Z. H. Sun [Su3] confirmed the conjecture that
$\sum_{k=0}^{p-1}\bi{2k}k^3/(-64)^k\eq0\ (\mo\ p^2)$ for any prime $p\eq5,7\ (\mo\ 8)$.
In [S11g] the author proved that
$$2n\bi{2n}n\ \bigg|\ \sum_{k=0}^{n-1}(4k+1)\bi{2k}k^3(-64)^{n-1-k}$$
for all $n=2,3,\ldots$.

(b) It is known (cf. [G5] and [G2]) that
$$\sum_{k=0}^\infty\f{4k+1}{(-64)^k}\bi{2k}k^3=\f2{\pi}\ \ \t{and}\ \ \sum_{k=1}^\infty\f{(4k-1)(-64)^k}{k^3\bi{2k}k^3}=-16G,$$
where $G:=\sum_{k=0}^\infty(-1)^k/(2k+1)^2$ is the Catalan constant. Ramanujan (cf. [Be, p.\,354]) found that
$$\sum_{k=0}^\infty(40k+3)\f{(1/2)_k(1/4)_k(3/4)_k}{k!^3 7^{4k}}=
\sum_{k=0}^\infty\f{40k+3}{28^{4k}}\bi{4k}{k,k,k,k}=\f{49}{3\pi\sqrt3}.$$

\proclaim{Conjecture A29 {\rm ([S11e])}} Let $p$ be an odd prime. If $p\eq1\pmod3$, then
$$\sum_{k=0}^{p-1}\f{9k+2}{108^k}\bi{2k}k^2\bi{3k}k\eq0\pmod{p^2}.$$
If $p\eq1,3\ (\mo\ 8)$, then
$$\sum_{k=0}^{p-1}\f{16k+3}{256^k}\bi{4k}{k,k,k,k}\eq0\ (\mo\ p^2).$$
If $p\eq1\ (\mo\ 4)$, then
$$\sum_{k=0}^{p-1}\f{36k+5}{12^{3k}}\bi{6k}{3k}\bi{3k}{k,k,k}\eq 0\ (\mo\ p^2).$$
\endproclaim
\Remark. A related conjecture of Rodriguez-Villegas [RV] confirmed by Mortenson [M2] and the author [S11i] together
states that for any prime $p>3$ we have
$$\sum_{k=0}^{p-1}\f{\bi{2k}k^2\bi{3k}k}{108^k}\eq\cases 4x^2-2p\ (\mo\ p^2)&\t{if}
\ (\f p3)=1\ \&\ p=x^2+3y^2\ (x,y\in\Z),
\\0\ (\mo\ p^2)&\t{if}\ p\eq2\ (\mo\ 3),\endcases$$
$$\sum_{k=0}^{p-1}\f{\bi{4k}{k,k,k,k}}{256^k}\eq\cases 4x^2-2p\ (\mo\ p^2)&\t{if}
\ (\f{-2}p)=1\ \&\ p=x^2+2y^2\ (x,y\in\Z),
\\0\ (\mo\ p^2)&\t{if}\ (\f{-2}p)=-1,\ \t{i.e.,}\ p\eq 5,7\ (\mo\ 8),\endcases$$
and
$$\align &\sum_{k=0}^{p-1}\f{\bi{6k}{3k}\bi{3k}{k,k,k}}{12^{3k}}=\sum_{k=0}^{p-1}\f{(6k)!}{(3k)!(k!)^3}1728^{-k}
\\\eq&\cases(\f p3)(4x^2-2p)\ (\mo\ p^2)&\t{if}\ 4\mid p-1\ \&\ p=x^2+y^2\ (2\nmid x,\ 2\mid y),
\\0\ (\mo\ p^2)&\t{if}\ p\eq3\ (\mo\ 4).\endcases
\endalign$$
Z. H. Sun [Su3] confirmed the congruences in Conj. A29 modulo $p$.

\proclaim{Conjecture A30 {\rm ([S11e, S-7])}} Let $p>5$ be a prime.
If $p>7$ then
$$\sum_{k=1}^{p-1}\f{\bi{2k}k}{k^3}\eq-\f 2{p^2}H_{p-1}-\f{13}{27}\sum_{k=1}^{p-1}\f1{k^3}\ (\mo\ p^4).$$
We also have
$$\gather\sum_{k=1}^{p-1}\f1{k^4\bi{2k}k}-\f{H_{p-1}}{p^3}\eq \f 7{54p}\sum_{k=1}^{p-1}\f1{k^3}
\eq-\f 7{45}pB_{p-5}\ (\mo\ p^2),
\\\sum_{k=1}^{(p-1)/2}\f{(-1)^k}{k^3\bi{2k}k}\eq-2B_{p-3}\ (\mo\ p),
\\\sum_{k=1}^{(p-1)/2}\f{(-1)^k}{k^2}\bi{2k}k\eq\f{56}{15}pB_{p-3}\ (\mo\ p^2),
\\\sum_{p/2<k<p}\f{\bi{2k}k^2}{k16^k}\eq-\f {21}2H_{p-1}\ (\mo\ p^4).
\endgather$$
\endproclaim
\Remark. It is known that $H_{p-1}/p^2\eq-B_{p-3}/3\ (\mo\ p)$ for any prime $p>3$
and $\sum_{k=1}^{p-1}1/k^3\eq-\f 65 p^2B_{p-5}\ (\mo\ p^3)$ for each prime $p>5$, and that
$$\sum_{k=1}^\infty\f{(-1)^k}{k^3\bi{2k}k}=-\f25\zeta(3)
\ \ \t{and}\ \ \sum_{k=1}^\infty\f1{k^4\bi{2k}k}=\f{17}{36}\zeta(4).$$
The first and the second congruences in Conj. A30 modulo $p$ have been confirmed by
K. Hessami Pilehrood and T. Hessami Pilehrood [HP3].
Also, [T1, Theorem 4.2] implies
that $\sum_{k=1}^{p-1}\f{(-1)^k}{k^2}\bi{2k}k\eq-\f 4{15}pB_{p-3}\ (\mo\ p^2)$
for any prime $p>5$. Tauraso [T2] proved that $\sum_{k=1}^{p-1}\bi{2k}k^2/(k16^k)\eq-2H_{(p-1)/2}\ (\mo\ p^3)$
for each prime $p>3$. The author could prove that the third and the fourth
congruences in Conj. A30 are equivalent. By [S11e],
$$\sum_{k=1}^{(p-1)/2}\f{\bi{2k}k}k\eq(-1)^{(p+1)/2}\,\f 83pE_{p-3}\ (\mo\ p^2)$$
and
$$\sum_{k=1}^{(p-1)/2}\f1{k^2\bi{2k}k}\eq(-1)^{(p-1)/2}\f 43E_{p-3}\ (\mo\ p)$$
for any prime $p>3$.  {\tt Mathematica} (version 7) yields
$$\sum_{k=1}^\infty\f{\bi{2k}k^2}{k16^k}=4\log2-\f{8G}{\pi}$$
where $G=\sum_{k=0}^\infty(-1)^k/(2k+1)^2$ is the Catalan constant.

\proclaim{Conjecture A31 {\rm ([S11b, S-7])}} Let $p$ be an odd prime. Then
$$\sum_{k=1}^{p-1}\f{\bi{2k}k}{k2^k}\eq-\f{H_{(p-1)/2}}2+\f7{16}p^2B_{p-3}\ (\mo\ p^3).$$
When $p>3$, we have
$$\sum_{k=1}^{p-1}\f{\bi{2k}k}{k3^k}\eq-2\sum^{p-1}\Sb k=1\\k\not\eq p\,(\mo\ 3)\endSb\f1k\ (\mo\ p^3).$$
If $p>5$, then
$$\sum_{k=1}^{p-1}\f{\bi{2k}k}{k^24^k}\eq-\f{H_{(p-1)/2}^2}2-\f 74\cdot\f{H_{p-1}}{p}\ (\mo\ p^3).$$
$$$$
\endproclaim
\Remark. The congruences were motivated by the following known identities (cf. [Ma] or using {\tt Mathematica}):
$$\sum_{k=1}^\infty\f{2^k}{k^2\bi{2k}k}=\f{\pi^2}8,\ \ \sum_{k=1}^\infty\f{3^k}{k^2\bi{2k}k}=\f29\pi^2,\ \ \sum_{k=1}^\infty\f{\bi{2k}k}{k^24^k}=\f{\pi^2-3\log^24}6.$$

\proclaim{Conjecture A32 {\rm ([S11b], [S-8])}} Let $p>5$ be a prime
and let $H_{p-1}=\sum_{k=1}^{p-1}1/k$.
Then
$$\align\sum_{k=0}^{(p-3)/2}\f{\bi{2k}k}{(2k+1)16^k}\eq&(-1)^{(p-1)/2}\l(\f{H_{p-1}}{12}+\f {3p^4}{160}B_{p-5}\r)\ (\mo\ p^5),
\\\sum_{k=0}^{(p-3)/2}\f{\bi{2k}k}{(2k+1)^3 16^k}\eq&(-1)^{(p-1)/2}\l(\f{H_{p-1}}{4p^2}+\f{p^2}{36}B_{p-5}\r)\ (\mo\ p^3).
\endalign$$
We also have
$$\gather\sum_{k=0}^{(p-3)/2}\f{\bi{2k}k}{(2k+1)^2(-16)^k}\eq\f{H_{p-1}}{5p}\ (\mo\ p^3),
\\\sum_{p/2<k<p}\f{\bi{2k}k}{(2k+1)^2(-16)^k}\eq-\f p4B_{p-3}\ (\mo\ p^2),
\\\sum_{k=0}^{(p-3)/2}\f{\bi{2k}k\bar H_k^{(2)}}{(2k+1)16^k}\eq\l(\f{-1}p\r)\f{H_{p-1}}{12p^2}\ (\mo\ p^2),
\endgather$$
where $\bar H_k^{(2)}:=\sum_{0<j\ls k}1/(2j-1)^2$.
\endproclaim
\Remark. On March 6, 2010 the author [S11b] proved the first congruence in Conj. A32 mod $p^2$.
It is known (cf. [Ma]) that
$$\sum_{k=0}^\infty\f{\bi{2k}k}{(2k+1)16^k}=\f{\pi}3\ \ \ \t{and}
\ \ \ \sum_{k=0}^\infty\f{\bi{2k}k}{(2k+1)^2(-16)^k}=\f{\pi^2}{10}$$
which can be easily proved by using $1/(2k+1)=\int_0^1 x^{2k}dx$.
In March 2010 the author suggested that $\sum_{k=0}^\infty\bi{2k}k/((2k+1)^316^k)=7\pi^3/216$
via a public message to Number Theory List,
and then Olivier Gerard pointed out there is a computer proof via certain math.
softwares like {\tt Mathematica} (Version 7).
It is also known that (see [S-8], and use {\tt Mathematica} 7)
$$\sum_{k=1}^\infty\f{\bi{2k}k\bar H_k^{(2)}}{(2k+1)16^k}=\f{\pi^3}{648}.$$
Using Morley's congruence and the identity
$$\sum_{k=0}^n\bi nk\bi{n+k}k\f{(-1)^k}{2k+1}=\f1{2n+1},$$
the author got the following congruence (for primes $p>3$)
$$\sum_{k=0}^{(p-3)/2}\f{\bi{2k}k^2}{(2k+1)16^k}\eq-2q_p(2)-p\,q_p(2)^2\ (\mo\ p^2).$$

\proclaim{Conjecture A33 {\rm ([S11e])}} Let $p$ be an odd prime and let $a\in\Z^+$.
If $p\eq1\ (\mo\ 4)$ or $a>1$, then
$$\sum_{k=0}^{\lfloor \f34p^a\rfloor}\f{\bi{2k}k^2}{16^k}\eq (-1)^{(p^a-1)/2}\ (\mo\ p^3).$$
If $p>3$, and $p\eq1,3\ (\mo\ 8)$ or $a>1$, then
$$\sum_{k=0}^{\lfloor \f r8p^a\rfloor}\f{\bi{2k}k^2}{16^k}\eq (-1)^{(p^a-1)/2}\ (\mo\ p^3)
\ \ \ \t{for}\ r=5,7.$$
\endproclaim
\Remark. The author [S11e] showed that
$\sum_{k=0}^{\lfloor p/2\rfloor}\bi{2k}k^2/16^k\eq (-1)^{(p-1)/2}+p^2E_{p-3}\ (\mo\ p^3)$ for any odd prime $p$.

\proclaim{Conjecture A34}  {\rm (i) ([S11e])} For each $n=2,3,\ldots$ we have
$$\gather\f1{2n\bi{2n}n}\sum_{k=0}^{n-1}(6k+1)\bi{2k}k^3256^{n-1-k}\in\Z,
\\\f1{2n\bi{2n}n}\sum_{k=0}^{n-1}(6k+1)\bi{2k}k^3(-512)^{n-1-k}\in\Z,
\\\f1{2n\bi{2n}n}\sum_{k=0}^{n-1}(42k+5)\bi{2k}k^34096^{n-1-k}\in\Z.
\endgather$$

{\rm (ii) ([S11e])} Let $p>3$ be a prime.
We have
$$\f1{p^a}\sum_{k=0}^{(p^a-1)/2}\f{42k+5}{4096^k}\bi{2k}k^3\eq\l(\f{-1}{p^a}\r)\l(5-\f34pH_{p-1}\r)\ (\mo\ p^4)$$
for all $a\in\Z^+$.
Also,
$$\align\sum_{k=0}^{p-1}\f{6k+1}{256^k}\bi{2k}k^3\eq& p\l(\f{-1}p\r)-p^3E_{p-3}\ (\mo\ p^4),
\\\sum_{k=0}^{(p-1)/2}\f{6k+1}{(-512)^k}\bi{2k}k^3\eq&p\l(\f{-2}p\r)+\f{p^3}4\l(\f2p\r)E_{p-3}\ (\mo\ p^4),
\\\sum_{k=0}^{p-1}\f{42k+5}{4096^k}\bi{2k}k^3\eq& 5p\l(\f{-1}p\r)-p^3E_{p-3}\ (\mo\ p^4).
\endalign$$

{\rm (iii) ([S-2])} For any prime $p>3$ we have
$$\sum_{k=0}^{p-1}\f{7k+1}{648^k}\bi{4k}{k,k,k,k}\eq p\l(\f{-1}p\r)-\f 53p^3E_{p-3}\ (\mo\ p^4).$$
Also, for $n=2,3,\ldots$ we have
$$\f1{2n(2n+1)\bi{2n}n}\sum_{k=0}^{n-1}(7k+1)\bi{4k}{k,k,k,k}648^{n-1-k}\in\Z$$
unless $2n+1$ is a power of $3$ in which case the denominator of the quotient is $3$.
\endproclaim
\Remark. Those congruences in part (ii) mod $p^3$ are van Hamme's conjectures (cf. [vH])
which are $p$-adic analogues of corresponding
Ramanujan series. For the congruence in (iii) the corresponding Ramanujan series is
$$\sum_{k=0}^\infty\f{7k+1}{648^k}\bi{4k}{k,k,k,k}=\f 9{2\pi}.$$

\proclaim{Conjecture A35 {\rm ([S11b, S11e])}} Let $p$ be an odd prime.
If $p\eq1\ (\mo\ 4)$, then
$$\sum_{k=0}^{p-1}\f{\bi{2k}k^3}{(-8)^k}\l(1-\f1{(-8)^k}\r)\eq0\ (\mo\ p^3).$$
If $p\eq3\ (\mo\ 4)$, then
$$\sum_{k=0}^{p-1}\f{\bi{2k}k^2}{8^k}\l(1+\f1{(-2)^k}\r)\eq0\ (\mo\ p^3)$$
and
$$\sum_{k=0}^{p-1}\bi{p-1}k\f{\bi{2k}k^2}{(-8)^k}\eq0\ (\mo\ p^2).$$
\endproclaim
\Remark. When $p>3$ is a prime congruent to 3 mod 4, the author [S11b] proved  that
$$\sum_{k=0}^{p-1}\bi{p-1}k\f{\bi{2k}k^3}{(-64)^k}\eq0\ (\mo\ p^2).$$

\proclaim{Conjecture A36 {\rm ([S11e])}} Let $p$ be an odd prime.

{\rm (i)} We have
$$\align \sum_{k=0}^{p-1}\f{\bi{2k}k^3}{(-8)^k}\eq&\l(\f{-2}p\r)\sum_{k=0}^{p-1}\f{\bi{2k}k^3}{(-512)^k}\ (\mo\ p^2),
\\\sum_{k=0}^{p-1}\f{\bi{2k}k^3}{16^k}\eq&\l(\f{-1}p\r)\sum_{k=0}^{p-1}\f{\bi{2k}k^3}{256^k}\ (\mo\ p^2).
\endalign$$
Moreover, if $p\eq1\ (\mo\ 4)$ then
$$\sum_{k=0}^{p-1}\f{\bi{2k}k^3}{(-8)^k}\eq\l(\f{2}p\r)\sum_{k=0}^{p-1}\f{\bi{2k}k^3}{(-512)^k}
\eq\sum_{k=0}^{p-1}\f{\bi{4k}{k,k,k,k}}{648^k}\ (\mo\ p^3);$$ if
$p\eq1\ (\mo\ 3)$ then
$$\sum_{k=0}^{p-1}\f{\bi{2k}k^3}{16^k}\eq\l(\f{-1}p\r)\sum_{k=0}^{p-1}\f{\bi{2k}k^3}{256^k}
\eq\sum_{k=0}^{p-1}\f{\bi{4k}{k,k,k,k}}{(-144)^k}\ (\mo\ p^3).$$

{\rm (ii)} If
$p\eq1,3\ (\mo\ 8)$, then
$$\sum_{k=0}^{p-1}\f{\bi{2k}k^3}{(-64)^k}
\eq\l(\f{-1}p\r)\sum_{k=0}^{p-1}\f{\bi{4k}{k,k,k,k}}{256^k}\ (\mo\ p^3);$$
if $p\eq1,2,4\ (\mo\ 7)$, then
$$\l(\f{-1}p\r)\sum_{k=0}^{p-1}\bi{2k}k^3\eq\sum_{k=0}^{p-1}\f{\bi{2k}k^3}{4096^k}
\eq\sum_{k=0}^{p-1}\f{\bi{4k}{k,k,k,k}}{(-3969)^k}\ (\mo\ p^3).$$
\endproclaim
\Remark. The author observed that
$$\align\l(\f xp\r)\sum_{k=0}^{p-1}\f{\bi{2k}k^3}{(-64x)^k}
\eq&\sum_{k=0}^{p-1}\bi{2k}k^3\l(\f x{-64}\r)^k
\\\eq&\l(\f{x+1}p\r)\sum_{k=0}^{p-1}\bi{4k}{k,k,k,k}\l(\f x{64(x+1)^2}\r)^k
\ (\mo\ p)\endalign$$ for any odd prime $p$ and $p$-adic integer $x\not\eq0,-1\ (\mo\ p)$;
this can be easily deduced by taking $n=(p-1)/2$ in the
known identity
$$\sum_{k=0}^n\bi nk^3x^k=\sum_{k=0}^n\bi nk\bi{n+k}k\bi{n-k}kx^k(1+x)^k.$$
Z. H. Sun [Su2] proved that $\sum_{k=0}^{p-1}\bi{2k}k^3/256^k\eq0\ (\mo\ p^2)$ for any prime $p\eq5\ (\mo\ 6)$,
and $\sum_{k=0}^{p-1}\bi{2k}k^3/(-512)^k\eq0\ (\mo\ p^2)$ for any prime $p\eq3\ (\mo\ 4)$. See also
Z. H. Sun [Su1] for his conjectures on
$\sum_{k=0}^{p-1}\bi{4k}{k,k,k,k}/m^k$ mod $p^2$ with $m=-144,\,648,\,-3969$ motivated by the author's papers [S11b] and [S11e].

\proclaim{Conjecture A37 {\rm (Discovered on March 2, 2010)}} Let $p$ be an odd prime.

{\rm (i)} If $p\eq1\ (\mo\ 4)$ then
$$\align\sum_{k=0}^{(p-1)/2}\f{\bi{2k}k^3}{(-8)^k}\sum_{k<j\ls 2k}\f1j
\eq&\f12\sum_{k=0}^{(p-1)/2}\f{\bi{2k}k^3}{64^k}\sum_{k<j\ls 2k}\f1j
\\\eq&\f 13\l(\f 2p\r)\sum_{k=0}^{(p-1)/2}\f{\bi{2k}k^3}{(-512)^k}\sum_{k<j\ls 2k}\f1j\ (\mo\ p^2);
\endalign$$
when $p\eq3\ (\mo\ 4)$ we have
$$\align\sum_{k=0}^{(p-1)/2}\f{\bi{2k}k^3}{(-8)^k}\sum_{k<j\ls 2k}\f1j
\eq&-\f 72\sum_{k=0}^{(p-1)/2}\f{\bi{2k}k^3}{64^k}\sum_{k<j\ls 2k}\f1j\ (\mo\ p^2),
\\\sum_{k=0}^{(p-1)/2}\f{\bi{2k}k^3}{64^k}\sum_{k<j\ls 2k}\f1j
\eq&-\l(\f 2p\r)\sum_{k=0}^{(p-1)/2}\f{\bi{2k}k^3}{(-512)^k}\sum_{k<j\ls 2k}\f1j\ (\mo\ p^2),
\endalign$$
and
$$\sum_{k=0}^{(p-1)/2}\f{\bi{2k}k^3}{m^k}\sum_{k<j\ls 2k}\f1j\eq0\ (\mo\ p)\ \t{for}\ m=-8,-512\ \ \t{if}\ p>3.$$

{\rm (ii)} If $p\eq1\ (\mo\ 3)$ then
$$\sum_{k=0}^{(p-1)/2}\f{\bi{2k}k^3}{16^k}\sum_{k<j\ls 2k}\f1j
\eq\f12\l(\f{-1}p\r)\sum_{k=0}^{(p-1)/2}\f{\bi{2k}k^3}{256^k}\sum_{k<j\ls 2k}\f1j\ (\mo\ p^2).$$
If $p\eq 2\ (\mo\ 3)$ then
$$\sum_{k=0}^{(p-1)/2}\f{\bi{2k}k^3}{16^k}\sum_{k<j\ls 2k}\f1j\eq0\ (\mo\ p)
\ \t{and}\ \sum_{k=0}^{(p-1)/2}\f{\bi{2k}k^3}{256^k}\sum_{k<j\ls 2k}\f1j\eq0\ (\mo\ p^2).$$

{\rm (iii)} If $p\eq 5,7\ (\mo\ 8)$, then
$$\sum_{k=0}^{(p-1)/2}\f{\bi{2k}k^3}{(-64)^k}\sum_{k<j\ls 2k}\f1j\eq0\ (\mo\ p).$$

{\rm (iv)} If $p>3$ and $p\eq3,5,6\pmod7$, then
$$\sum_{k=0}^{(p-1)/2}\bi{2k}k^3\sum_{k<j\ls 2k}\f1j\eq0\pmod{p^2}.$$
\endproclaim

\Remark.  During their attempt to prove Conjecture A1, M. Jameson and K. Ono
realized that $\sum_{k=0}^{(p-1)/2}\bi{2k}k^3\sum_{k<j\ls 2k}1/j\eq0\pmod{p}$ for any prime $p>3$
but they did not have a proof of this observation.
When $p>3$ is a prime with $p\eq 3\ (\mo\ 4)$, the author was able to show that
$$\sum_{k=0}^{(p-1)/2}\f{\bi{2k}k^3}{64^k}H_{2k}\eq\sum_{k=0}^{(p-1)/2}\f{\bi{2k}k^3}{64^k}H_{k}\eq0\ (\mo\ p).$$

\medskip

Now we introduce polynomials
$$S_n(x):=\sum_{k=0}^n\bi nk^4x^k\quad (n=0,1,2,\ldots).$$

\proclaim{Conjecture A38 {\rm ([S-12])}} Let $p$ be an odd prime.

{\rm (i)} We have
$$\align&\sum_{k=0}^{p-1}S_k(12)
\\\eq&\cases(-1)^{[3\mid x]}(4x^2-2p)\pmod{p^2}&\t{if}\ p\eq1\pmod{12}\ \&\ p=x^2+y^2\ (4\mid x-1),
\\(\f{xy}3)4xy\pmod{p^2}&\t{if}\ p\eq5\pmod{12}\ \&\ p=x^2+y^2\ (4\mid x-1),
\\0\pmod{p^2}&\t{if}\ p\eq3\pmod 4.\endcases
\endalign$$
And
$$\sum_{k=0}^{p-1}(4k+3)S_k(12)\eq p\l(1+2\l(\f 3p\r)\r)\pmod {p^2}.$$
Moreover,
$$\f1{n}\sum_{k=0}^{n-1}(4k+3)S_k(12)\in\Z
\qquad\t{for all}\ n=1,2,3,\ldots.$$

{\rm (ii)} We have
$$\align&\sum_{k=0}^{p-1}S_k(-20)
\\\eq&\cases(-1)^{[5\mid x]}(4x^2-2p)\pmod{p^2}&\t{if}\ (\f{-1}p)=(\f p5)=1\ \&\ p=x^2+y^2\ (2\nmid x),
\\(\f{x}5)4xy\pmod{p^2}&\t{if}\ (\f{-1}p)=-(\f p5)=1,\  p=x^2+y^2\ (2\nmid x,\ 5\mid xy-1),
\\0\pmod{p^2}&\t{if}\ p\eq3\pmod 4.\endcases
\endalign$$
And
$$\sum_{k=0}^{p-1}(6k+5)S_k(-20)\eq p\l(\f{-1}p\r)\l(2+3\l(\f {-5}p\r)\r)\pmod {p^2}.$$
Moreover,
$$\f1{n}\sum_{k=0}^{n-1}(6k+5)S_k(-20)\in\Z
\qquad\t{for all}\ n=1,2,3,\ldots.$$
\endproclaim

\proclaim{Conjecture A39 {\rm ([S-12])}} Let $p$ be an odd prime.

{\rm (i)} We have
$$\align&\sum_{k=0}^{p-1}S_k(36)
\\\eq&\cases 4x^2-2p\pmod{p^2}&\t{if}\ (\f 2p)=(\f p3)=(\f p5)=1\ \&\ p=x^2+30y^2,
\\12x^2-2p\pmod{p^2}&\t{if}\ (\f p3)=1,\ (\f 2p)=(\f p5)=-1\ \&\ p=3x^2+10y^2,
\\8x^2-2p\pmod{p^2}&\t{if}\ (\f 2p)=1,\ (\f p3)=(\f p5)=-1\ \&\ p=2x^2+15y^2,
\\2p-20x^2\pmod{p^2}&\t{if}\ (\f p5)=1,\ (\f 2p)=(\f p3)=-1\ \&\ p=5x^2+6y^2,
\\0\pmod{p^2}&\t{if}\ (\f{-30}p)=-1.\endcases
\endalign$$
And
$$\sum_{k=0}^{p-1}(8k+7)S_k(36)\eq p\l(\f p{15}\r)\l(3+4\l(\f{-6}p\r)\r)\pmod {p^2}.$$
We also have
$$\f1{n}\sum_{k=0}^{n-1}(8k+7)S_k(36)\in\Z
\qquad\t{for all}\ n=1,2,3,\ldots.$$

{\rm (ii)} We have
$$\align&\sum_{k=0}^{p-1}S_k(196)
\\\eq&\cases 4x^2-2p\pmod{p^2}&\t{if}\ (\f 2p)=(\f p5)=(\f p7)=1\ \&\ p=x^2+70y^2,
\\8x^2-2p\pmod{p^2}&\t{if}\ (\f p7)=1,\ (\f 2p)=(\f p5)=-1\ \&\ p=2x^2+35y^2,
\\2p-20x^2\pmod{p^2}&\t{if}\ (\f p5)=1,\ (\f 2p)=(\f p7)=-1\ \&\ p=5x^2+14y^2,
\\28x^2-2p\pmod{p^2}&\t{if}\ (\f 2p)=1,\ (\f p5)=(\f p7)=-1\ \&\ p=7x^2+10y^2,
\\0\pmod{p^2}&\t{if}\ (\f{-70}p)=-1.\endcases
\endalign$$
And
$$\sum_{k=0}^{p-1}(120k+109)S_k(196)\eq p\l(\f p{7}\r)\l(49+60\l(\f{-14}p\r)\r)\pmod {p^2}.$$
We also have
$$\f1{n}\sum_{k=0}^{n-1}(120k+109)S_k(196)\in\Z
\qquad\t{for all}\ n=1,2,3,\ldots.$$

{\rm (iii)} We have
$$\align&\sum_{k=0}^{p-1}S_k(-324)
\\\eq&\cases 4x^2-2p\pmod{p^2}&\t{if}\ (\f {-1}p)=(\f p5)=(\f p{17})=1\ \&\ p=x^2+85y^2,
\\2x^2-2p\pmod{p^2}&\t{if}\ (\f {p}{17})=1,\ (\f {-1}p)=(\f p5)=-1\ \&\ 2p=x^2+85y^2,
\\2p-20x^2\pmod{p^2}&\t{if}\ (\f {-1}p)=1,\ (\f p5)=(\f p{17})=-1\ \&\ p=5x^2+17y^2,
\\2p-10x^2\pmod{p^2}&\t{if}\ (\f p{5})=1,\ (\f {-1}p)=(\f p{17})=-1\ \&\ 2p=5x^2+17y^2,
\\0\pmod{p^2}&\t{if}\ (\f{-85}p)=-1.\endcases
\endalign$$
Provided $p>3$ we have
$$\sum_{k=0}^{p-1}(34k+31)S_k(-324)\eq p\l(\f p5\r)\l(17+14\l(\f{-1}p\r)\r)\pmod {p^2}.$$
Moreover,
$$\f1{n}\sum_{k=0}^{n-1}(34k+31)S_k(-324)\in\Z
\qquad\t{for all}\ n=1,2,3,\ldots.$$

{\rm (iv)} We have
$$\align&\sum_{k=0}^{p-1}S_k(1296)
\\\eq&\cases 4x^2-2p\pmod{p^2}&\t{if}\ (\f {-2}p)=(\f p5)=(\f p{13})=1\ \&\ p=x^2+130y^2,
\\8x^2-2p\pmod{p^2}&\t{if}\ (\f {-2}p)=1,\ (\f p5)=(\f p{13})=-1\ \&\ p=2x^2+65y^2,
\\2p-20x^2\pmod{p^2}&\t{if}\ (\f p5)=1,\ (\f {-2}p)=(\f p{13})=-1\ \&\ p=5x^2+26y^2,
\\2p-40x^2\pmod{p^2}&\t{if}\ (\f p{13})=1,\ (\f {-2}p)=(\f p5)=-1\ \&\ p=10x^2+13y^2,
\\0\pmod{p^2}&\t{if}\ (\f{-130}p)=-1.\endcases
\endalign$$
Provided $p>3$ we have
$$\sum_{k=0}^{p-1}(130k+121)S_k(1296)\eq p\l(\f{-2} p\r)\l(56+65\l(\f{-26}p\r)\r)\pmod {p^2}.$$
Moreover,
$$\f1{n}\sum_{k=0}^{n-1}(130k+121)S_k(1296)\in\Z
\qquad\t{for all}\ n=1,2,3,\ldots.$$

{\rm (v)} We have
$$\align&\sum_{k=0}^{p-1}S_k(5776)
\\\eq&\cases 4x^2-2p\pmod{p^2}&\t{if}\ (\f {2}p)=(\f p5)=(\f p{19})=1\ \&\ p=x^2+190y^2,
\\8x^2-2p\pmod{p^2}&\t{if}\ (\f {2}p)=1,\ (\f p5)=(\f p{19})=-1\ \&\ p=2x^2+95y^2,
\\2p-20x^2\pmod{p^2}&\t{if}\ (\f p{19})=1,\ (\f {2}p)=(\f p{5})=-1\ \&\ p=5x^2+38y^2,
\\2p-40x^2\pmod{p^2}&\t{if}\ (\f p{5})=1,\ (\f {2}p)=(\f p{19})=-1\ \&\ p=10x^2+19y^2,
\\0\pmod{p^2}&\t{if}\ (\f{-190}p)=-1.\endcases
\endalign$$
And
$$\sum_{k=0}^{p-1}(816k+769)S_k(5776)\eq p\l(\f p{95}\r)\l(361+408\l(\f{p}{19}\r)\r)\pmod {p^2}.$$
Moreover,
$$\f1{n}\sum_{k=0}^{n-1}(816k+769)S_k(5776)\in\Z
\qquad\t{for all}\ n=1,2,3,\ldots.$$
\endproclaim
\Remark. The reader may consult [S-12] for more conjectures of this type.
I'd like to offer \$300 (300 US dollars) for the first correct proof of part (i) of Conj. A39.

\proclaim{Conjecture A40 {\rm (Discovered on April 5, 2010)}}
{\rm (i)} For any odd prime $p$ we have
$$\sum_{k=0}^{p-1}\f{28k^2+18k+3}{(-64)^k}\bi{2k}k^4\bi{3k}k\eq 3p^2-\f 72p^5B_{p-3}\ (\mo\ p^6)$$
and
$$\sum_{k=0}^{(p-1)/2}\f{28k^2+18k+3}{(-64)^k}\bi{2k}k^4\bi{3k}k\eq 3p^2+6\l(\f{-1}p\r)p^4E_{p-3}\ (\mo\ p^5).$$

{\rm (ii)} For any integer $n>1$, we have
$$\sum_{k=0}^{n-1}(28k^2+18k+3)\bi{2k}k^4\bi{3k}k(-64)^{n-1-k}\eq0\ \ \(\mo\ (2n+1)n^2\bi{2n}n^2\).$$
Also,
$$\sum_{k=1}^\infty\f{(28k^2-18k+3)(-64)^k}{k^5\bi{2k}k^4\bi{3k}k}=-14\zeta(3).$$
\endproclaim

\proclaim{Conjecture A41 {\rm (Discovered on April 5, 2010)}} Let $p$ be an odd prime.

{\rm (i)} If $p\not=3$, then
$$\sum_{k=0}^{p-1}\f{10k^2+6k+1}{(-256)^k}\bi{2k}k^5\eq p^2-\f 76 p^5B_{p-3}\ (\mo\ p^6)$$
and
$$\sum_{k=0}^{(p-1)/2}\f{10k^2+6k+1}{(-256)^k}\bi{2k}k^5\eq p^2+\f 73 p^5B_{p-3}\ (\mo\ p^6).$$

{\rm (ii)} If $p\not=5$, then
$$\sum_{k=0}^{p-1}\f{74k^2+27k+3}{4096^{k}}\bi{2k}k^4\bi{3k}k\eq 3p^2+7p^5B_{p-3}\ (\mo\ p^6)$$
and
$$\sum_{k=0}^{(p-1)/2}\f{74k^2+27k+3}{4096^k}\bi{2k}k^4\bi{3k}k\eq 3p^2-\f 94p^3H_{p-1}\ (\mo\ p^7).$$
\endproclaim
\Remark. By [G3, Identity 8] and [G4], we have
$$\sum_{k=1}^\infty\f{(10k^2-6k+1)(-256)^k}{k^5\bi{2k}k^5}=-28\zeta(3)$$
and $$\sum_{k=0}^\infty\f{74k^2+27k+3}{4096^k}\bi{2k}k^4\bi{3k}k=\f{48}{\pi^2}.$$

\proclaim{Conjecture A42 {\rm (Discovered on April 6, 2010)}}
{\rm (i)} For any prime $p\not=2,5$ we have
$$\sum_{k=0}^{p-1}\f{21k^3+22k^2+8k+1}{256^k}\bi{2k}k^7\eq p^3\ (\mo\ p^8)$$
and
$$\sum_{k=0}^{(p-1)/2}\f{168k^3+76k^2+14k+1}{2^{20k}}\bi{2k}k^7\eq \l(\f{-1}p\r)p^3\ (\mo\ p^8).$$

{\rm (ii)} For any integer $n>1$, we have
$$\sum_{k=0}^{n-1}(21k^3+22k^2+8k+1)\bi{2k}k^7256^{n-1-k}\eq0\ \ \(\mo\ 2n^3\bi{2n}n^3\)$$
and
$$\sum_{k=0}^{n-1}(168k^3+76k^2+14k+1)\bi{2k}k^72^{20(n-1-k)}\eq0\ \ \(\mo\ 2n^3\bi{2n}n^3\).$$
\endproclaim
\Remark. (a) B. Gourevich and Guillera (see [G1, Section 4]) conjectured
$$\sum_{k=0}^\infty\f{168k^3+76k^2+14k+1}{2^{20k}}\bi{2k}k^7=\f{32}{\pi^3}$$
and
$$\sum_{k=1}^\infty\f{(21k^3-22k^2+8k-1)256^k}{k^7\bi{2k}k^7}=\f{\pi^4}8$$
respectively. Zudilin [Zu] suggested that for any odd prime $p$ we might have
$$\sum_{k=0}^{p-1}\f{168k^3+76k^2+14k+1}{2^{20k}}\bi{2k}k^7\eq \l(\f{-1}p\r)p^3\ (\mo\ p^7),$$
which is weaker than the second congruence in Conj. A42(i).

(b) Let $a_1=2$ and
$$(2n+1)^3a_n=32n^3a_{n-1}+(21n^3+22n^2+8n+1)\bi{2n-1}n^4\ \ \t{for}\ n=2,3,\ldots.$$
Then for any $n\in\Z^+$ we have
$$\align a_n=&\f1{16(2n+1)^3\bi{2n}n^3}\sum_{k=0}^n(21k^3+22k^2+8k+1)256^{n-k}\bi{2k}k^7
\\=&\f1{2(n+1)^3\bi{2(n+1)}{n+1}^3}\sum_{k=0}^n(21k^3+22k^2+8k+1)256^{n-k}\bi{2k}k^7.
\endalign$$
The author created the sequence $\{a_n\}_{n>0}$ at OEIS as A176477 (cf. [S]).
We not only conjectured that $a_n\in\Z^+$ for all $n=1,2,3,\ldots$ but also guessed that $a_n$ is odd if and only
if $n=2^k$ for some $k\in\Z^+$.  We have a similar conjecture related to the last congruence in Conj. A42.

\proclaim{Conjecture A43 {\rm ([S11e])}}
{\rm (i)} For any prime $p>3$ we have
$$\sum_{n=0}^{p-1}\f{3n^2+n}{16^n}\sum_{k=0}^n\bi nk^2\bi{2k}k\bi{2(n-k)}{n-k}\eq-4p^4q_p(2)+6p^5q_p(2)^2\pmod{p^6}.$$

{\rm (ii)} Let $m>1$ be an integer. Then
$$a_m:=\f1{2m^3(m-1)}\sum_{n=0}^{m-1}(3n^2+n)16^{m-1-n}\sum_{k=0}^n\bi nk^2\bi{2k}k\bi{2(n-k)}{n-k}\in\Z.$$
Moreover, $a_m$ is odd if and only if $m$ is a power of two.
\endproclaim
\Remark. The reader may consult [S-11, S-12] for more congruences and identities involving similar sums.

\proclaim{Conjecture A44 {\rm ([S-11])}} For any prime $p>5$, we have
$$\align&\sum_{n=0}^{p-1}\f{357n+103}{2160^n}\bi{2n}n\sum_{k=0}^n\bi nk\bi{n+2k}{2k}\bi{2k}k(-324)^{n-k}
\\\quad\qquad&\eq p\l(\f{-1}p\r)\l(54+49\l(\f p{15}\r)\r)\pmod{p^2},
\endalign$$
and
$$\align&\sum_{n=0}^{p-1}\f{\bi{2n}n}{2160^n}\sum_{k=0}^n\bi nk\bi{n+2k}{2k}\bi{2k}k(-324)^{n-k}
\\\eq&\cases 4x^2-2p\,(\mo\ p^2)&\t{if}\, (\f{-1}p)=(\f p3)=(\f p5)=(\f p7)=1,\, p=x^2+105y^2,
\\2x^2-2p\,(\mo\ p^2)&\t{if}\, (\f{-1}p)=(\f p7)=1,\,(\f p3)=(\f p5)=-1,\, 2p=x^2+105y^2,
\\2p-12x^2\,(\mo\ p^2)&\t{if}\,(\f{-1}p)=(\f p3)=(\f p5)=(\f p7)=-1,\, p=3x^2+35y^2,
\\2p-6x^2\,(\mo\ p^2)&\t{if}\, (\f{-1}p)=(\f p7)=-1,\, (\f p3)=(\f p5)=1,\, 2p=3x^2+35y^2,
\\20x^2-2p\,(\mo\ p^2)&\t{if}\, (\f{-1}p)=(\f p5)=1,\,(\f p3)=(\f p7)=-1,\, p=5x^2+21y^2,
\\10x^2-2p\,(\mo\ p^2)&\t{if}\, (\f{-1}p)=(\f p3)=1,\, (\f p5)=(\f p7)=-1,\, 2p=5x^2+21y^2,
\\28x^2-2p\,(\mo\ p^2)&\t{if}\, (\f{-1}p)=(\f p5)=-1,\,(\f p3)=(\f p7)=1,\, p=7x^2+15y^2,
\\14x^2-2p\,(\mo\ p^2)&\t{if}\, (\f{-1}p)=(\f p3)=-1,\, (\f p5)=(\f p7)=1,\, 2p=7x^2+15y^2,
\\0\pmod{p^2}&\t{if}\ (\f{-105}p)=-1.
\endcases
\endalign$$
Also,
$$\sum_{n=0}^\infty\f{357n+103}{2160^n}\bi{2n}n\sum_{k=0}^n\bi nk\bi{n+2k}{2k}\bi{2k}k(-324)^{n-k}=\f{90}{\pi}.$$
\endproclaim
\Remark. The quadratic field $\Q(\sqrt{-105})$ has class number eight.
I'd like to offer \$1050 for the first correct proof of all the congruences in Conj. A44, and \$90 for the first rigorous proof
of the series for $90/\pi$ in Conj. A44.

\proclaim{Conjecture A45 {\rm ([S11e])}} For any prime $p>3$, we have
$$\sum_{k=0}^{p-1}\f{\bi{6k}{3k}\bi{3k}k}{432^k}\eq\l(\f{-1}p\r)-\f{25}9p^2E_{p-3}\ (\mo\ p^3)$$
and
$$\sum_{k=0}^{p-1}\f{\bi{6k}{3k}C_k^{(2)}}{432^k}\eq\l(\f p3\r)\ (\mo\ p^2),$$
where $C_k^{(2)}$ stands for the second-order Catalan number $\bi{3k}k/(2k+1)=\bi{3k}k-2\bi{3k}{k-1}$.
\endproclaim
\Remark. A related conjecture of Rodriguez-Villegas [RV] proved by Mortenson [M2] states that if $p>3$ is a prime then
$$\sum_{k=0}^{p-1}\f{\bi{6k}{3k}\bi{3k}k}{432^k}=\sum_{k=0}^{p-1}\f{(6k)!}{k!(2k)!(3k)!}(2^43^3)^{-k}\eq\l(\f{-1}p\r)\ (\mo\ p^2).$$
Using the Gosper algorithm we find the identity
$$\sum_{k=0}^n(36k+5)\bi{6k}{3k}\bi{3k}k432^{n-k}=(6n+1)(6n+5)\bi{6n}{3n}\bi{3n}n\ \ (n\in\N)$$
which implies that
$$\sum_{k=0}^{p-1}\f{36k+5}{432^k}\bi{6k}{3k}\bi{3k}k\eq 5p^2\ (\mo\ p^3)\quad \ \t{for any prime}\ p>3.$$

\proclaim{Conjecture A46 {\rm ([S11e])}} Let $p>3$ be a prime.

{\rm (i)} We have
$$\sum_{k=0}^{p-1}\f{\bi{3k}{k,k,k}}{24^k}
\eq\cases\bi{2(p-1)/3}{(p-1)/3}\ (\mo\ p^2)&\t{if}\ p\eq1\ (\mo\ 3),
\\p/\bi{2(p+1)/3}{(p+1)/3}\ (\mo\ p^2)&\t{if}\ p\eq 2\ (\mo\ 3).\endcases$$

{\rm (ii)} When $p\eq1\pmod3$ and $4p=x^2+27y^2$ with $x\eq2\pmod 3$, we may determine $x$ mod $p^2$ in the following way:
$$\sum_{k=0}^{p-1}\f{k+2}{24^k}\bi{3k}{k,k,k}\eq x\pmod{p^2}.$$

{\rm (iii)} We also have
$$\sum_{k=0}^{p-1}\f{\bi{3k}kC_k}{24^k}\eq\f12\bi{2(p-(\f p3))/3}{(p-(\f p3))/3}\ (\mo\ p).$$
\endproclaim
\Remark. (a) I'd like to offer \$27 (27 US dollars) for the first correct proof of part (ii) of Conj. 46.

(b) It is known (cf. [HW]) that for any prime $p\eq1\pmod 3$ with $4p=x^2+27y^2$ we have
$\bi{2(p-1)/3}{(p-1)/3}\eq (\f x3)(\f px-x)\pmod{p^2}.$
The author [S-9] showed that for any prime $p>3$ we have
$$\align \sum_{k=0}^{p-1}\f{\bi{3k}{k,k,k}}{24^k}\eq&\l(\f p3\r)\sum_{k=0}^{p-1}\f{\bi{3k}{k,k,k}}{(-216)^k}\ (\mo\ p^2),
\\ \sum_{k=0}^{p-1}\f{k\bi{3k}{k,k,k}}{24^k}\eq&9\l(\f p3\r)\sum_{k=0}^{p-1}\f{k\bi{3k}{k,k,k}}{(-216)^k}\ (\mo\ p^2).
\endalign$$
Recently Z. H. Sun [Su3] confirmed  Conj. A46(i) modulo $p$.

\proclaim{Conjecture A47 {\rm ([S11e])}} Let $p>3$ be a prime. If $p\eq1\ (\mo\ 4)$ and
   $p=x^2+y^2$ with $x\eq1\ (\mo\ 4)$ and $y\eq0\ (\mo\ 2)$, then
 $$\gather\sum_{k=0}^{p-1}\f{\bi{2k}k\bi{4k}{2k}}{72^k}\eq \l(\f 6p\r)\l(2x-\f p{2x}\r)\ (\mo\ p^2),
 \\\sum_{k=0}^{p-1}\f{1-k}{72^k}\bi{2k}k\bi{4k}{2k}\eq\l(\f 6p\r)x\ (\mo\ p^2).
 \endgather$$
 If $p\eq3\pmod 4$, then
 $$\sum_{k=0}^{p-1}\f{\bi{2k}k\bi{4k}{2k}}{72^k}\eq\l(\f 6p\r)\f{2p}{3\bi{(p+1)/2}{(p+1)/4}}\pmod{p^2}.$$
 \endproclaim

 \proclaim{Conjecture A48 {\rm ([S11e])}} Let $p>3$ be a prime. Then
 $$\sum_{k=0}^{p-1}\f{\bi{2k}k\bi{4k}{2k+1}}{48^k}\eq0\ (\mo\ p^2).$$
If $p\eq1\ (\mo\ 3)$ and $p=x^2+3y^2$ with $x\eq1\ (\mo\ 3)$, then
 $$\sum_{k=0}^{p-1}\f{\bi{2k}k\bi{4k}{2k}}{48^k}\eq 2x-\f p{2x}\ (\mo\ p^2),$$
 and we may determine $x$ mod $p^2$ via the congruence
 $$\sum_{k=0}^{p-1}\f{k+1}{48^k}\bi{2k}k\bi{4k}{2k}\eq x\ (\mo\ p^2).$$
 If $p\eq2\pmod 3$, then
 $$\sum_{k=0}^{p-1}\f{\bi{2k}k\bi{4k}{2k}}{48^k}\eq\f{3p}{2\bi{(p+1)/2}{(p+1)/6}}\pmod{p^2}.$$
 \endproclaim
\Remark. The author [S-9] proved that for any prime $p>3$ we have
$$\sum_{k=0}^{p-1}\f{\bi{2k}k\bi{4k}{2k}}{(-192)^k}
\eq\l(\f{-2}p\r)\sum_{k=0}^{p-1}\f{\bi{2k}k\bi{4k}{2k}}{48^k}\ (\mo\ p^2)$$
and
$$\sum_{k=0}^{p-1}\f{k\bi{2k}k\bi{4k}{2k}}{(-192)^k}
\eq\f14\l(\f{-2}p\r)\sum_{k=0}^{p-1}\f{k\bi{2k}k\bi{4k}{2k}}{48^k}\ (\mo\ p^2).$$
I'd like to offer \$48 for the first correct proof of the first congruence in Conj. A48.

\proclaim{Conjecture A49 {\rm ([S11e])}} Let $p>3$ be a prime. If
$(\f p7)=1$ and $p=x^2+7y^2$ with $(\f x7)=1$, then
 $$\gather\sum_{k=0}^{p-1}\f{\bi{2k}k\bi{4k}{2k}}{63^k}\eq \l(\f p3\r)\l(2x-\f p{2x}\r)\ (\mo\ p^2),
 \\\sum_{k=0}^{p-1}\f{k+8}{63^k}\bi{2k}k\bi{4k}{2k}\eq 8\l(\f p3\r)x\ (\mo\ p^2).
  \endgather$$
If $(\f p7)=-1$, then
$$\sum_{k=0}^{p-1}\f{\bi{2k}k\bi{4k}{2k}}{63^k}\eq\sum_{k=0}^{p-1}\f{\bi{2k}k\bi{4k}{2k}^2}{63^k}\eq0\ (\mo\ p).$$
 \endproclaim
 \Remark. Z. H. Sun [Su2] made progress on Conjectures 47-49
 by obtaining $\sum_{k=0}^{p-1}\bi{4k}{2k}\bi{2k}k/m^k$ mod $p$
 for $m=48,\, 63,\, 72$ (see also Conj. B16).

\proclaim{Conjecture A50 {\rm ([S11e], [S-7])}} Let $p$ be an odd prime. For any $a\in\Z^+$ we have
$$\sum_{k=0}^{p^a-1}\f{\bi{2k}k\bi{4k}{2k}}{64^k}\eq\l(\f{-2}{p^a}\r)-\l(\f{-2}{p^{a-1}}\r)\f{3p^2}{16}E_{p-3}\l(\f14\r)\pmod{p^3}$$
and
$$\sum_{k=0}^{p^a-1}\f{\bi{2k}kC_{2k}}{64^k}\eq\l(\f{-1}{p^a}\r)-\l(\f{-1}{p^{a-1}}\r)3p^2E_{p-3}\ (\mo\ p^3).$$
Also,
$$\align\sum_{k=1}^{p-1}\f{\bi{2k}k\bi{4k}{2k}}{k64^k}\eq&-3H_{(p-1)/2}+\f 74p^2B_{p-3}\ (\mo\ p^3),
\\\sum_{k=1}^{(p-1)/2}\f{\bi{2k}k\bi{4k}{2k}}{k64^k}\eq&-3H_{(p-1)/2}-2\l(\f{-1}p\r)pE_{p-3}\ (\mo\ p^2),
\\p\sum_{k=1}^{(p-1)/2}\f{64^k}{k^3\bi{2k}k\bi{4k}{2k}}\eq&32\l(\f{-1}p\r)E_{p-3}\ (\mo\ p),
\\p\sum_{k=1}^{(p-1)/2}\f{64^k}{(2k-1)k^2\bi{2k}k\bi{4k}{2k}}\eq& 16\l(pE_{p-3}-\l(\f{-1}p\r)q_p(2)\r)\ (\mo\ p^2).
\endalign$$
\endproclaim
\Remark. Mortenson [M2] proved the following conjecture of
Rodriguez-Villegas [RV]: For any odd prime $p$ we have
$$\sum_{k=0}^{p-1}\f{\bi{2k}k\bi{4k}{2k}}{64^k}\eq\l(\f{-2}p\r)\
(\mo\ p^2).$$

\proclaim{Conjecture A51 {\rm ([S11b])}} Let $p$ be an odd prime.
If $p\eq1\ (\mo\ 3)$, then
 $$\sum_{k=0}^{(p-1)/2}\f{kC_k^3}{16^k}\eq2p-2\ (\mo\ p^2).$$
If $p\eq1\ (\mo\ 4)$, then
$$\sum_{k=0}^{(p-1)/2}\f{C_k^3}{64^k}\eq8\ (\mo\ p^2).$$
\endproclaim
\Remark. The author [11b] determined $\sum_{k=0}^{p-1}C_k^3/64^k$ modulo any odd prime $p$.

\proclaim{Conjecture A52 {\rm ([S11e])}} Let $p>3$ be a prime. Then for any $a\in\Z^+$ we have
$$\sum_{k=0}^{p^a-1}\f{\bi{3k}{k,k,k}}{27^k}\eq\l(\f{p^a}3\r)-\l(\f{p^{a-1}}3\r)\f{p^2}3B_{p-2}\l(\f13\r)\pmod{p^3}$$
and
$$\sum_{k=0}^{p^a-1}\f{\bi{2k}kC_k^{(2)}}{27^k}\eq\l(\f {p^a}3\r)-\f23\l(\f{p^{a-1}}3\r)p^2B_{p-2}\l(\f13\r)\ (\mo\ p^2),$$
where
$$C_k^{(2)}=\f{\bi{3k}k}{2k+1}=\bi{3k}k-2\bi{3k}{k-1}$$
is a second-order Catalan number (of the first kind).
Furthermore,
$$\sum_{k=0}^{p-1}(4k+1)\f{\bi{2k}kC_k^{(2)}}{27^k}\eq\l(\f p3\r)\ (\mo\ p^4).$$
\endproclaim
\Remark. Mortenson [M2] proved the following conjecture of Rodriguez-Villegas [RV] for primes $p>3$:
$$\sum_{k=0}^{p-1}\f{\bi{3k}{k,k,k}}{27^k}=\sum_{k=0}^{p-1}\f{\bi{2k}{k}\bi{3k}k}{27^k}\eq\l(\f p3\r)\ (\mo\ p^2).$$

\proclaim{Conjecture A53 {\rm ([S11e])}} Let $p>3$ be a prime. Then
$$\sum_{k=0}^{p-1}\f{C_kC_k^{(2)}}{27^k}\eq 2\l(\f p3\r)-p\ (\mo\ p^2)$$
and
$$\sum_{k=0}^{p-1}\f{C_k\bar C_k^{(2)}}{27^k}\eq-7\ (\mo\ p),$$
where
$$\bar C_k^{(2)}=\f 2{k+1}\bi{3k}k=2\bi{3k}k-\bi{3k}{k+1}$$
is a second-order Catalan number of the second kind.
Hence
$$\sum_{k=1}^{p-1}\f{\bi{2k}{k-1}\bi{3k}{k-1}}{27^k}\eq\l(\f p3\r)-p\ (\mo\ p^2)$$
and
$$\sum_{k=1}^{p-1}\f{\bi{2k}{k+1}\bi{3k}{k+1}}{27^k}\eq2\l(\f p3\r)-7\ (\mo\ p).$$
\endproclaim
\Remark. Note that
$$\bi{2k}{k-1}\bi{3k}{k-1}=\l(\bi{2k}k-C_k\r)\f{\bi{3k}k-C_k^{(2)}}2$$
and
$$\bi{2k}{k+1}\bi{3k}{k+1}=\l(\bi{2k}k-C_k\r)\l(2\bi{3k}k-\bar C_k^{(2)}\r).$$
{\tt Mathematica} yields that
$$\sum_{k=0}^\infty\f{C_k\bar C_k^{(2)}}{27^k}=\f{81\sqrt3}{4\pi}-9.$$

\proclaim{Conjecture A54 {\rm ([S-1])}} Let $p$ be an odd prime. Then
$$\sum_{k=1}^{p-1}\f{2^k}k\bi{3k}k\eq-3p\,q_p^2(2)\ (\mo\ p^2),$$
and
$$p\sum_{k=1}^{p-1}\f1{k2^k\bi{3k}k}\eq\cases0\ (\mo\ p^2)&\t{if}\ p\eq1\ (\mo\ 4),
\\-3/5\ (\mo\ p^2)&\t{if}\ p\eq3\ (\mo\ 4).\endcases$$
When $p>3$ we have
$$p\sum_{k=1}^{p-1}\f1{k^22^k\bi{3k}k}\eq-\f{q_p(2)}2-\f p4q_p(2)^2\ (\mo\ p^2).$$
Also,
$$\sum_{k=0}^{p-1}(25k+3)k2^k\bi{3k}k\eq 6\l(\f{-1}p\r)-18p\ (\mo\ p^2),$$
$$p\sum_{k=0}^{(p-1)/2}\f{25k-3}{2^k\bi{3k}k}\eq \l(\f{-1}p\r)-\l(\f 2p\r)\f {5p}2\ (\mo\ p^2)$$
and
$$2p\sum_{k=0}^{p-1}\f{25k-3}{2^k\bi{3k}k}\eq3\l(\f{-1}p\r)+(E_{p-3}-9)p^2\ (\mo\ p^4).$$
\endproclaim
\Remark. Gosper announced in 1974 that $\sum_{k=0}^\infty(25k-3)/(2^k\bi{3k}k)=\pi/2$.
In [ZPS] Zhao, Pan and Sun proved that
$\sum_{k=1}^{p-1}\f{2^k}k\bi{3k}k\eq0\ (\mo\ p)$ for any odd prime $p$.

\proclaim{Conjecture A55 {\rm ([S-6])}} Let
$p$ be an odd prime and let $a\in\Z^+$.

{\rm (i)} If $p^a\eq1,2\ (\mo\ 5)$, or $a>1$ and $p\not\eq3\
(\mo\ 5)$,
$$\sum_{k=0}^{\lfloor\f 45p^a\rfloor}(-1)^k\bi{2k}k\eq\l(\f5{p^a}\r)\ (\mo\ p^2).$$
If $p^a\eq1,3\ (\mo\ 5)$, or $a>1$ and $p\not\eq2\ (\mo\ 5)$, then
$$\sum_{k=0}^{\lfloor\f 35p^a\rfloor}(-1)^k\bi{2k}k\eq\l(\f5{p^a}\r)\ (\mo\ p^2).$$
Thus, if $p^a\eq1\ (\mo\ 5)$ then
$$\sum_{\f35p^a<k<\f 45p^a}(-1)^k\bi{2k}k\eq0\ (\mo\ p^2).$$

{\rm (ii)} If $p\eq1,7\ (\mo\ 10)$ or $a>2$, then
$$\sum_{k=0}^{\lfloor\f 7{10}p^a\rfloor}\f{\bi{2k}k}{(-16)^k}\eq\l(\f {p^a}5\r)\ (\mo\ p^2).$$
If $p\eq1,3\ (\mo\ 10)$ or $a>2$, then
$$\sum_{k=0}^{\lfloor\f 9{10}p^a\rfloor}\f{\bi{2k}k}{(-16)^k}\eq\l(\f {p^a}5\r)\ (\mo\ p^2).$$

{\rm (iii)} If $p\eq1\ (\mo\ 3)$ or $a>1$, then
$$\sum_{k=0}^{\lfloor\f56p^a\rfloor}\f{\bi{2k}k}{16^k}\eq\l(\f 3{p^a}\r)\ (\mo\ p^2).$$
For any nonnegative integer
$n$ we have
$$\f1{(2n+1)^2\bi{2n}n}\sum_{k=0}^n\f{\bi{2k}k}{16^k}\eq\cases1\ (\mo\ 9)&\t{if}\ 3\mid n,\\4\ (\mo\ 9)&\t{if}\ 3\nmid n.
\endcases$$. Also,
$$\f1{3^{2a}}\sum_{k=0}^{(3^a-1)/2}\f{\bi{2k}k}{16^k}\eq(-1)^a10\ (\mo\ 27)$$
for every $a=1,2,3,\ldots.$
\endproclaim
\Remark. Let $\{F_n\}_{n\gs0}$ be the Fibonacci sequence given by $F_0=0,\ F_1=1$ and $F_{n+1}=F_n+F_{n-1}\ (n=1,2,3,\ldots)$.
H. Pan and the author [PS] proved that
if $p\not=2,5$ is a prime and $a$ is a positive integer then
 $$\sum_{k=0}^{p^a-1}(-1)^k\bi{2k}k
\eq\l(\f {p^a}5\r)\l(1-2F_{p^a-(\f {p^a}5)}\r)\ (\mo\ p^3),$$
which is Conjecture 3.1 of [ST]. The author [S-6] proved that for any prime $p\not=2,5$ and $a\in\Z^+$ we have
$$\sum_{k=0}^{(p^a-1)/2}\f{\bi{2k}k}{(-16)^k}\eq\l(\f{p^a}5\r)\l(1+\f{F_{p^a-(\f{p^a}5)}}2\r)\ (\mo\ p^3).$$
He [S-6] also showed that $\sum_{k=0}^{\lfloor p^a/2\rfloor}\binom{2k}k/16^k\eq(\f 3{p^a})
\ (\mo\ p^2)$ for any odd prime $p$ and $a\in\Z^+$.

\proclaim{Conjecture A56 {\rm ([S11e])}} Let $p>3$ be a prime.
If $p\eq7\ (\mo\ 12)$ and $p=x^2+3y^2$ with $y\eq1\ (\mo\ 4)$, then
$$\sum_{k=0}^{p-1}\l(\f k3\r)\f{\bi{2k}k^2}{(-16)^k}\eq (-1)^{(p-3)/4}\l(4y-\f p{3y}\r)\ (\mo\ p^2)$$
and
$$\sum_{k=0}^{p-1}\l(\f k3\r)\f{k\bi{2k}k^2}{(-16)^k}\eq (-1)^{(p+1)/4}y\ (\mo\ p^2).$$
If $p\eq1\ (\mo\ 12)$, then
$$\sum_{k=0}^{p-1}\bi{p-1}k\l(\f k3\r)\f{\bi{2k}k^2}{16^k}\eq0\ (\mo\ p^2).$$
\endproclaim

Recall that the Pell sequence $\{P_n\}_{n\gs0}$ and its companion $\{Q_n\}_{n\gs0}$ are defined as follows:
$$P_0=0,\ P_1=1,\ \t{and}\ P_{n+1}=2P_n+P_{n-1}\ (n=1,2,3,\ldots);$$
$$Q_0=2,\ Q_1=2,\ \t{and}\ Q_{n+1}=2Q_n+Q_{n-1}\ (n=1,2,3,\ldots).$$

\proclaim{Conjecture A57 {\rm ([S-4])}} {\rm (i)} Let $p$ be a prime
with $p\eq1,3\ (\mo\ 8)$. Write $p=x^2+2y^2$ with $x,y\in\Z$ and $x\eq1,3\ (\mo\ 8)$. Then
$$\sum_{k=0}^{p-1}\f{P_k}{(-8)^k}\bi{2k}k^2\eq \cases0\ (\mo\ p^2)&\t{if}\ p\eq1\ (\mo\ 8),
\\(-1)^{(p-3)/8}(p/(2x)-2x)\ (\mo\ p^2)&\t{if}\ p\eq3\ (\mo\ 8).
\endcases$$
Also,
$$\sum_{k=0}^{p-1}\f{kP_k}{(-8)^k}\bi{2k}k^2\eq\f{(-1)^{(x+1)/2}}2\l(x+\f p{2x}\r)\ (\mo\ p^2).$$

{\rm (ii)} If $p\eq7\ (\mo\ 8)$ is a prime, then
$$\sum_{k=0}^{p-1}\bi{p-1}k\f{P_k}{8^k}\bi{2k}k^2\eq0\ (\mo\ p^2).$$
\endproclaim

\proclaim{Conjecture A58 {\rm ([S-4])}} Let $p$ be an odd prime.

{\rm (i)} If $p\eq3\ (\mo\ 8)$ and $p=x^2+2y^2$ with $y\eq1,3\ (\mo\ p)$, then
$$\sum_{k=0}^{p-1}\f{P_k}{32^k}\bi{2k}k^2\eq (-1)^{(y-1)/2}\l(2y-\f p{4y}\r)\ (\mo\ p^2).$$

{\rm (ii)} Suppose that $p\eq1,3\ (\mo\ 8)$, $p=x^2+2y^2$ with $x\eq1,3\ (\mo\ 8)$
and also $y\eq1,3\ (\mo\ 8)$ when $p\eq3\ (\mo\ 8)$.
Then
$$\sum_{k=0}^{p-1}\f{kP_k}{32^k}\bi{2k}k^2\eq\cases(-1)^{(p-1)/8}(p/(4x)-x/2)\ (\mo\ p^2)&\t{if}\ p\eq1\ (\mo\ 8),
\\(-1)^{(y+1)/2}y\ (\mo\ p^2)&\t{if}\ p\eq3\ (\mo\ 8).
\endcases$$
\endproclaim

\proclaim{Conjecture A59 {\rm ([S-4])}} Let $p$ be a prime with
$p\eq1,3\ (\mo\ 8)$. Write $p=x^2+2y^2$ with $x,y\in\Z$ and $x\eq1,3\ (\mo\ 8)$. Then we have
$$\sum_{k=0}^{p-1}\f{Q_k}{(-8)^k}\bi{2k}k^2\eq (-1)^{(x-1)/2}\l(4x-\f px\r)\ (\mo\ p^2)$$
and
$$\sum_{k=0}^{p-1}\f{kQ_k}{(-8)^k}\bi{2k}k^2\eq\cases0\ (\mo\ p^2)&\t{if}\ p\eq1\ (\mo\ 8),
\\(-1)^{(p-3)/8}2(x+p/x)\ (\mo\ p^2)&\t{if}\ p\eq3\ (\mo\ 8).
\endcases$$
\endproclaim

\proclaim{Conjecture A60 {\rm ([S-4])}} Let $p$ be an odd prime.

{\rm (i)} When $p\eq1\ (\mo\ 8)$ and $p=x^2+2y^2$ with $x,y\in\Z$ and $x\eq1,3\ (\mo\ 8)$, we have
$$\sum_{k=0}^{p-1}\f{Q_k}{32^k}\bi{2k}k^2\eq (-1)^{(p-1)/8}\l(4x-\f px\r)\ (\mo\ p^2).$$

{\rm (ii)} If $p\eq1,3\ (\mo\ 8)$ and $p=x^2+2y^2$ with $x\eq1,3\ (\mo\ 8)$
and also $y\eq1,3\ (\mo\ 8)$ when $p\eq3\ (\mo\ 8)$,
then
$$\sum_{k=0}^{p-1}\f{k\bi{2k}k^2}{32^k}Q_k\eq\cases(-1)^{(p-1)/8}(p/x-2x)\ (\mo\ p^2)&\t{if}\ p\eq1\ (\mo\ 8),
\\(-1)^{(y+1)/2}2y\ (\mo\ p^2)&\t{if}\ p\eq 3\ (\mo\ 8).\endcases$$
\endproclaim

\proclaim{Conjecture A61 {\rm ([S-4])}} Let $p>3$ be a prime.
If $p\eq7\ (\mo\ 12)$ and $p=x^2+3y^2$ with $y\eq1\ (\mo\ 4)$, then
$$\sum_{k=0}^{p-1}\f{u_k(4,1)}{4^k}\bi{2k}k^2\eq (-1)^{(p+1)/4}\l(4y-\f p{3y}\r)\ (\mo\ p^2)$$
and
$$\sum_{k=0}^{p-1}\f{ku_k(4,1)}{4^k}\bi{2k}k^2\eq(-1)^{(p-3)/4}\l(6y-\f{7p}{3y}\r)\ (\mo\ p^2).$$
If $p\eq1\ (\mo\ 12)$, then
$$\sum_{k=0}^{p-1}\f{u_k(4,1)}{4^k}\bi{2k}k^2\eq0\ (\mo\ p^2)$$
\endproclaim

\proclaim{Conjecture A62 {\rm ([S-4])}} Let $p>3$ be a prime.
If $p\eq7\ (\mo\ 12)$ and $p=x^2+3y^2$ with $y\eq1\ (\mo\ 4)$, then
$$\sum_{k=0}^{p-1}\f{u_k(4,1)}{64^k}\bi{2k}k^2\eq 2y-\f p{6y}\ (\mo\ p^2)$$
and
$$\sum_{k=0}^{p-1}\f{ku_k(4,1)}{64^k}\bi{2k}k^2\eq y\ (\mo\ p^2).$$
\endproclaim

\proclaim{Conjecture A63 {\rm ([S-4])}} Let $p$ be an odd prime.

{\rm (i)} If $p\eq1\ (\mo\ 12)$ and $p=x^2+3y^2$ with $x\eq1\ (\mo\ 3)$, then
$$\sum_{k=0}^{p-1}\f{v_k(4,1)}{4^k}\bi{2k}k^2\eq(-1)^{(p-1)/4+(x-1)/2}\l(4x-\f p{x}\r)\ (\mo\ p^2)$$
and
$$\sum_{k=0}^{p-1}\f{v_k(4,1)}{64^k}\bi{2k}k^2\eq(-1)^{(x-1)/2}\l(4x-\f p{x}\r)\ (\mo\ p^2);$$
also
$$\sum_{k=0}^{p-1}\f{kv_k(4,1)}{4^k}\bi{2k}k^2\eq (-1)^{(p-1)/4+(x+1)/2}\l(4x-\f{2p}x\r)\ (\mo\ p^2)$$
and
$$\sum_{k=0}^{p-1}\f{kv_k(4,1)}{64^k}\bi{2k}k^2\eq (-1)^{(x-1)/2}\l(2x-\f{p}x\r)\ (\mo\ p^2).$$

{\rm (ii)} If $p\eq7\ (\mo\ 12)$ and $p=x^2+3y^2$ with $y\eq1\ (\mo\ 4)$, then
$$\sum_{k=0}^{p-1}\f{v_k(4,1)}{4^k}\bi{2k}k^2\eq (-1)^{(p-3)/4}\l(12y-\f p{y}\r)\ (\mo\ p^2),$$
$$\sum_{k=0}^{p-1}\f{kv_k(4,1)}{4^k}\bi{2k}k^2\eq(-1)^{(p+1)/4}\l(20y-\f{8p}y\r)\ (\mo\ p^2)$$
and
$$\sum_{k=0}^{p-1}\f{kv_k(4,1)}{64^k}\bi{2k}k^2\eq 4y\ (\mo\ p^2).$$

{\rm (iii)} If $p\eq 11\ (\mo\ 12)$, then
$$\sum_{k=0}^{p-1}\bi{p-1}k\f{v_k(4,1)}{(-4)^k}\bi{2k}k^2\eq0\ (\mo\ p^2).$$
\endproclaim
\Remark. Recently Z. H. Sun [Su2] confirmed the first congruences in Conjectures 56-62 modulo $p$,
and the first and the second congruences in Conj. A63(i) and the first congruence in Conj. A63(ii)
modulo $p$. See also Conj. B12 for his progress on some of the author's conjectured congruences.

\medskip

Recall Ap\'ery numbers are those integers
$$A_n=\sum_{k=0}^n\bi nk^2\bi{n+k}k^2\ (n\in\N)$$
which play a central role in Ap\'ery's proof of the irrationality of $\zeta(3)=\sum_{n=1}^\infty1/n^3$.
We also define Ap\'ery polynomials by
$$A_n(x)=\sum_{k=0}^n\bi nk^2\bi{n+k}k^2x^k\ \ (n=0,1,2,\ldots).$$
Note that $A_n(1)=A_n$.

\proclaim{Conjecture A64}
{\rm (i) ([S11j])} Let $p$ be an odd prime. Then
$$\align&\sum_{k=0}^{p-1}A_k
\\\eq&\cases4x^2-2p\ (\mo\ p^2)&\t{if}\ p\eq1,3\ (\mo\ 8)\ \t{and}\ p=x^2+2y^2\ (x,y\in\Z),
\\0\ (\mo\ p^2)&\t{if}\ p\eq5,7\ (\mo\ 8);\endcases
\endalign$$
and
$$\align&\sum_{k=0}^{p-1}(-1)^kA_k
\\\eq&\cases4x^2-2p\ (\mo\ p^2)&\t{if}\ p\eq1\ (\mo\ 3)\ \t{and}\ p=x^2+3y^2\ (x,y\in\Z),
\\0\ (\mo\ p^2)&\t{if}\ p\eq2\ (\mo\ 3).\endcases
\endalign$$

{\rm (ii) ([S11j])} Let $p>3$ be a prime. If $p\eq1\pmod3$, then
$$\sum_{k=0}^{p-1}(-1)^kA_k\eq\sum_{k=0}^{p-1}\f{\bi{2k}k^3}{16^k}\pmod{p^3}.$$
If $p\eq1,3\pmod8$, then
$$\sum_{k=0}^{p-1}A_k\eq\sum_{k=0}^{p-1}\f{\bi{4k}{k,k,k,k}}{256^k}\pmod{p^3}.$$

{\rm (iii) ([S11j])} Let $p>3$ be a prime. If $x$ belongs to the set
$$\align &\{1,-4,9,-48,81,-324,2401,9801,-25920,-777924,96059601\}
\\&\quad\qquad\bigcup\l\{\f{81}{256},-\f9{16},\f{81}{32},-\f{3969}{256}\r\}
\endalign$$
and $x\not\eq0\pmod p$, then we must have
$$\sum_{k=0}^{p-1}A_k(x)\eq\l(\f xp\r)\sum_{k=0}^{p-1}\f{\bi{4k}{k,k,k,k}}{(256x)^k}\pmod{p^2}.$$
\endproclaim
\Remark. Let $p$ be an odd prime.
In [S11j] the author proved that
$$\sum_{k=0}^{p-1}(-1)^kA_k(x)\eq\sum_{k=0}^{p-1}\f{\bi{2k}k^3}{16^k}x^k\pmod{p^2},$$
and that
$$\sum_{k=0}^{p-1}A_k(x)\eq\l(\f xp\r)\sum_{k=0}^{p-1}\f{\bi{4k}{k,k,k,k}}{(256x)^k}\pmod{p}$$
for any $p$-adic integer $x\not\eq0\pmod p$. Thus the two congruences in Conj. A64(i) hold modulo $p$
and also $\sum_{k=0}^{p-1}(-1)^kA_k\eq0\pmod{p^2}$ if $p\eq2\ (\mo\ 3)$.
For those $x=-4,9,-48,81,-324,2401,9801,-25920,-777924$, 96059601, 81/256, the author
had conjectures on $\sum_{k=0}^{p-1}\bi{4k}{k,k,k,k}/(256x)^k$ mod $p^2$ (see A3, A13-19, A21, A24, A25, A28).
Motivated by this, Z. H. Sun [Su1] guessed
$\sum_{k=0}^{p-1}\bi{4k}{k,k,k,k}/(256x)^k$ mod $p^2$ for $x=-\f9{16},\f{81}{32},-\f{3969}{256}$.

\medskip

\proclaim{Conjecture A65} {\rm (i) ([S11j])} For any prime $p>3$ we have
$$\sum_{k=0}^{p-1}(2k+1)A_k\eq p-\f 72p^2H_{p-1}\ (\mo\ p^6),$$
$$\sum_{k=0}^{p-1}(2k+1)^3A_k\eq p^3+4p^4H_{p-1}+\f 65p^8B_{p-5}\ (\mo\ p^9)$$
and
$$\sum_{k=0}^{p-1}(2k+1)^3(-1)^kA_k\eq-\f p3\l(\f p3\r)\ (\mo\ p^3).$$

{\rm (ii) (Discovered on Oct. 2, 2011)} Let $p>5$ be a prime. If $(\f p{15})=-1$, then
$$\gather \sum_{k=0}^{p-1}A_ku_k(7,1)\eq\sum_{k=0}^{p-1}A_kv_k(7,1)\eq0\ (\mo\ p^2),
\\\sum_{k=0}^{p-1}kA_ku_k(7,1)\eq\f p{90}\l(25\l(\f p3\r)+27\r)\pmod{p^2},
\\\sum_{k=0}^{p-1}kA_kv_k(7,1)\eq-\f p2\l(5\l(\f p3\r)+3\r)\pmod{p^2}.
\endgather$$
If $p\eq1,4\ (\mo\ 15)$ and $p=x^2+15y^2\ (x,y\in\Z)$, then
$$\gather\sum_{k=0}^{p-1}A_ku_k(7,1)\eq0\pmod{p^3},
\\\sum_{k=0}^{p-1}kA_ku_k(7,1)\eq\f{3p-4x^2}{45}\pmod{p^2},
\\\sum_{k=0}^{p-1}A_kv_k(7,1)\eq 8x^2-2p\pmod{p^2},
\\\sum_{k=0}^{p-1}(2k+1)A_kv_k(7,1)\eq2p\pmod{p^3}.
\endgather$$
If $p\eq2,8\ (\mo\ 15)$ and $p=3x^2+5y^2\ (x,y\in\Z)$, then
$$\gather\sum_{k=0}^{p-1}A_ku_k(7,1)\eq2p-12x^2\pmod{p^2},
\\\sum_{k=0}^{p-1}(45k+19)A_ku_k(7,1)\eq26p\pmod{p^3},
\\\sum_{k=0}^{p-1}A_kv_k(7,1)\eq 84x^2-14p\pmod{p^2},
\\\sum_{k=0}^{p-1}(7k+3)A_kv_k(7,1)\eq-28p\pmod{p^3}.
\endgather$$

{\rm (iii) (Discovered on Oct. 2, 2011)} Let $p>5$ be a prime. If $p\eq3\ (\mo\ 4)$, then
$$\gather \sum_{k=0}^{p-1}(-1)^kA_ku_k(14,1)\eq\sum_{k=0}^{p-1}(-1)^kA_kv_k(14,1)\eq0\ (\mo\ p^2),
\\\sum_{k=0}^{p-1}(-1)^kkA_ku_k(14,1)\eq-\f p{48}\l(15\l(\f p3\r)+16\r)\pmod{p^2},
\\\sum_{k=0}^{p-1}(-1)^kkA_kv_k(14,1)\eq p\l(5\l(\f p3\r)+4\r)\pmod{p^2}.
\endgather$$
If $p\eq1\ (\mo\ 12)$ and $p=x^2+9y^2\ (x,y\in\Z)$, then
$$\gather\sum_{k=0}^{p-1}(-1)^kA_ku_k(14,1)\eq0\pmod{p^3},
\\\sum_{k=0}^{p-1}(-1)^kkA_ku_k(14,1)\eq\f{3p-4x^2}{48}\pmod{p^2},
\\\sum_{k=0}^{p-1}(-1)^kA_kv_k(14,1)\eq 8x^2-4p\pmod{p^2},
\\\sum_{k=0}^{p-1}(-1)^k(2k+1)A_kv_k(14,1)\eq2p\pmod{p^3}.
\endgather$$
If $p\eq5\ (\mo\ 12)$ and $p=x^2+y^2\ (x,y\in\Z)$, then
$$\gather\sum_{k=0}^{p-1}(-1)^kA_ku_k(14,1)\eq-4xy\l(\f{xy}3\r)\pmod{p^2},
\\\sum_{k=0}^{p-1}(-1)^k(48k+17)A_ku_k(14,1)\eq31p\pmod{p^3},
\\\sum_{k=0}^{p-1}(-1)^kA_kv_k(14,1)\eq 56xy\l(\f{xy}3\r)\pmod{p^2},
\\\sum_{k=0}^{p-1}(-1)^k(14k+5)A_kv_k(14,1)\eq-126p\pmod{p^3}.
\endgather$$
\endproclaim

\Remark. The author [S11j] proved that $n\mid \sum_{k=0}^{n-1}(2k+1)A_k(x)$ for all $n\in\Z^+$ and $x\in\Z$
and that $\sum_{k=0}^{p-1}(2k+1)A_k\eq p+\f 76p^4B_{p-3}\ (\mo\ p^5)$ for any prime $p>3$.
Motivated by the author's work in [S11j], Guo and Zeng [GZ1] proved that
$n^3\mid \sum_{k=0}^{n-1}(2k+1)^3A_k$ for all $n\in\Z^+$ and  $\sum_{k=0}^{p-1}(2k+1)^3A_k\eq p^3\ (\mo\ p^6)$
for any prime $p>3$.
\medskip

Those integers
$$D_n=\sum_{k=0}^n\bi nk\bi{n+k}k=\sum_{k=0}^n\bi{n+k}{2k}\bi{2k}k\ (n\in\N)$$
are called central Delannoy numbers; they arise naturally in many enumeration problems in combinatorics.

\proclaim{Conjecture A66 {\rm ([S11f])}} Let $p>3$ be a prime. Then
$$\sum_{k=1}^{p-1}\f{D_k}k\eq-q_p(2)+p\,q_p(2)^2\ (\mo\ p^2).$$
Also, $$\sum_{n=1}^{p-1}D_nS_n\eq-2pH_{(p-1)/2}+\f 43p^3B_{p-3}\eq-2p\sum_{k=1}^{p-1}\f{3+(-1)^k}k\ (\mo\ p^4)$$
and
$$\sum_{n=1}^{(p-1)/2}D_nS_n\eq\cases 4x^2\ (\mo\ p)&\t{if}\ p\eq1\ (\mo\ 4)\ \&\ p=x^2+y^2\ (2\nmid x),
\\0\ (\mo\ p)&\t{if}\ p\eq3\ (\mo\ 4),\endcases$$
where
$$S_n=\sum_{k=0}^n\bi {n+k}{2k}C_k=\sum_{k=0}^n\f1{k+1}\bi nk\bi{n+k}k$$
is the $n$th Schr\"oder number.
\endproclaim
\Remark.  For any prime $p>3$, the author ([S11f], [S11j]) proved that
$$\sum_{k=0}^{p-1}D_k\eq\l(\f{-1}p\r)-p^2E_{p-3}\ (\mo\ p^3),
\ \sum_{k=1}^{p-1}\f{D_k}{k^2}\eq2\l(\f{-1}p\r)E_{p-3}\ (\mo\ p),$$
and that
$$\sum_{k=1}^{p-1}\f{D_k}k\eq -q_p(2)\ (\mo\ p),\ \ \sum_{k=1}^{p-1}\f{D_k^2}{k^2}\eq-2q_p(2)^2\pmod p.$$
\medskip

Just like $A_n(x)=\sum_{k=0}^n\bi nk^2\bi{n+k}k^2x^k$ we define
$$D_n(x)=\sum_{k=0}^n\bi nk\bi{n+k}kx^k.$$
Actually $D_n((x-1)/2)$ coincides with the Legendre polynomial $P_n(x)$ of degree $n$.

\proclaim{Conjecture A67 {\rm ([S11j])}} {\rm (i)} For any $n\in\Z$ the numbers
$$s(n)=\f1{n^2}\sum_{k=0}^{n-1}(2k+1)(-1)^kA_k\l(\f14\r)$$
and
$$t(n)=\f1{n^2}\sum_{k=0}^{n-1}(2k+1)(-1)^kD_k\l(-\f14\r)^3$$
are rational numbers with denominators $2^{2\nu_2(n!)}$ and $2^{3(n-1+\nu_2(n!))-\nu_2(n)}$ respectively.
Moreover, the numerators of $s(1),s(3),s(5),\ldots$ are congruent to $1$ modulo $12$ and the numerators
of $s(2),s(4),s(6),\ldots$ are congruent to $7$ modulo $12$.
If $p$ is an odd prime and $a\in\Z^+$, then
$$s(p^a)\eq t(p^a)\eq1\ (\mo\ p).$$
For $p=3$ and $a\in\Z^+$ we have
$$s(3^a)\eq4\ (\mo\ 3^2)\ \ \ \t{and}\ \ \ t(3^a)\eq-8\ (\mo\ 3^5).$$

{\rm (ii)} Let $p$ be a prime. For any positive integer $n$ and $p$-adic integer $x$, we have
$$\nu_p\(\f1n\sum_{k=0}^{n-1}(2k+1)(-1)^kA_k\l(x\r)\)\gs\min\{\nu_p(n),\nu_p(4x-1)\}$$
and
$$\nu_p\(\f1n\sum_{k=0}^{n-1}(2k+1)(-1)^kD_k\l(x\r)^3\)\gs\min\{\nu_p(n),\nu_p(4x+1)\}.$$
\endproclaim

\proclaim{Conjecture A68 {\rm ([S11k])}}
 If $p$ is an odd prime and $x\not\eq0,-1\pmod p$ is an integer, then
$$\sum_{k=0}^{p-1}(2k+1)(-1)^kD_k(x)^3\eq p\l(\f{4x+1}p\r)\ (\mo\ p^2),$$
and $$\sum_{k=0}^{p-1}(2k+1)D_k(x)^4\eq p\ (\mo\ p^2).$$
\endproclaim

\proclaim{Conjecture A69} Let $p>3$ be a prime.

{\rm (i)} {\rm ([S11f])} We have
$$\align&\sum_{k=0}^{p-1}D_k(-3)^3=\sum_{k=0}^{p-1}(-1)^kD_k(2)^3
\\\eq&\sum_{k=0}^{p-1}(-1)^k D_k\l(-\f14\r)^3
\eq\l(\f{-2}p\r)\sum_{k=0}^{p-1}(-1)^kD_k\l(\f18\r)^3
\\\eq&\cases(\f{-1}p)(4x^2-2p)\ (\mo\ p^2)&\t{if}\ p\eq1\ (\mo\ 3)\ \&\ p=x^2+3y^2\, (x,y\in\Z),
\\0\ (\mo\ p^2)&\t{if}\ p\eq2\ (\mo\ 3).\endcases
\endalign$$
Also,
$$\align&\l(\f{-1}p\r)\sum_{k=0}^{p-1}(-1)^kD_k\l(\f12\r)^3
\\\eq&\cases4x^2-2p\ (\mo\ p^2)&\t{if}\ p\eq1,7\ (\mo\ 24)\ \t{and}\ p=x^2+6y^2,
\\8x^2-2p\ (\mo\ p^2)&\t{if}\ p\eq5,11\ (\mo\ 24)\ \t{and}\ p=2x^2+3y^2,
\\0\ (\mo\ p^2)&\t{if}\ (\f {-6}p)=-1.\endcases
\endalign$$
And
$$\align&\sum_{k=0}^{p-1}D_k(3)^3=\sum_{k=0}^{p-1}(-1)^kD_k(-4)^3\eq\l(\f{-5}p\r)\sum_{k=0}^{p-1}(-1)^kD_k\l(-\f1{16}\r)^3
\\\eq&\cases4x^2-2p\ (\mo\ p^2)&\t{if}\ p\eq1,4\ (\mo\ 15)\ \t{and}\ p=x^2+15y^2,
\\12x^2-2p\ (\mo\ p^2)&\t{if}\ p\eq2,8\ (\mo\ 15)\ \t{and}\ p=3x^2+5y^2,
\\0\ (\mo\ p^2)&\t{if}\ (\f p{15})=-1.\endcases
\endalign$$

{\rm (ii) (Discovered on Sept. 30, 2011)} If $(\f{-6}p)=-1$, then
$$\sum_{k=0}^{p-1}D_k^3u_k(6,1)\eq\sum_{k=0}^{p-1}D_k^3v_k(6,1)\eq0\pmod{p^2}.$$
If $p\eq1,7\ (\mo\ 24)$ and $p=x^2+6y^2\ (x,y\in\Z)$, then
$$\gather\sum_{k=0}^{p-1}D_k^3u_k(6,1)\eq0\pmod{p^2},
\\\sum_{k=0}^{p-1}kD_k^3u_k(6,1)\eq-\f{11}{96}x^2\pmod{p^2},
\\\sum_{k=0}^{p-1}D_k^3v_k(6,1)\eq8x^2-4p\pmod{p^2},
\\\sum_{k=0}^{p-1}(2k+1)D_k^3v_k(6,1)\eq-\f p4\pmod{p^2}.
\endgather$$
If $p\eq5,11\ (\mo\ 24)$ and $p=2x^2+3y^2\ (x,y\in\Z)$, then
$$\gather\sum_{k=0}^{p-1}D_k^3u_k(6,1)\eq8x^2-2p\pmod{p^2},
\\\sum_{k=0}^{p-1}(128k+53)D_k^3u_k(6,1)\eq30p\pmod{p^3},
\\\sum_{k=0}^{p-1}D_k^3v_k(6,1)\eq12p-48x^2\pmod{p^2},
\\\sum_{k=0}^{p-1}(144k+61)D_k^3v_k(6,1)\eq-186p\pmod{p^2}.
\endgather$$
\endproclaim
\Remark. It is known that $(-1)^nD_n(x)=D_n(-x-1)$ (see [S11f, Remark 1.2]).

\medskip

Recall that
$$T_n:=[x^n](1+x+x^2)^n=\sum_{k=0}^{\lfloor n/2\rfloor}\bi n{2k}\bi{2k}k$$
is called a central trinomial coefficient.
And  those numbers
$$M_n=\sum_{k=0}^{\lfloor n/2\rfloor}\bi n{2k}C_k\ \ (n=0,1,2,\ldots)$$
are called Motzkin numbers. H. Q. Cao and the author [CS] showed that
$$T_{p-1}\eq\l(\f p3\r)3^{p-1}\ (\mo\ p^2)\qquad\t{for any prime}\ p>3.$$

\proclaim{Conjecture A70 {\rm ([S11k])}}
For any prime $p>3$, we have
$$\align\sum_{k=0}^{p-1}M_k^2&\eq(2-6p)\l(\f p3\r)\ (\mo\ p^2),
\\ \sum_{k=0}^{p-1}kM_k^2&\eq(9p-1)\l(\f p3\r)\ (\mo\ p^2),
\\\sum_{k=0}^{p-1}T_kM_k\eq&\f 43\l(\f p3\r)+\f p6\l(1-9\l(\f p3\r)\r)\ (\mo\ p^2),
\\\sum_{k=0}^{p-1}\f{T_kH_k}{3^k}\eq& \f{3+(\f p3)}2-p\l(1+\l(\f p3\r)\r)\ (\mo\ p^2).
\endalign$$
\endproclaim
\Remark. The author [S11k] proved that $\sum_{k=0}^{p-1}T_k^2\eq(\f{-1}p)\ (\mo\ p)$
for any odd prime $p$.
\medskip

For $b,c\in\Z$ and $n\in\N$ we define
$$T_n(b,c):=[x^n](x^2+bx+c)^n=\sum_{k=0}^{\lfloor n/2\rfloor}\bi n{2k}\bi {2k}kb^{n-2k}c^k$$
and
$$M_n(b,c):=\sum_{k=0}^{\lfloor n/2\rfloor}\bi n{2k}C_kb^{n-2k}c^k.$$

\proclaim{Conjecture A71}
{\rm (i) ([S11k])} Let $b$ and $c$ be integers. For any $n\in\Z^+$ we have
$$\sum_{k=0}^{n-1}T_k(b,c)M_k(b,c)(b^2-4c)^{n-1-k}\eq0\ (\mo\ n).$$
If $p$ is an odd prime not dividing $c(b^2-4c)$, then
$$\sum_{k=0}^{p-1}\f{T_k(b,c)M_k(b,c)}{(b^2-4c)^k}\eq\f{pb^2}{2c}\(\l(\f{b^2-4c}p\r)-1\)\ (\mo\ p^2).$$

{\rm (ii) ([S11k])} Let $p>3$ be a prime. Then
$$\sum_{k=0}^{p-1}\f{T_k(3,3)M_k(3,3)}{(-3)^k}\eq\cases2p^2\ (\mo\ p^3)&\t{if}\ p\eq1\ (\mo\ 3),
\\p^3-p^2-3p\ (\mo\ p^4)&\t{if}\ p\eq2\ (\mo\ 3).\endcases$$
\endproclaim

\Remark. The author [S11k] proved that if $p$ is an odd prime not dividing $c(b^2-4c)$ then
$$\sum_{k=0}^{p-1}\f{T_k(b,c)M_k(b,c)}{(b^2-4c)^k}\eq0\ (\mo\ p)$$
and $$\sum_{k=0}^{P-1}\f{T_k(b,c^2)H_k}{(b-2c)^k}\eq1+\l(\f{b^2-4c}p\r)+\f b{2c}\l(\l(\f{b^2-4c}p\r)-1\r)\ (\mo\ p).$$
He also showed that if $p$ is an odd prime not dividing $c(b^2-4c^2)$ then
$$\sum_{k=0}^{p-1}\f{T_k(b,c^2)M_k(b,c^2)}{(b-2c)^{2k}}\eq\f{4b}{b+2c}\l(\f {b^2-4c^2}p\r)\ (\mo\ p).$$

\proclaim{Conjecture A72 {\rm ([S11k])}} Let $b,c\in\Z$. For any $n\in\Z^+$ we have
$$\sum_{k=0}^{n-1}(8ck+4c+b)T_k(b,c^2)^2(b-2c)^{2(n-1-k)}\eq0\ (\mo\ n).$$
If $p$ is an odd prime not dividing $b(b-2c)$, then
$$\sum_{k=0}^{p-1}(8ck+4c+b)\f{T_k(b,c^2)^2}{(b-2c)^{2k}}\eq p(b+2c)\l(\f{b^2-4c^2}p\r)\ (\mo\ p^2).$$
\endproclaim
\Remark. The author [S11k] showed that $$\f1n\sum_{k=0}^{n-1}(2k+1)T_k\ 3^{n-1-k}=\sum_{k=0}^{n-1}\bi{n-1}k(-1)^{n-1-k}(k+1)\bi{2k}k$$
for all $n=1,2,3,\ldots$ and that if $p>3$ is a prime then
$$\sum_{k=0}^{p-1}(2k+1)\f{T_k}{3^k}\eq\f p3\l(\f p3\r)+\f{p^2}3\l(1+\l(\f p3\r)\r)\ (\mo\ p^3).$$
He also proved that for any $b,c\in\Z$ and odd prime $p\nmid b-2c$ we have
$$\sum_{k=0}^{p-1}\f{T_k(b,c^2)^2}{(b-2c)^{2k}}\eq\l(\f{-c^2}p\r)\ (\mo\ p).$$
The author (cf. [S]) added the sequence
$$\f1n\sum_{k=0}^{n-1}(8k+5)T_k^2\ \ \ (n=1,2,3,\ldots)$$
as A179100 at Sloane's OEIS.
\medskip

\proclaim{Conjecture A73 {\rm ([S11k])}} Let $p$ be an odd prime.
We have
$$\sum_{k=0}^{p-1}\f{T_k(2,2)^2}{4^k}-\sum_{k=0}^{p-1}\f{\bi{2k}k^2}{8^k}
\eq\cases0\ (\mo\ p^3)&\t{if}\ p\eq1\ (\mo\ 4),\\0\ (\mo\ p^2)&\t{if}\ p\eq3\ (\mo\ 4).\endcases$$
If $p>3$, then
$$\sum_{k=0}^{p-1}\f{T_k(4,1)^2}{4^k}\eq\sum_{k=0}^{p-1}\f{T_k(4,1)^2}{36^k}\eq\l(\f{-1}p\r)\ (\mo\ p^2).$$
\endproclaim
\Remark. The author [S11k] proved that for any prime $p>3$ we have
$$\sum_{k=0}^{p-1}\f{T_k(6,-3)^2}{48^k}\eq\l(\f{-1}p\r)+\f{p^2}3E_{p-3}\ (\mo\ p^3),$$
$$\sum_{k=0}^{p-1}\f{T_k(2,-1)^2}{8^k}\eq\l(\f{-2}p\r)\ (\mo\ p^2),
\ \sum_{k=0}^{p-1}\f{T_k(2,-3)^2}{16^k}\eq\l(\f{p}3\r)\ (\mo\ p^2).$$

\proclaim{Conjecture A74 {\rm ([S11k])}} Let $p>3$ be a prime.

{\rm (i)} We have
$$\align&\l(\f 3p\r)\sum_{k=0}^{p-1}\f{T_k(2,3)^3}{8^k}\eq\sum_{k=0}^{p-1}\f{T_k(2,3)^3}{(-64)^k}
\\\eq&\sum_{k=0}^{p-1}\f{T_k(2,9)^3}{(-64)^k}\eq\l(\f 3p\r)\sum_{k=0}^{p-1}\f{T_k(2,9)^3}{512^k}
\\\eq&\cases4x^2-2p\ (\mo\ p^2)&\t{if}\ p\eq1,7\ (\mo\ 24)\ \t{and}\ p=x^2+6y^2,
\\2p-8x^2\ (\mo\ p^2)&\t{if}\ p\eq5,11\ (\mo\ 24)\ \t{and}\ p=2x^2+3y^2,
\\0\ (\mo\ p^2)&\t{if}\ (\f {-6}p)=-1.\endcases
\endalign$$
And
$$\align \sum_{k=0}^{p-1}(3k+2)\f{T_k(2,3)^3}{8^k}\eq&p\l(3\l(\f 3p\r)-1\r)\pmod{p^2},
\\\sum_{k=0}^{p-1}(3k+1)\f{T_k(2,3)^3}{(-64)^k}\eq&p\l(\f{-2}p\r)\pmod{p^3}.
\endalign$$
When $(\f{-6}p)=1$ we have
$$\sum_{k=0}^{p-1}(72k+47)\f{T_k(2,9)^3}{(-64)^k}\eq 42p\ (\mo\ p^2)$$
and
$$\sum_{k=0}^{p-1}(72k+25)\f{T_k(2,9)^3}{512^k}\eq 12p\l(\f 3p\r)\ (\mo\ p^2).$$
Also,
$$\sum_{k=0}^{n-1}(3k+2)T_k(2,3)^38^{n-1-k}\eq0\ (\mo\ 2n)$$
and
$$\sum_{k=0}^{n-1}(3k+1)T_k(2,3)^3(-64)^{n-1-k}\eq0\ (\mo\ n)$$
for every positive integer $n$.

{\rm (ii)} We have
$$\align&\l(\f {2}p\r)\sum_{k=0}^{p-1}\f{T_k(18,49)^3}{8^{3k}}\eq\sum_{k=0}^{p-1}\f{T_k(18,49)^3}{16^{3k}}
\\\eq&\cases 4x^2-2p\pmod{p^2}&\t{if}\ p\eq1\pmod4\ \&\ p=x^2+y^2\ (2\nmid x,\ 2\mid y),
\\0\pmod{p^2}&\t{if}\ p\eq3\pmod4.\endcases
\endalign$$
And
$$\align&\l(\f {-1}p\r)\sum_{k=0}^{p-1}\f{T_k(10,49)^3}{(-8)^{3k}}\eq\l(\f 6p\r)\sum_{k=0}^{p-1}\f{T_k(10,49)^3}{12^{3k}}
\\\eq&\cases 4x^2-2p\pmod{p^2}&\t{if}\ p\eq1,3\pmod8\ \&\ p=x^2+2y^2\ (x,y\in\Z),
\\0\pmod{p^2}&\t{if}\ (\f{-2}p)=-1,\ \t{i.e.},\ p\eq5,7\pmod8.\endcases
\endalign$$
Also,
$$\align\sum_{k=0}^{p-1}(7k+4)\f{T_k(10,49)^3}{(-8)^{3k}}\eq&\f p{14}\l(\f 2p\r)\l(65-9\l(\f p3\r)\r)\pmod{p^2},
\\\sum_{k=0}^{p-1}(7k+3)\f{T_k(10,49)^3}{12^{3k}}\eq&\f{3p}{28}\l(13+15\l(\f p3\r)\r)\pmod{p^2}.
\endalign$$
For each $n=1,2,3,\ldots$ we have
$$\sum_{k=0}^{n-1}(7k+4)T_k(10,49)^3(-8^3)^{n-1-k}\pmod{4n}$$
and
$$\sum_{k=0}^{n-1}(7k+3)T_k(10,49)^3(12^3)^{n-1-k}\pmod{n}.$$
\endproclaim

\proclaim{Conjecture A75 {\rm (Discovered on Oct. 2, 2011)}}
Let $p$ be an odd prime.

{\rm (i)} When $p>5$ we have
$$\align&\sum_{k=0}^{p-1}\l(\f{T_k(38,21^2)}{(-16)^k}\r)^3\eq\l(\f{-5}p\r)\sum_{k=0}^{p-1}\l(\f{T_k(38,21^2)}{20^k}\r)^3
\\\eq&\cases 4x^2-2p\pmod{p^2}&\t{if}\ (\f 2p)=(\f p3)=(\f p5)=1\ \&\ p=x^2+30y^2,
\\12x^2-2p\pmod{p^2}&\t{if}\ (\f p3)=1,\ (\f 2p)=(\f p5)=-1\ \&\ p=3x^2+10y^2,
\\2p-8x^2\pmod{p^2}&\t{if}\ (\f 2p)=1,\ (\f p3)=(\f p5)=-1\ \&\ p=2x^2+15y^2,
\\20x^2-2p\pmod{p^2}&\t{if}\ (\f p5)=1,\ (\f 2p)=(\f p3)=-1\ \&\ p=5x^2+6y^2,
\\\p da_{p,7}\pmod{p^2}&\t{if}\ (\f{-30}p)=-1.\endcases
\endalign$$
Also,
$$\sum_{k=0}^{p-1}(28k+15)\f{T_k^3(38,21^2)}{(-16)^{3k}}\eq \f p7\l(124-19\l(\f p3\r)\r)\pmod {p^2}.$$

{\rm (i)} When $p\not=7$ we have
$$\align&\sum_{k=0}^{p-1}\l(\f{T_k(110,57^2)}{32^k}\r)^3\eq\l(\f{-14}p\r)\sum_{k=0}^{p-1}\l(\f{T_k(110,57^2)}{(-28)^k}\r)^3
\\\eq&\cases4x^2-2p\pmod{p^2}&\t{if}\ (\f{-2}p)=(\f p3)=(\f p{7})=1\ \&\ p=x^2+42y^2\ (x,y\in\Z),
\\8x^2-2p\pmod{p^2}&\t{if}\ (\f{p}7)=1,\,(\f {-2}p)=(\f p{3})=-1\ \&\ p=2x^2+21y^2\ (x,y\in\Z),
\\12x^2-2p\pmod{p^2}&\t{if}\ (\f{-2}{p})=1,\,(\f p3)=(\f p7)=-1\ \&\ p=3x^2+14y^2\ (x,y\in\Z),
\\24x^2-2p\pmod{p^2}&\t{if}\ (\f{p}3)=1,\,(\f {-2}p)=(\f p{7})=-1\ \&\ p=6x^2+7y^2\ (x,y\in\Z),
\\p\da_{p,19}\pmod{p^2}&\t{if}\ (\f{-42}p)=-1.\endcases
\endalign$$
Also,
$$\sum_{k=0}^{p-1}(684k+329)\f{T_k^3(110,57^2)}{2^{15k}}\eq \f p{19}\l(5160\l(\f{-2}p\r)+1091\r)\pmod{p^2}.$$
\endproclaim

\proclaim{Conjecture A76 {(\rm Discovered on Oct. 4, 2011)}} Let $p>3$ be a prime. Then
$$\align&\sum_{k=0}^{p-1}\f{\bi{2k}kT_k^2(3,-3)}{(-108)^k}
\\\eq&\cases 4x^2-2p\pmod{p^2}&\t{if}\ (\f{-1}p)=(\f{p}3)=(\f p{7})=1\ \&\ p=x^2+21y^2,
\\12x^2-2p\pmod{p^2}&\t{if}\ (\f{-1}p)=(\f p7)=-1,\ (\f p3)=1\ \&\ p=3x^2+7y^2,
\\2x^2-2p\pmod{p^2}&\t{if}\ (\f {-1}p)=(\f p{3})=-1,\ (\f p{7})=1\ \&\ 2p=x^2+21y^2,
\\6x^2-2p\pmod{p^2}&\t{if}\ (\f {-1}p)=1,\ (\f p{3})=(\f p7)=-1,\ \&\ 2p=3x^2+7y^2,
\\0\pmod{p^2}&\t{if}\ (\f{-21}p)=-1.\endcases
\endalign$$
We also have
$$\sum_{k=0}^{p-1}\f{56k+19}{(-108)^k}\bi{2k}kT_{k}^2(3,-3)
\eq \f p2\l(21\l(\f{p}7\r)+17\r)\pmod{p^2}.$$
\endproclaim

\proclaim{Conjecture A77 {(\rm Discovered on Oct. 5, 2011)}} Let $p$ be an odd prime.

{\rm (i)} We have
$$\align&\sum_{k=0}^{p-1}\f{\bi{2k}kT_k^2(7,12)}{4^k}
\\\eq&\cases 4x^2-2p\ (\mo\ p^2)&\t{if}\ p\eq1\ (\mo\ 12)\ \&\ p=x^2+9y^2\ (x,y\in\Z),
\\-4xy(\f{xy}3)\ (\mo\ p^2)&\t{if}\ p\eq5\ (\mo\ 12)\ \&\ p=x^2+y^2\ (x,y\in\Z),
\\0\ (\mo\ p^2)&\t{if}\ p\eq3\ (\mo\ 4).
\endcases\endalign$$
If $p\not=3$, then
$$\sum_{k=0}^{p-1}\f{\bi{2k}kT_k^2(7,12)}{4^k}
\eq\sum_{k=0}^{p-1}\f{\bi{2k}kT^2_{2k}(3,3)}{36^k}\ \ \l(\mo\ p^{(5+(\f{-1}p))/2}\r).$$

{\rm (ii)} We have
$$\align&\sum_{k=0}^{p-1}\f{\bi{2k}kT_{2k}^2(9,20)}{4^k}
\\\eq&\cases 4x^2-2p\ (\mo\ p^2)&\t{if}\ p\eq 1,9\ (\mo\ 20)\ \&\ p=x^2+25y^2\ (x,y\in\Z),
\\(\f x5)4xy\ (\mo\ p^2)&\t{if}\ p\eq 13,17\ (\mo\ 20)\ \&\ p=x^2+y^2\ (2\nmid x,\ 5\mid xy-1),
\\0\ (\mo\ p^2)&\t{if}\ p\eq3\ (\mo\ 4).
\endcases
\endalign$$
If $p\not=11$, then
$$\sum_{k=0}^{p-1}\f{\bi{2k}kT_{2k}^2(9,20)}{4^k}\eq\sum_{k=0}^{p-1}\f{\bi{2k}kT_k^2(19,-20)}{22^{2k}}\pmod{p^2}.$$
\endproclaim
\Remark. Note that $T_k(7,12)=D_k(3)$ and $T_k(9,20)=D_k(4)$ for all $k=0,1,2,\ldots$.

\proclaim{Conjecture A78 {\rm ([S-10])}} Let $p>3$ be a prime.

{\rm (i)} We have
$$\align\sum_{k=0}^{p-1}\f{\bi{2k}k}{12^k}T_k\eq&\l(\f p3\r)\f{3^{p-1}+3}4\pmod{p^2},
\\\sum_{k=0}^{(p-1)/2}\f{\bi{2k}k}{16^k}T_{2k}(4,1)\eq&1\pmod{p^2},
\\\sum_{k=0}^{(p-1)/2}\f{C_k}{16^k}T_{2k}(4,1)\eq&\f 43\l(\l(\f 3p\r)-p\l(\f{-1}p\r)\r)\pmod{p^2},
\\\sum_{k=0}^{(p-1)/2}\f{\bi{2k}k}{4^k}T_{2k}(3,4)\eq&\l(\f{-1}p\r)\f{7-3^p}4\pmod{p^2},
\\\sum_{k=0}^{(p-1)/2}\f{\bi{2k}k}{16^k}T_{2k}(8,9)\eq&\l(\f 3p\r)\pmod{p^2}.
\endalign$$

{\rm (ii)} We have
$$\align&\sum_{k=0}^{(p-1)/2}\f{\bi{2k}k}{16^k}T_{2k}(2,3)\eq\sum_{k=0}^{p-1}\f{\bi{2k}k}{16^k}T_{2k}(4,-3)
\\\eq&\cases(\f{-1}p)(2x-\f p{2x})\pmod{p^2}&\t{if}\ p\eq1\pmod3\ \&\ p=x^2+3y^2\ (3\mid x-1),
\\0\pmod p&\t{if}\ p\eq2\pmod3.\endcases
\endalign$$
Also,
$$\sum_{k=0}^{(p-1)/2}\f{\bi{2k}k}{4^k}T_{2k}(1,-3)\eq\cases(-1)^{xy/2}2x\pmod p&\t{if}\ p=x^2+3y^2\ (3\mid x-1),
\\0\pmod{p}&\t{if}\ p\eq2\pmod 3;\endcases$$
and
$$\align&\sum_{k=0}^{(p-1)/2}\f{\bi{2k}k}{16^k}T_{2k}(4,3)
\\\eq&\cases(-1)^{(p-1)/4-\lfloor x/6\rfloor}2x\pmod{p}&\t{if}\ p\eq1\pmod{12}\ \&\ p=x^2+3y^2\ (4\mid x-1),
\\(-1)^{y/2-1}(\f{xy}3)2y\pmod{p}&\t{if}\ p\eq5\pmod{12}\ \&\ p=x^2+3y^2\ (4\mid x-1),
\\0\pmod{p}&\t{if}\ p\eq3\pmod4.\endcases
\endalign$$
\endproclaim

\proclaim{Conjecture A79} {\rm (i) ([S-12])} Let $p>3$ be a prime. Then
$$\align&\sum_{k=0}^{p-1}(-1)^k\bi{2k}k^2T_k
\\\eq&\cases 4x^2-2p\pmod{p^2}&\t{if}\ p\eq1,4\pmod{15}\ \&\ p=x^2+15y^2\ (x,y\in\Z),
\\2p-12x^2\pmod{p^2}&\t{if}\ p\eq2,8\pmod{15}\ \&\ p=3x^2+5y^2\ (x,y\in\Z),
\\0\pmod{p^2}&\t{if}\ (\f p{15})=-1,\ \t{i.e.},\ p\eq 7,11,13,14\pmod{15}.\endcases
\endalign$$
And
$$\sum_{k=0}^{p-1}(105k+44)(-1)^k\bi{2k}k^2T_k\eq p\l(20+24\l(\f p3\r)(2-3^{p-1})\r)\pmod{p^3}.$$
Also, we have
$$\f1{2n\bi{2n}n}\sum_{k=0}^{n-1}(-1)^{n-1-k}(105k+44)\bi{2k}k^2T_k\in\Z^+\quad\t{for all}\ n=1,2,3,\ldots.$$

{\rm (ii) ([S-10])} If $p\eq1,4\pmod{15}$ and $p=x^2+15y^2$ with $x,y\in\Z$, then
$$P_{(p-1)/2}(7\sqrt{-15}\pm16\sqrt{-3})\eq\l(\f{-\sqrt{-15}}p\r)\l(\f x{15}\r)\l(2x-\f p{2x}\r)\pmod{p^2},$$
where $P_n(x)$ denotes the Legendre polynomial $\sum_{k=0}^n\bi nk\bi{n+k}k(\f{x-1}2)^k$.

{\rm (iii) ([S-10])} Let $p>5$ be a prime. If $p\eq1,4\ (\mo\ 15)$ and $p=x^2+15y^2\ (x,y\in\Z)$ with $x\eq1\ (\mo\ 3)$, then
$$\align\sum_{k=0}^{p-1}\f{k\bi{2k}k\bi{3k}k}{27^k}F_k\eq&\f 2{15}\l(\f px-2x\r)\pmod{p^2},
\\\sum_{k=0}^{p-1}\f{\bi{2k}k\bi{3k}k}{27^k}L_k\eq& 4x-\f px\pmod{p^2}
\endalign$$
and
$$\sum_{k=0}^{p-1}(3k+2)\f{\bi{2k}k\bi{3k}k}{27^k}L_k\eq 4x\pmod{p^2},$$
where $F_0,F_1,F_2,\ldots$ are Fibonacci numbers and $L_0,L_1,L_2,\ldots$ are Lucas numbers.
If $p\eq2,8\ (\mo\ 15)$ and $p=3x^2+5y^2\ (x,y\in\Z)$ with $y\eq1\ (\mo\ 3)$, then
$$\sum_{k=0}^{p-1}\f{\bi{2k}k\bi{3k}k}{27^k}F_k\eq\f p{5y}-4y\pmod{p^2}$$
and
$$\sum_{k=0}^{p-1}\f{k\bi{2k}k\bi{3k}k}{27^k}F_k\eq\sum_{k=0}^{p-1}\f{k\bi{2k}k\bi{3k}k}{27^k}L_k\eq\f 43y\pmod{p^2}.$$
\endproclaim
\Remark. Let $p>5$ be a prime. By the theory of binary quadratic forms (cf. [C]),
if $p\eq1,4\ (\mo\ 15)$ then $p=x^2+15y^2$ for some $x,y\in\Z$; if $p\eq 2,8\ (\mo\ 15)$ then
$p=5x^2+3y^2$ for some $x,y\in\Z$. The author could show that  for any prime $p>3$ we have
$$\sum_{k=0}^{p-1}\f{\bi{2k}k\bi{3k}k}{27^k}F_k\eq0\ (\mo\ p^2)\ \ \t{if}\ \ p\eq1\ (\mo\ 3),$$
and
$$\sum_{k=0}^{p-1}\f{\bi{2k}k\bi{3k}k}{27^k}L_k\eq0\ (\mo\ p^2)\ \ \t{if}\ \ p\eq2\ (\mo\ 3).$$

\proclaim{Conjecture A80} Let $p>3$ be a prime.

 {\rm (i) ([S-10])}  We have
 $$\align&\sum_{k=0}^{p-1}\f{\bi{2k}k^2T_k(10,1)}{(-64)^k}\eq
 \l(\f p3\r)\sum_{k=0}^{p-1}\f{\bi{2k}k^2T_{2k}(6,1)}{256^k}\eq\sum_{k=0}^{p-1}\f{\bi{2k}k^2T_{2k}(6,1)}{1024^k}
 \\\eq&\cases
 4x^2-2p\ (\mo\ p^2)&\t{if}\ p\eq1,7\ (\mo\ 24)\ \&\ p=x^2+6y^2\ (x,y\in\Z),
 \\8x^2-2p\ (\mo\ p^2)&\t{if}\ p\eq5,11\ (\mo\ 24)\ \&\ p=2x^2+3y^2\ (x,y\in\Z),
 \\0\ (\mo\ p^2)&\t{if}\ (\f {-6}p)=-1,\ i.e.,\ p\eq 13,17,19,23\ (\mo\ 24);
 \endcases\endalign$$
 Also,
 $$\sum_{k=0}^{p-1}(3k+1)\f{\bi{2k}k^2T_k(10,1)}{(-64)^k}\eq\f
 p4\l(3\l(\f p3\r)+1\r)\pmod{p^2}$$
 and
$$\sum_{k=0}^{n-1}(3k+1)\bi{2k}k^2T_k(10,1)(-64)^{n-1-k}\eq0\ \l(\mo\ 2n\bi{2n}n\r)$$
for all $n=2,3,\ldots$. If $(\f{-6}p)=1$, then
 $$\sum_{k=0}^{p-1}(16k+5)\f{\bi{2k}k^2T_{2k}(6,1)}{256^k}\eq\f
 83p\l(\f p3\r)\pmod{p^2}$$
 and
$$\sum_{k=0}^{p-1}(16k+3)\f{\bi{2k}k^2T_{2k}(6,1)}{1024^k}\eq
 -\f23p\pmod{p^2}.$$

 {\rm (ii) (Discovered on Sept. 29, 2011)} If $p=x^2+6y^2$ with $x\eq1\ (\mo\ 3)$, then
$$\align \sum_{k=0}^{p-1}\f{\bi{2k}k\bi{3k}k}{108^k}ku_k(4,2)\eq&\f 16\l(2x-\f px\r)\pmod{p^2},
\\\sum_{k=0}^{p-1}\f{\bi{2k}k\bi{3k}k}{108^k}v_k(4,2)\eq&4x-\f px\pmod{p^2},
\\\sum_{k=0}^{p-1}(3k-1)\f{\bi{2k}k\bi{3k}k}{108^k}v_k(4,2)\eq&-2x\pmod{p^2}.
\endalign$$
If $p=2x^2+3y^2$ with $x\eq1\ (\mo\ 3)$, then
$$\align \sum_{k=0}^{p-1}\f{\bi{2k}k\bi{3k}k}{108^k}u_k(4,2)\eq&2x-\f p{4x}\pmod{p^2},
\\\sum_{k=0}^{p-1}\f{\bi{2k}k\bi{3k}k}{108^k}ku_k(4,2)\eq&\f x3\pmod{p^2},
\\\sum_{k=0}^{p-1}\f{\bi{2k}k\bi{3k}k}{108^k}kv_k(4,2)\eq&\f 43x\pmod{p^2}.
\endalign$$
\endproclaim

\proclaim{Conjecture 81 {\rm ([S-10])}} Let $p>7$ be a prime.

{\rm (i)} We have
$$\align&\sum_{k=0}^{p-1}\f{\bi{2k}k\bi{3k}kT_k(18,1)}{512^k}
\eq\l(\f{10}p\r)\sum_{k=0}^{p-1}\f{\bi{2k}k\bi{3k}kT_{3k}(6,1)}{(-512)^k}
\\\eq&
\cases x^2-2p\pmod{p^2}&\t{if}\ (\f p5)=(\f p7)=1\ \&\ 4p=x^2+35y^2,
\\ 2p-5x^2\pmod{p^2}&\t{if}\ (\f p5)=(\f p7)=-1\ \&\ 4p=5x^2+7y^2,
\\0\pmod{p^2}&\t{if}\ (\f p{35})=-1.\endcases
\endalign$$
And
$$\align\sum_{k=0}^{p-1}(35k+9)\f{\bi{2k}k\bi{3k}kT_k(18,1)}{512^k}
\eq&\f{9p}2\l(7-5\l(\f p5\r)\r)\pmod{p^2},
\\\sum_{k=0}^{p-1}(35k+9)\f{\bi{2k}k^2T_{3k}(6,1)}{(-512)^k}\eq&\f {9p}{32}\l(\f 2p\r)\l(25+7\l(\f{p}7\r)\r)\pmod{p^2}.
\endalign$$

{\rm (ii)} Suppose that $(\f p5)=(\f p{7})=1$ and write $4p=x^2+35y^2$ with $x,y\in\Z$. If $p\eq1\ (\mo\ 3)$, then
$$\sum_{k=0}^{p-1}\f{\bi{2k}{k}\bi{3k}k}{3456^k}(64+27\sqrt{5}\pm \sqrt{-35})^k\eq\l(\f x3\r)\l(2x-\f p{2x}\r)\pmod{p^2}.$$
If $p\eq2\ (\mo\ 3)$, then
$$\sum_{k=0}^{p-1}\f{\bi{2k}{k}\bi{3k}k}{3456^k}(64+27\sqrt{5}\pm \sqrt{-35})^k\eq \pm\sqrt{-35}\l(\f y3\r)\l(y-\f p{35y}\r)\pmod{p^2}.$$
\endproclaim

\proclaim{Conjecture 82 {\rm ([S-10])}} Let $p>3$ be a prime. Then
$$\align&\sum_{k=0}^{p-1}\f{\bi{2k}k^2 T_{2k}(18,1)}{256^k}
\\\eq&\cases 4x^2-2p\pmod{p^2}&\t{if}\ (\f 2p)=(\f p3)=(\f p5)=1\ \&\ p=x^2+30y^2,
\\12x^2-2p\pmod{p^2}&\t{if}\ (\f p3)=1,\ (\f 2p)=(\f p5)=-1\ \&\ p=3x^2+10y^2,
\\2p-8x^2\pmod{p^2}&\t{if}\ (\f 2p)=1,\ (\f p3)=(\f p5)=-1\ \&\ p=2x^2+15y^2,
\\2p-6x^2\pmod{p^2}&\t{if}\ (\f p5)=1,\ (\f 2p)=(\f p3)=-1\ \&\ 2p=3x^2+10y^2,
\\0\pmod{p^2}&\t{if}\ (\f{-30}p)=-1.\endcases
\endalign$$
And
$$\align&\sum_{k=0}^{p-1}\f{\bi{2k}k^2 T_{2k}(30,1)}{256^k}
\\\eq&\cases 4x^2-2p\pmod{p^2}&\t{if}\ (\f {-2}p)=(\f p3)=(\f p7)=1\ \&\ p=x^2+42y^2,
\\12x^2-2p\pmod{p^2}&\t{if}\ (\f {-2}p)=1,\ (\f p3)=(\f p7)=-1\ \&\ p=3x^2+14y^2,
\\2p-8x^2\pmod{p^2}&\t{if}\ (\f p7)=1,\ (\f {-2}p)=(\f p3)=-1\ \&\ p=2x^2+21y^2,
\\2p-6x^2\pmod{p^2}&\t{if}\ (\f p3)=1,\ (\f {-2}p)=(\f p7)=-1\ \&\ 2p=3x^2+14y^2,
\\0\pmod{p^2}&\t{if}\ (\f{-42}p)=-1.\endcases
\endalign$$
\endproclaim

\proclaim{Conjecture 83 {\rm ([S-10])}} Let $p>3$ be a prime. When $p\not=13,17$, we have
$$\align&\sum_{k=0}^{p-1}\f{\bi{2k}k\bi{3k}kT_k(102,1)}{102^{3k}}
\\\eq&\cases 4x^2-2p\pmod{p^2}&\t{if}\ (\f 2p)=(\f p3)=(\f p{13})=1\ \&\ p=x^2+78y^2,
\\2p-8x^2\pmod{p^2}&\t{if}\ (\f 2p)=1,\ (\f p3)=(\f p{13})=-1\ \&\ p=2x^2+39y^2,
\\12x^2-2p\pmod{p^2}&\t{if}\ (\f p{13})=1,\ (\f 2p)=(\f p3)=-1\ \&\ p=3x^2+26y^2,
\\2p-24x^2\pmod{p^2}&\t{if}\ (\f p3)=1,\ (\f 2p)=(\f p{13})=-1\ \&\ p=6x^2+13y^2,
\\0\pmod{p^2}&\t{if}\ (\f{-78}p)=-1.\endcases
\endalign$$
Provided $p\not=11,17$, we have
$$\align&\sum_{k=0}^{p-1}\f{\bi{2k}k\bi{3k}kT_k(198,1)}{198^{3k}}
\\\eq&\cases 4x^2-2p\pmod{p^2}&\t{if}\ (\f {2}p)=(\f p3)=(\f p{17})=1\ \&\ p=x^2+102y^2,
\\2p-8x^2\pmod{p^2}&\t{if}\ (\f p{17})=1,\ (\f {2}p)=(\f p{3})=-1\ \&\ p=2x^2+51y^2,
\\12x^2-2p\pmod{p^2}&\t{if}\ (\f p{3})=1,\ (\f {2}p)=(\f p{17})=-1\ \&\ p=3x^2+34y^2,
\\2p-24x^2\pmod{p^2}&\t{if}\ (\f {2}p)=1,\ (\f p3)=(\f p{17})=-1\ \&\ p=6x^2+17y^2,
\\0\pmod{p^2}&\t{if}\ (\f{-102}p)=-1.\endcases
\endalign$$
\endproclaim

\proclaim{Conjecture 84 {\rm ([S-10])}} Let $p\not=2,5,19$ be a prime. We have
$$\align&\sum_{k=0}^{p-1}\f{\bi{2k}k^2T_{2k}(5778,1)}{1216^{2k}}
\\\eq&\cases 4x^2-2p\pmod{p^2}&\t{if}\ (\f{2}p)=(\f{p}5)=(\f p{19})=1\ \&\ p=x^2+190y^2,
\\8x^2-2p\pmod{p^2}&\t{if}\ (\f{2}p)=1,\ (\f p5)=(\f p{19})=-1\ \&\ p=2x^2+95y^2,
\\2p-20x^2\pmod{p^2}&\t{if}\ (\f 2p)=(\f p{5})=-1,\ (\f p{19})=1\ \&\ p=5x^2+38y^2,
\\2p-40x^2\pmod{p^2}&\t{if}\ (\f 2p)=(\f p{19})=-1,\ (\f p{5})=1\ \&\ p=10x^2+19y^2,
\\0\pmod{p^2}&\t{if}\ (\f{-190}p)=-1,\endcases
\endalign$$
and
$$\sum_{k=0}^{p-1}(57720k+24893)\f{\bi{2k}k^2T_{2k}(5778,1)}{1216^{2k}}
\eq p\l(11548+13345\l(\f p{95}\r)\r)\pmod{p^2}.$$
Provided $p\not=17$ we have
$$\sum_{k=0}^{p-1}\f{\bi{2k}k^2T_{2k}(5778,1)}{439280^{2k}}
\eq\l(\f p5\r)\sum_{k=0}^{p-1}\f{\bi{2k}k^2T_{2k}(5778,1)}{1216^{2k}}\pmod{p^2}$$
and
$$\sum_{k=0}^{p-1}(57720k+3967)\f{\bi{2k}k^2T_{2k}(5778,1)}{439280^{2k}}
\eq p\l(\f p{19}\r)\l(3983-16\l(\f p{95}\r)\r)\pmod{p^2}.$$
\endproclaim

\proclaim{Conjecture 85 {\rm ([S-10])}} Let $p>5$ be a prime. Then
$$\align&\sum_{k=0}^{p-1}\f{\bi{2k}k^2T_{2k}(198,1)}{224^{2k}}\eq\l(\f p7\r)\sum_{k=0}^{p-1}\f{\bi{2k}k^2T_{2k}(322,1)}{48^{4k}}
\\\eq&\cases 4x^2-2p\pmod{p^2}&\t{if}\ (\f 2p)=(\f p5)=(\f p7)=1\ \&\ p=x^2+70y^2,
\\8x^2-2p\pmod{p^2}&\t{if}\ (\f p7)=1,\ (\f 2p)=(\f p5)=-1\ \&\ p=2x^2+35y^2,
\\2p-20x^2\pmod{p^2}&\t{if}\ (\f p5)=1,\ (\f 2p)=(\f p7)=-1\ \&\ p=5x^2+14y^2,
\\28x^2-2p\pmod{p^2}&\t{if}\ (\f 2p)=1,\ (\f p5)=(\f p7)=-1\ \&\ p=7x^2+10y^2,
\\0\pmod{p^2}&\t{if}\ (\f{-70}p)=-1.\endcases
\endalign$$
Also,
$$\align&\sum_{k=0}^{p-1}\f{\bi{2k}k^2T_{2k}(322,1)}{(-2^{10}3^4)^k}
\\\eq&\cases 4x^2-2p\pmod{p^2}&\t{if}\ (\f {-1}p)=(\f p5)=(\f p{17})=1\ \&\ p=x^2+85y^2,
\\2p-2x^2\pmod{p^2}&\t{if}\ (\f {p}{17})=1,\ (\f {-1}p)=(\f p5)=-1\ \&\ 2p=x^2+85y^2,
\\2p-20x^2\pmod{p^2}&\t{if}\ (\f {-1}p)=1,\ (\f p5)=(\f p{17})=-1\ \&\ p=5x^2+17y^2,
\\10x^2-2p\pmod{p^2}&\t{if}\ (\f p{5})=1,\ (\f {-1}p)=(\f p{17})=-1\ \&\ 2p=5x^2+17y^2,
\\0\pmod{p^2}&\t{if}\ (\f{-85}p)=-1.\endcases
\endalign$$
And
$$\align&\sum_{k=0}^{p-1}\f{\bi{2k}k^2T_{2k}(1298,1)}{24^{4k}}
\\\eq&\cases 4x^2-2p\pmod{p^2}&\t{if}\ (\f {-2}p)=(\f p5)=(\f p{13})=1\ \&\ p=x^2+130y^2,
\\8x^2-2p\pmod{p^2}&\t{if}\ (\f {-2}p)=1,\ (\f p5)=(\f p{13})=-1\ \&\ p=2x^2+65y^2,
\\2p-20x^2\pmod{p^2}&\t{if}\ (\f p5)=1,\ (\f {-2}p)=(\f p{13})=-1\ \&\ p=5x^2+26y^2,
\\2p-40x^2\pmod{p^2}&\t{if}\ (\f p{13})=1,\ (\f {-2}p)=(\f p5)=-1\ \&\ p=10x^2+13y^2,
\\0\pmod{p^2}&\t{if}\ (\f{-130}p)=-1.\endcases
\endalign$$
\endproclaim

\proclaim{Conjecture A86 {\rm (Parts (i)-(ii) and (iii)-(iv) were discovered on August 16 and June 18, 2011 respectively)}}
 {\rm (i)} For any prime $p>3$ we have
$$\align&\sum_{k=0}^{p-1}\f{\bi{4k}{2k}\bi{2k}kT_k(2702,1)}{384^{2k}}
\\\eq&\cases4x^2-2p\ (\mo\ p^2)&\t{if}\ (\f{-1}p)=(\f p3)=(\f p{11})=1,\ p=x^2+33y^2,
\\2p-2x^2\ (\mo\ p^2)&\t{if}\ (\f{-1}p)=1,\ (\f p3)=(\f p{11})=-1,\ 2p=x^2+33y^2,
\\12x^2-2p\ (\mo\ p^2)&\t{if}\ (\f{p}{11})=1,\ (\f {-1}p)=(\f p{3})=-1,\ p=3x^2+11y^2,
\\2p-6x^2\ (\mo\ p^2)&\t{if}\ (\f{p}{3})=1,\ (\f {-1}p)=(\f p{11})=-1,\ 2p=3x^2+11y^2,
\\0\ (\mo\ p^2)&\t{if}\ (\f{-33}p)=-1,\endcases
\endalign$$
and also
$$\sum_{k=0}^{p-1}\f{220k+119}{384^{2k}}\bi{4k}{2k}\bi{2k}kT_k(2702,1)
\eq\f p2\l(\f p3\r)\l(55+183\l(\f{-1}p\r)\r)\ (\mo\ p^2).$$

{\rm (ii)} Let $p>3$ be a prime with $p\not=7$. Then
$$\align&\l(\f p7\r)\sum_{k=0}^{p-1}\f{\bi{4k}{2k}\bi{2k}kT_k(115598,1)}{2688^{2k}}
\\\eq&\cases4x^2-2p\ (\mo\ p^2)&\t{if}\ (\f{-1}p)=(\f p3)=(\f p{19})=1,\ p=x^2+57y^2,
\\2p-2x^2\ (\mo\ p^2)&\t{if}\ (\f{-1}p)=1,\ (\f p3)=(\f p{19})=-1,\ 2p=x^2+57y^2,
\\12x^2-2p\ (\mo\ p^2)&\t{if}\ (\f{p}{3})=1,\ (\f {-1}p)=(\f p{19})=-1,\ p=3x^2+19y^2,
\\2p-6x^2\ (\mo\ p^2)&\t{if}\ (\f{p}{19})=1,\ (\f {-1}p)=(\f p{3})=-1,\ 2p=3x^2+19y^2,
\\0\ (\mo\ p^2)&\t{if}\ (\f{-57}p)=-1,\endcases
\endalign$$
and also
$$\sum_{k=0}^{p-1}\f{260k+513}{2688^{2k}}\bi{4k}{2k}\bi{2k}kT_k(115598,1)
\eq\f p2\l(\f {21}p\r)\l(961+65\l(\f{-1}p\r)\r)\ (\mo\ p^2).$$

{\rm (iii)} Let $p>3$ be a prime. Then
$$\align&\sum_{n=0}^{p-1}\bi{2n}n\sum_{k=0}^n\f{\bi nk\bi{n+2k}{2k}\bi{2k}k}{64^k}
\\\eq&\cases x^2-2p\pmod{p^2}&\t{if}\ (\f p3)=(\f p{17})=1\ \&\ 4p=x^2+51y^2\ (x,y\in\Z),
\\3x^2-2p\pmod{p^2}&\t{if}\ (\f p3)=(\f p{17})=-1\ \&\ 4p=3x^2+17y^2\ (x,y\in\Z),
\\0\pmod{p^2}&\t{if}\ (\f p{51})=-1,\endcases
\endalign$$
and also
$$\sum_{n=0}^{p-1}(17n+9)\bi{2n}n\sum_{k=0}^n\f{\bi nk\bi{n+2k}{2k}\bi{2k}k}{64^k}\eq\f p3\l(34\l(\f p{17}\r)-7\r)\pmod{p^2}.$$

{\rm (iv)} Let $p>5$ be a prime. Then
$$\align&\sum_{n=0}^{p-1}\f{\bi{2n}n}{400^n}\sum_{k=0}^n\bi nk\bi{n+2k}{2k}\bi{2k}k196^{n-k}
\\\eq&\cases x^2-2p\pmod{p^2}&\t{if}\ (\f 2p)=(\f p{11})=1\ \&\ p=x^2+22y^2\ (x,y\in\Z),
\\8x^2-2p\pmod{p^2}&\t{if}\ (\f 2p)=(\f p{11})=-1\ \&\ p=2x^2+11y^2\ (x,y\in\Z),
\\0\pmod{p^2}&\t{if}\ (\f {-22}p)=-1,\endcases
\endalign$$
and also
$$\sum_{n=0}^{p-1}\f{33n+19}{400^n}\bi{2n}n\sum_{k=0}^n\bi nk\bi{n+2k}{2k}\bi{2k}k196^{n-k}\eq19p\pmod{p^2}.$$
\endproclaim
\Remark. There are many similar congruences.
The quadratic fields $\Q(\sqrt{-33})$ and $\Q(\sqrt{-57})$ have class number four.

\proclaim{Conjecture A87 {\rm ([S-13])}} For any positive integer $n$, we have
$$8(2n-1)\bi{3n}n\ \bigg|\ \sum_{k=0}^n\bi{6k}{3k}\bi{3k}k\bi{6(n-k)}{3(n-k)}\bi{3(n-k)}{n-k}.$$
\endproclaim

\proclaim{Conjecture A88 {\rm ([S11h])}} Let $k$ and $l$ be positive integers.
If $(ln+1)\mid\bi{kn+ln}{kn}$ for all sufficiently large positive integers $n$, then
each prime factor of $k$ divides $l$. In other words, if $k$ has a prime factor not dividing $l$ then
there are infinitely many positive integers $n$ such that  $(ln+1)\nmid\bi{kn+ln}{kn}$.
\endproclaim
\Remark. The author [S11h] noted that if $k$ and $l$ are positive integers then
$\bi{kn+ln}{kn}\eq0\ (\mo\ (ln+1)/(k,ln+1))$ for all $n\in\Z^+$.
\medskip

\proclaim{Conjecture A89 {\rm ([S11g])}} For $n=0,1,2,\ldots$ set
$$s_n=\f{\bi{6n}{3n}\bi{3n}n}{2(2n+1)\bi{2n}n}.$$
Then, for any prime $p>3$ we have
 $$\sum_{k=1}^{p-1}\f{s_k}{108^k}\eq\cases 0\ (\mo\ p)&\t{if}\ p\eq\pm1\ (\mo\ 12),
 \\-1\ (\mo\ p)&\t{if}\ p\eq\pm5\ (\mo\ 12).\endcases$$
 Also, there are positive integers $t_1,t_2,t_3,\ldots$ such that
$$\sum_{k=0}^\infty s_kx^{2k+1}+\f1{24}-\sum_{k=1}^\infty t_kx^{2k}
=\f{\cos(\f23\arccos(6\sqrt 3x))}{12}$$
for all real $x$ with $|x|\ls 1/(6\sqrt3)$. Moreover, $t_p\eq-2\ (\mo\ p)$ for any prime $p$.
\endproclaim
\Remark. The author [S11g] showed that $s_n\in\Z$ for all $n=1,2,3,\ldots$.
Using {\tt Mathematica} the author found that
  $$\sum_{k=0}^\infty s_kx^k=\f{\sin(\f23\arcsin(6\sqrt{3x}))}{8\sqrt{3x}}\ \ \l(0<x<\f1{108}\r)$$
  and in particular
  $$\sum_{k=0}^\infty\f{s_k}{108^k}=\f{3\sqrt3}8.$$
\medskip

Recall that the Fibonacci sequence $\{F_n\}_{n\gs0}$ and the Lucas sequence $\{L_n\}_{n\gs0}$ are defined  by
$$F_0=0,\ F_1=1,\ \t{and}\ F_{n+1}=F_n+F_{n-1}\ (n=1,2,3,\ldots),$$
and
$$L_0=2,\ L_1=1,\ \t{and}\ L_{n+1}=L_n+L_{n-1}\ (n=1,2,3,\ldots).$$

\proclaim{Conjecture A90 {\rm  (Discovered on Oct. 27, 2010)}} Let $p\not=2,5$ be a prime and set $q:=F_{p-(\f p5)}/p$. Then
$$\align p\sum_{k=1}^{p-1}\f{F_{2k}}{k^2\bi{2k}k}\eq&-\l(\f p5\r)\l(\f 32q+\f 54p\, q^2\r)\ (\mo\ p^2),
\\ p\sum_{k=1}^{p-1}\f{L_{2k}}{k^2\bi{2k}k}\eq&-\f 52q-\f {15}4p\, q^2\ (\mo\ p^2),
\\\sum_{k=0}^{(p-3)/2}\f{\bi{2k}kF_{2k+1}}{(2k+1)16^k}\eq&(-1)^{(p+1)/2}\l(\f p5\r)\l(\f 12q+\f 58p\, q^2\r)\ (\mo\ p^2),
\\\sum_{k=0}^{(p-3)/2}\f{\bi{2k}kL_{2k+1}}{(2k+1)16^k}\eq&(-1)^{(p+1)/2}\l(\f 52q+\f 58p\, q^2\r)\ (\mo\ p^2).
\endalign$$
\endproclaim
\Remark. The conjecture was motivated by the following new identities observed by the author
on Oct. 27, 2010:
$$\gather \sum_{k=1}^\infty\f{F_{2k}}{k^2\bi{2k}k}=\f{4\pi^2}{25\sqrt 5},\ \ \sum_{k=1}^\infty\f{L_{2k}}{k^2\bi{2k}k}=\f{\pi^2}5,
\\\sum_{k=0}^\infty\f{\bi{2k}kF_{2k+1}}{(2k+1)16^k}=\f{2\pi}{5\sqrt5},\ \ \sum_{k=0}^\infty\f{\bi{2k}kL_{2k+1}}{(2k+1)16^k}=\f{2\pi}5.
\endgather$$
In fact, they can be obtained by putting $x=(\sqrt5\pm1)/2$ in the identities
$$\arcsin \f x2=\sum_{k=0}^\infty\f{\bi{2k}k}{(2k+1)4^k}\l(\f x2\r)^{2k+1}
\ \t{and}\ \ \sum_{k=1}^\infty\f{x^{2k}}{k^2\bi{2k}k}=2\arcsin^2\f x2.$$

\proclaim{Conjecture A91 {\rm ([S-5])}} For any $n\in\Z^+$ we have
$$\f{(-1)^{\lfloor n/5\rfloor-1}}{(2n+1)n^2\bi{2n}n}\sum_{k=0}^{n-1}F_{2k+1}\bi{2k}k\eq
\cases6\ (\mo\ 25)&\t{if}\ n\eq0\ (\mo\ 5),\\4\ (\mo\ 25)&\t{if}\ n\eq1\ (\mo\ 5),
\\1\ (\mo\ 25)&\t{if}\ n\eq2,4\ (\mo\ 5),\\9\ (\mo\ 25)&\t{if}\ n\eq3\ (\mo\ 5).\endcases$$
Also, if $a,b\in\Z^+$ and $a\gs b$ then the sum
$$\f1{5^{2a}}\sum_{k=0}^{5^a-1}F_{2k+1}\bi{2k}k$$
modulo $5^b$ only depends on $b$.
\endproclaim
\Remark. In [S-5] the author proved that if $p\not=2,5$ is a prime then
$$\sum_{k=0}^{p-1}F_{2k}\bi{2k}k\eq(-1)^{\lfloor
p/5\rfloor}\l(1-\l(\f p5\r)\r)\ (\mo\ p^2)$$ and
$$\sum_{k=0}^{p-1}F_{2k+1}\bi{2k}k\eq(-1)^{\lfloor p/5\rfloor}\l(\f p5\r)\ (\mo\ p^2).$$
\medskip

Recall that the usual $q$-analogue of $n\in\N$ is given by
$$[n]_q=\f{1-q^n}{1-q}=\sum_{0\ls k<n}q^k$$
which tends to $n$ as $q\to1$.
For any $n,k\in\N$ with $n\gs k$,
$$\M nk=\f{\prod_{0<r\ls n}[r]_q}{(\prod_{0<s\ls k}[s]_q)
(\prod_{0<t\ls n-k}[t]_q)}$$
is a natural extension of the usual binomial coefficient
$\bi nk$. A $q$-analogue of Fibonacci numbers introduced by I. Schur [Sc]
is defined as follows:
$$F_0(q)=0,\ F_1(q)=1,\ \t{and}\ F_{n+1}(q)=F_n(q)+q^nF_{n-1}(q)\ (n=1,2,3,\ldots).$$

\proclaim{Conjecture A92  {\rm ([S-5])}} Let $a$ and $m$ be positive integers.
 Then, in the ring $\Z[q]$, we have the following congruence
$$\sum_{k=0}^{5^am-1}q^{-2k(k+1)}\M{2k}kF_{2k+1}(q)\eq0\ (\mo\ [5^a]_q^2).$$
\endproclaim

\proclaim{Conjecture A93 {\rm ([S-5])}} For any $n\in\Z^+$ we have
$$\f{(-1)^{n-1}}{n^2(n+1)\bi{2n}n}\sum_{k=0}^{n-1}u_{k+1}(4,1)\bi{2k}k\eq
\cases1\ (\mo\ 9)&\t{if}\ n\eq0,2\ (\mo\ 9),\\4\ (\mo\ 9)&\t{if}\ n\eq5,6\ (\mo\ 9),
\\-2\ (\mo\ 9)&\t{otherwise}.\endcases$$
Also, if $a,b\in\Z^+$ and $a\gs b-1$ then the sum
$$\f1{3^{2a}}\sum_{k=0}^{3^a-1}u_{k+1}(4,1)\bi{2k}k$$
modulo $3^b$ only depends on $b$.
\endproclaim
\Remark. In [S-5] the author proved that if $p>3$ is a prime then
$$\sum_{k=0}^{p-1}u_k(4,1)\bi{2k}k\eq2\l(\l(\f p3\r)-\l(\f{-1}p\r)\r)\ (\mo\ p^2)$$
and
$$\sum_{k=0}^{p-1}u_{k+1}(4,1)\bi{2k}k\eq\l(\f p3\r)\ (\mo\ p^2).$$

\proclaim{Conjecture A94 {\rm ([S-7, S-8])}} For $k\in\N$ set $H_k^{(2)}=\sum_{0<j\ls k}1/j^2$.
If $p>3$ is a prime, then
$$\align \sum_{k=1}^{p-1}\f{\bi{2k}k}kH^{(2)}_k\eq&\f{2H_{p-1}}{3p^2}+\f{76}{135}p^2B_{p-5}\ (\mo\ p^3),
\\\sum_{k=1}^{p-1}\f{\bi{2k}kH_k^{(2)}}{k2^k}\eq&-\f3{16}\cdot\f{H_{p-1}}{p^2}+\f{479}{1280}p^2B_{p-5}\ (\mo\ p^3),
\\\sum_{k=1}^{p-1}\f{\bi{2k}kH_k^{(2)}}{k3^k}\eq&-\f 89\cdot\f{H_{p-1}}{p^2}+\f{268}{1215}p^2B_{p-5}\ (\mo\ p^3),
\\\sum_{k=1}^{p-1}\f{\bi{2k}k}{k4^k}H_k^{(2)}\eq&-\f 32\cdot\f{H_{p-1}}{p^2}+\f 7{80}p^2B_{p-5}\ (\mo\ p^3).
\endalign$$
\endproclaim
\Remark. {\tt Mathematica 7} yields that
$$\sum_{k=1}^\infty\f{2^kH_{k-1}^{(2)}}{k^2\bi{2k}k}=\f{\pi^4}{384}
\ \ \ \t{and}\ \ \ \sum_{k=1}^\infty\f{3^kH_{k-1}^{(2)}}{k^2\bi{2k}k}=\f{2\pi^4}{243}.$$
Also,
$$\sum_{k=1}^\infty\f{\bi{2k}k}{k4^k}H_k^{(2)}=\f32\zeta(3)\ \ \ \t{and}
\ \ \ \sum_{k=1}^\infty\f{4^kH_{k-1}^{(2)}}{k^2\bi{2k}k}=\f{\pi^4}{24}.$$

\proclaim{Conjecture A95 {\rm ([S11e])}} Let $p$ be an odd prime and
let $h\in\Z$ with $2h-1\eq0\ (\mo\ p)$. If $a\in\Z^+$ and $p^a>3$,
then
$$\sum_{k=0}^{p^a-1}\bi{hp^a-1}k\bi{2k}k\l(-\f h2\r)^k\eq0\ (\mo\ p^{a+1}).$$
Also, for any $n\in\Z^+$ we have
$$\nu_p\(\sum_{k=0}^{n-1}\bi{hn-1}k\bi{2k}k\l(-\f h2\r)^k\)\gs\nu_p(n).$$
\endproclaim

\proclaim{Conjecture A96 {\rm ([S11e])}} Let $m\in\Z$ with $m\eq1\
(\mo\ 3)$. Then
$$\nu_3\(\f1n\sum_{k=0}^{n-1}\bi{n-1}k(-1)^k\f{\bi{2k}k}{m^k}\)\gs \min\{\nu_3(n),\nu_3(m-1)\}-1$$
for every $n\in\Z^+$.  Furthermore,
 $$\f1{3^a}\sum_{k=0}^{3^a-1}\bi{3^a-1}k(-1)^k\f{\bi{2k}k}{m^k}\eq-\f{m-1}3\ (\mo\ 3^{\nu_3(m-1)})$$
for each integer $a>\nu_3(m-1)$. Also,
$$\sum_{k=0}^{3^a-1}\bi{3^a-1}k(-1)^k\bi{2k}k\eq -3^{2a-1}\ (\mo\ 3^{2a})\ \ \t{for every}\ a=2,3,\ldots.$$
 \endproclaim

\proclaim{Conjecture A97 {\rm ([S11a])}} For any odd prime $p$ and
positive integer $n$ we have
 $$\nu_p\(\sum_{k=0}^{n-1}\bi{(p-1)k}{k,\ldots,k}\)\gs\nu_p\l(n\bi{2n}n\r).$$
 \endproclaim

\Remark. The author [S11a] proved that an integer $p>1$ is a prime if and only if
$$\sum_{k=0}^{p-1}\bi{(p-1)k}{k,\ldots,k}\eq 0\ (\mo\ p).$$ He also showed that if $n\in\Z^+$
is a multiple of a prime $p$ then
$$\sum_{k=0}^{n-1}\bi{(p-1)k}{k,\ldots,k}\eq 0\ (\mo\ p).$$

\proclaim{Conjecture A98 {\rm ([S-3])}} Let $m\gs2$ and $r$ be integers.
And let $p>r$ be an odd prime not dividing $m$.

 {\rm (i)} If $m>2$, $m\not\eq r\ (\mo\ 2)$, and $p\eq r\ (\mo\ m)$ with $r\gs-m/2$, then
 $$\sum_{k=0}^{p-1}(-1)^{km}\bi{r/m}k^m\eq0\ (\mo\ p^3).$$

 {\rm (ii)} If $p\eq r\ (\mo\ 2m)$ with $r\gs-m$, then
 $$\sum_{k=0}^{p-1}(-1)^k\bi{r/m}k^{2n+1}\eq0\ (\mo\ p^2)\quad\t{for all}\ n=1,\ldots,m-1.$$

 {\rm (iii)} For any prime $p$ and positive integer $n$, we have
  $$\nu_p\(\sum_{k=0}^{n-1}\bi{-1/(p+1)}k^{p+1}\)\gs c_p\l\lfloor\f{\nu_p(n)+1}2\r\rfloor,$$
 where
 $$c_p=\cases1&\t{if}\ p=2,\\3&\t{if}\ p=3,\\5&\t{if}\ p\gs5.\endcases$$
\endproclaim
\Remark. The author [S-3] proved part (ii) in the case $n=1$.
He [S10j] also showed that for any prime $p>3$ we have
$$\sum_{k=0}^{p-1}\bi{-1/(p+1)}k^{p+1}\eq 0\ (\mo\ p^5)$$
and
$$\sum_{k=0}^{p-1}(-1)^{km}\bi{p/m-1}k^m\eq0\ (\mo\ p^4)\quad\t{for any integer}\ m\not\eq0\ (\mo\ p).$$
He  conjectured that there are no composite numbers $n$ satisfying the congruence
$$\sum_{k=0}^{n-1}\bi{-1/(n+1)}k^{n+1}\eq 0\ (\mo\ n^5).$$

\proclaim{Conjecture A99 {\rm (Discovered in 2007)}} Let $p$ be a prime and let $l,n\in\N$ and $r\in\Z$.
If $n$ or $r$ is not divisible by $p$ then we have
$$\align&\nu_p\(\sum_{k\eq r\,(\mo\ p)}\bi nk(-1)^k\bi{(k-r)/p}l\)
\\\gs&\l\lfloor\f{n-lp-1}{p-1}\r\rfloor
+\nu_p\(\bi{\lfloor(n-l-1)/(p-1)\rfloor}l\).
\endalign$$
\endproclaim
\Remark. D. Wan [W] proved that the inequality holds if the last term on the right-hand side
is omitted (see also Sun and Wan [SW]).

\proclaim{Conjecture A100} {\rm (i) (raised on Nov. 2, 2009 via a message to Number Theory List)}
If $n>1$ is an odd integer satisfying the Morley congruence
$$\bi{n-1}{(n-1)/2}\eq(-1)^{(n-1)/2}4^{n-1}\ (\mo\ n^3),$$
then $n$ must be a prime.

{\rm (ii) ([S10])} If an odd integer $n>1$ satisfies the congruence
$$\sum_{k=0}^{n-1}\f{\bi{2k}k}{2^k}\eq(-1)^{(n-1)/2}\ (\mo\ n^2),$$
then $n$ must be a prime.

{\rm (iii) ([S11j])} Set $a_n:=\sum_{k=0}^n\bi nk^2C_k$ for $n=0,1,2,\ldots$.
Any integer $n>1$ satisfying $a_1+\cdots+a_{n-1}\eq0\ (\mo\ n^2)$ must be a prime.
\endproclaim
\Remark. (a) In 1895 Morley [Mo] showed that $\bi{p-1}{(p-1)/2}\eq(-1)^{\f{p-1}2}4^{p-1}\, (\mo\ p^3)$
for any prime $p>3$. The author has verified part (i) of Conj. A100 for $n<10^4$.

(b) The author [S10] proved that if $p$ is an odd prime then
$$\sum_{k=0}^{p-1}\f{\bi{2k}k}{2^k}\eq(-1)^{(p-1)/2}\ (\mo\ p^2).$$
And he verified part (ii) of Conj. A100 for $n<10^4$ via Mathematica. On the author's request,
Qing-Hu Hou at Nankai Univ. finished the verification for $n<10^5$.

(c) The author [S11j] proved that $\sum_{k=1}^{p-1}a_k\eq0\ (\mo\ p^2)$ for any odd prime $p$,
and he also verified part (iii) of Conj. A100 for $n\ls70,000$.

 \bigskip\medskip

\heading{Part B. Conjectures that have been confirmed}\endheading

\proclaim{Conjecture B1 {\rm (raised in an early version of [S11e], and confirmed by
Kasper Andersen)}} For any positive integer $n$, the arithmetic mean
$$s_n:=\f1n\sum_{k=0}^{n-1}(21k+8)\bi{2k}k^3$$
is always an integer divisible by $4\bi{2n}n$.
\endproclaim
\Remark. The author created the sequence $\{s_n\}_{n\gs1}$
at OEIS as A173774 (cf. [S]).
On Feb. 11, 2010, Andersen proved the conjecture by noting that $t_n:=s_n/(4\bi{2n}n)$ coincides with
$$r_n:=\sum_{k=0}^{n-1}\bi{n+k-1}k^2.$$

\proclaim{Conjecture B2 {\rm (raised in [S11b], and confirmed by
Zhi-Hong Sun [Su1])}} Let $p$ be an odd prime. Then
$$\sum_{k=0}^{p-1}\l((-2)^{-k}-4^{-k}\r)\bi{2k}k^2\eq0\ (\mo\ p)$$
and
$$\sum_{k=0}^{p-1}\f{k\bi{2k}k^2}{16^k}\eq\f{(-1)^{(p+1)/2}}4\ (\mo\ p^2).$$
If $p\eq1\ (\mo\ 4)$ and
$p=x^2+y^2$ with $x\eq1\ (\mo\ 4)$ and $y\eq0\ (\mo\ 2)$, then
$$\sum_{k=0}^{p-1}\f{\bi{2k}k^2}{8^k}\eq\sum_{k=0}^{p-1}\f{\bi{2k}k^2}{(-16)^k}
\eq(-1)^{(p-1)/4}\l(2x-\f p{2x}\r)\ (\mo\ p^2)$$
and
$$\sum_{k=0}^{p-1}\f{\bi{2k}k^2}{32^k}\eq2x-\f p{2x}\ (\mo\ p^2).$$
If $p\eq 3\ (\mo\ 4)$ then
$$\sum_{k=0}^{p-1}\f{\bi{2k}k^2}{32^k}\eq0\ (\mo\ p^2).$$
\endproclaim
\Remark. The author proved those congruences modulo $p$ except the first one.

\proclaim{Conjecture B3 {\rm (raised in 2009, and confirmed by Roberto Tauraso)}} Let $p$ be an odd prime.
Then
$$\sum_{k=0}^{(p-1)/2}\f{C_k^2}{16^k}\eq12p^2-4\ (\mo\ p^3)\ \t{and}\
\sum_{k=0}^{(p-1)/2}\f{kC_k^2}{16^k}\eq4-10p^2\ (\mo\ p^3).$$
Also,
$$\sum_{k=0}^{p-1}\f{\bi{4k}{2k}C_k}{64^k}\eq p\ (\mo\ p^2),$$
and
$$\sum_{k=0}^{(p-1)/2}\f{\bi{4k}{2k}C_k}{64^k}\eq(-1)^{(p-1)/2}\ \f23p\ (\mo\ p^2)\ \ \t{provided}\ p>3.$$
If $p>3$, then
$$\sum_{k=0}^{p-1}\f{\bi{3k}kC_k}{27^k}\eq p\ (\mo\ p^2)\ \t{and}\
\sum_{k=0}^{(p-1)/2}\f{\bi{3k}kC_k}{27^k}\eq\f p2\l(\f p3\r)\ (\mo\ p^2).$$
\endproclaim
\Remark. The author [S11b] showed that $\sum_{k=0}^{p-1}C_k^2/16^k\eq-3\ (\mo\ p)$ for any odd prime $p$,
and his PhD student Yong Zhang proved the first and the second congruences mod $p^2$. {\tt Mathematica} yields that
$$\sum_{k=2}^\infty\f{27^k}{(k-1)k^2\bi{3k}{k,k,k}}=\f{81}4-3\sqrt3\pi.$$

\proclaim{Conjecture B4 {\rm (raised in 2009, and confirmed by the author's PhD student Yong Zhang)}}
Let $p$ be an odd prime. Then
$$\sum_{k=0}^{(p-1)/2}\f{\bi{2k}{k+1}^2}{16^k}\eq (-1)^{(p-1)/2}-4+p^2(8+E_{p-3})\ (\mo\ p^3).$$
If $p>3$, then
$$\sum_{k=0}^{(p-1)/2}\f{C_kC_{k+1}}{16^k}\eq8\ (\mo\ p^2).$$
If $p\eq1\ (\mo\ 4)$, then
$$\sum_{k=0}^{p-1}\f{\bi{2k}k\bi{2k}{k+1}}{8^k}\eq0\ (\mo\ p).$$
If $p\eq3\ (\mo\ 4)$, then
$$\sum_{k=0}^{p-1}\f{C_kC_{k+1}}{(-16)^k}\eq-10\ (\mo\ p).$$
\endproclaim
\Remark. As for the first congruence in Conjecture B5, the author proved the congruence mod $p$
and then his PhD student Yong Zhang showed the congruence mod $p^2$. Following the author's recent method
in [S11e], Zhang confirmed the congruence with the help of the software {\tt Sigma}.

\proclaim{Conjecture B5 {\rm (raised in [S11a], and confirmed by the author's PhD student Yong Zhang)}}
Let $m\in\Z$ with $m\eq1\ (\mo\ 3)$. Then
 $$\nu_3\(\f1n\sum_{k=0}^{n-1}\f{\bi{2k}k}{m^k}\)\gs \min\{\nu_3(n),\nu_3(m-1)-1\}$$
for every $n\in\Z^+$.  Furthermore,
 $$\f1{3^a}\sum_{k=0}^{3^a-1}\f{\bi{2k}k}{m^k}\eq\f{m-1}3\ (\mo\ 3^{\nu_3(m-1)})$$
 for any integer $a\gs\nu_3(m-1)$.
 \endproclaim

\proclaim{Conjecture B6 {\rm (raised in [S11e], and confirmed by Hao Pan and the author [PS])}}
Let $p$ be an odd prime and let $a\in\Z^+$.
If $p\eq1\ (\mo\ 4)$ or $a>1$, then
$$\sum_{k=0}^{\lfloor\f34p^a\rfloor}\f{\bi{2k}k}{(-4)^k}\eq\l(\f2{p^a}\r)\ (\mo\ p^2).$$
\endproclaim

\proclaim{Conjecture B7 {\rm (raised in [S11e], and confirmed by the author's former student Hui-Qin Cao)}}
If $p$ is a prime
with $p\eq11\ (\mo\ 12)$, then
$$\sum_{k=0}^{p-1}\l(\f k3\r)\f{\bi{2k}k^2}{(-16)^k}\eq0\ (\mo\ p).$$
\endproclaim

\proclaim{Conjecture B8 {\rm (raised in [S11c], and confirmed by R. Me\v strovi\'c [Me])}} For any prime $p>3$, we have
$$\sum_{k=1}^{p-1}\f{H_k^2}{k^2}\eq\f 45pB_{p-5}\ (\mo\ p^2).$$
\endproclaim
\Remark. Motivated by the identity $\sum_{k=1}^\infty H_k^2/k^2=17\pi^4/360$,
the author [S11c] proposed the conjecture and proved that $\sum_{k=1}^{p-1}H_k^2/k^2\eq0\pmod p$ for any prime $p>5$.
Other conjectures in [S11c] were confirmed by the author and Li-Lu Zhao [S-Z].
\medskip

\proclaim{Conjecture B9 {\rm (raised in [S11d], and confirmed by Hui-Qin Cao and the author [CS])}}
 Let $p>3$ be a prime. Then
$$\sum_{k=0}^{p-1}\bi{p-1}k\bi{2k}k((-1)^k-(-3)^{-k})\eq\l(\f p3\r)(3^{p-1}-1)\ (\mo\ p^3).$$
If $p\eq\pm1\ (\mo\ 12)$, then
$$\sum_{k=0}^{p-1}\bi{p-1}k\bi{2k}k(-1)^ku_k(4,1)\eq(-1)^{(p-1)/2}u_{p-1}(4,1)\ (\mo\ p^3).$$
If $p\eq\pm1\ (\mo\ 8)$, then
$$\sum_{k=0}^{p-1}\bi{p-1}k\bi{2k}k\f{u_k(4,2)}{(-2)^k}\eq(-1)^{(p-1)/2}u_{p-1}(4,2)\ (\mo\ p^3).$$
\endproclaim
\Remark. Note that
$$u_k(4,2)=\cases2^{k/2}P_k&\t{if}\ k\ \t{is even},\\2^{(k-3)/2}Q_k&\t{if}\ k\ \t{is odd}.\endcases$$
The author [S-4] showed that
$$\sum_{k=0}^{(p-1)/2}\f{u_k(4,2)}{16^k}\bi{2k}k\eq\f{(-1)^{\lfloor (p-4)/8\rfloor}}2\l(1-\l(\f 2p\r)\r)\ (\mo\ p^2)$$
and
$$\sum_{k=0}^{(p-1)/2}\f{v_k(4,2)}{16^k}\bi{2k}k\eq2(-1)^{\lfloor p/8\rfloor}\l(\f{-1}p\r)\ (\mo\ p^2).$$

\medskip

Recall that the $n$th Bell number $b_n$ denote the number of partitions of a set of cardinality $n$.
Bell numbers are also given by $b_0=1$ and $b_{n+1}=\sum_{k=0}^n\bi nk b_k\ (n=1,2,3,\ldots)$.

\proclaim{Conjecture B10 {\rm (discovered on July 17, 2010, and confirmed by the author and D. Zagier [SZ])}}
For any positive integer $n$ there is a unique
integer $s(n)$ such that for any prime $p$ not dividing $n$ we have
$$\sum_{k=0}^{p-1}\f{b_k}{(-n)^k}\eq s(n)\ (\mo\ p).$$
In particular,
$$\gather s(2)=1,\ s(3)=2,\ s(4)=-1,\ s(5)=10,\ s(6)=-43,
\\\ s(7)=266,\ s(8)=-1853,\ s(9)=14834,\ s(10)=-133495.\endgather$$
\endproclaim
\Remark. The author and D. Zagier [SZ] proved that $s(n)=(-1)^{n-1}D_{n-1}+1$ for all $n=1,2,3,\ldots$,
where $D_m$ denotes the derangement number $m!\sum_{k=0}^m(-1)^k/k!$
(the number of fixed-point-free permutations of a set of cardinality $m$).

\medskip

\proclaim{Conjecture B11 {\rm (raised in an early version of [S-6],
and actually confirmed by F. Rodriguez-Villegas [RV1])}} Let $p$ be an odd prime and
let $m\eq 4\ (\mo\ p)$. Then
$$\nu_p\(\sum_{k=0}^n\f{\bi{2k}k}{m^k}\)\gs \nu_p\l((2n+1)\binom{2n}n\r)\quad\t{for any}\ n\in\Z^+.$$
 Moreover, if $p>3$ then
$$\f1{p^a}\sum_{k=0}^{(p^a-1)/2}\f{\bi{2k}k}{m^k}\eq(-1)^{(p^a-1)/2}\ (\mo\ p).$$
\endproclaim

\proclaim{Conjecture B12 {\rm (raised in [S11d], and confirmed
by Z. H. Sun [Su2])}} Let $p$ be an odd prime. Then
$$\align\sum_{k=0}^{p-1}\f{P_k}{(-8)^k}\bi{2k}k^2\eq&0\ (\mo\ p)\ \ \ \t{if}\ p\eq 5\ (\mo\ 8),
\\\sum_{k=0}^{p-1}\f{P_k}{32^k}\bi{2k}k^2\eq&0\ (\mo\ p)\ \ \ \t{if}\ p\eq 7\ (\mo\ 8),
\\\sum_{k=0}^{p-1}\f{Q_k}{(-8)^k}\bi{2k}k^2\eq&0\ (\mo\ p)\ \ \ \t{if}\ p\eq5,7\ (\mo\ 8),
\\\sum_{k=0}^{p-1}\f{Q_k}{32^k}\bi{2k}k^2\eq&0\ (\mo\ p)\ \ \ \t{if}\ p\eq 5\ (\mo\ 8).
\endalign$$
Also,
$$\gather\sum_{k=0}^{p-1}\f{u_k(4,1)}{4^k}\bi{2k}k^2\eq0\ (\mo\ p)\ \ \ \t{if}\ p\eq 2\ (\mo\ 3),
\\\sum_{k=0}^{p-1}\f{u_k(4,1)}{64^k}\bi{2k}k^2\eq0\ (\mo\ p)\ \ \ \t{if}\ p\eq 11\ (\mo\ 12),
\\\sum_{k=0}^{p-1}\f{v_k(4,1)}{4^k}\bi{2k}k^2\eq\sum_{k=0}^{p-1}\f{v_k(4,1)}{64^k}\bi{2k}k^2\eq0\ (\mo\ p)
 \ \t{if}\ p\eq 5\ (\mo\ 12).
 \endgather$$
\endproclaim

\proclaim{Conjecture B13 {\rm (raised in [S11e], and confirmed
by Z. H. Sun [Su2])}} Let $p$ be an odd prime. Then
$$\sum_{k=0}^{p-1}\f{\bi{2k}k^3}{(-8)^k}\eq\cases 4x^2-2p\ (\mo\ p^2)&\t{if}\ 4\mid p-1\ \&\ p=x^2+y^2\ (2\nmid x),
\\0\ (\mo\ p^2)&\t{if}\ p\eq3\ (\mo\ 4).\endcases$$
Also, when $p\eq1\ (\mo\ 4)$ we have
$$\sum_{k=0}^{(p-1)/2}\f{4k+1}{64^k}\bi{2k}k^3\eq0\ (\mo\ p^2).$$
\endproclaim

\proclaim{Conjecture B14 {\rm (raised in [S11b, S11d, S11e, S-8], and confirmed by S. Mattarei and R. Tauraso [MT])}}
Let $p>3$ be a prime. Then
$$\align\sum_{k=1}^{p-1}\f{(-2)^k}{k^2}\bi{2k}k\eq&-2q_p(2)^2\ (\mo\ p),
\\ p\sum_{k=1}^{p-1}\f{2^k}{k\bi{2k}k}\eq& \l(\f{-1}p\r)-1-p\,q_p(2)+p^2E_{p-3}\ (\mo\ p^3),
\\p\sum_{k=1}^{p-1}\f{2^k}{k^2\bi{2k}k}\eq&-q_p(2)+\f{p^2}{16}B_{p-3}\ (\mo\ p^3),
\\p\sum_{k=1}^{p-1}\f{3^k}{k^2\bi{2k}k}\eq&-\f 32q_p(3)+\f 49p^2B_{p-3}\ (\mo\ p^3),
\\\sum_{k=1}^{p-1}\f{4^k}{k^2\bi{2k}k}+\f{4q_p(2)}p\eq&-2q^2_p(2)+pB_{p-3}\ (\mo\ p^2),
\\\sum_{k=1}^{p-1}(-2)^k\bi{2k}kH_k^{(2)}\eq&\f23q_p(2)^2\ (\mo\ p),
\\\sum_{k=1}^{p-1}(-1)^k\bi{2k}kH_k^{(2)}\eq&\f 52\l(\f p5\r)\f {F^2_{p-(\f p5)}}{p^2}\ (\mo\ p)\ \ \ \t{if}\ p>5.
\endalign$$
\endproclaim
\Remark. It is known that
$$\sum_{k=1}^\infty\f{2^k}{k\bi{2k}k}=\f{\pi}2,\ \ \sum_{k=1}^\infty\f{2^k}{k^2\bi{2k}k}=\f{\pi^2}8,
\ \ \sum_{k=1}^\infty\f{3^k}{k^2\bi{2k}k}=\f29\pi^2.$$
Glaisher (cf. [Ma])
got the formula $\sum_{k=1}^\infty 4^k/(k^2\binom{2k}k)=\pi^2/2$. Let $p>3$ be a prime.
During March 6-7, 2010 the author [S11e] showed that
$$\sum_{k=1}^{(p-1)/2}\f{4^k}{k^2\bi{2k}k}\eq(-1)^{(p-1)/2}\,4E_{p-3}\ (\mo\ p).$$

\proclaim{Conjecture B15 {\rm (raised in [S11e], and confirmed
by J. Guillera [G6])}} We have the identity
$$\sum_{k=1}^\infty\f{(11k-3)64^k}{k^3\bi{2k}k^2\bi{3k}k}=8\pi^2.$$
\endproclaim

\proclaim{Conjecture B16 {\rm (raised in [S11e], and confirmed by Zhi-Hong Sun [Su4])}}

{\rm (i)} Let $p$ be a prime with
$p\eq1,3\ (\mo\ 8)$. Write $p=x^2+2y^2$ with $x\eq1\ (\mo\ 4)$. Then
 $$\sum_{k=0}^{p-1}\f{\bi{4k}{2k}\bi{2k}k}{128^k}\eq (-1)^{\lfloor(p+5)/8\rfloor}\l(2x-\f p{2x}\r)\ (\mo\ p^2).$$

{\rm (ii)} Let $p\eq1\ (\mo\ 4)$ be a prime.
Write $p=x^2+y^2$ with $x\eq1\ (\mo\ 4)$ and $y\eq0\ (\mo\ 2)$. Then
$$\sum_{k=0}^{p-1}\f{\bi{6k}{3k}\bi{3k}k}{864^k}\eq\cases(-1)^{\lfloor x/6\rfloor}(2x-p/(2x))\ (\mo\ p^2)&\t{if}\ p\eq1\ (\mo\ 12),
\\(\f{xy}3)(2y-p/(2y))\ (\mo\ p^2)&\t{if}\ p\eq5\ (\mo\ 12).
\endcases$$
\endproclaim
\Remark. The author [S-9] proved that
$\sum_{k=0}^{p-1}\bi{4k}{2k}\bi{2k}k/128^k\eq0\ (\mo\ p^2)$
for any prime $p\eq5,7\ (\mo\ 8)$,
and $\sum_{k=0}^{p-1}k\bi{4k}{2k}\bi{2k}k/128^k\eq 0\ (\mo\ p^2)$ for any prime $p\eq1,3\ (\mo\ 8)$.
He also showed that $\sum_{k=0}^{p-1}\bi{6k}{3k}\bi{3k}k/864^k\eq0\ (\mo\ p^2)$ for any prime $p\eq3\ (\mo\ 4)$, and
$\sum_{k=0}^{p-1}k\bi{6k}{3k}\bi{3k}k/864^k\eq0\ (\mo\ p^2)$ for any prime $p\eq1\ (\mo\ 4)$.
The author conjectured that
$$\f1{5^{a}}\sum_{k=0}^{5^a-1}\f{k\bi{6k}{3k}\bi{3k}k}{864^k}\eq75\pmod{125}\quad\t{for all}\ a=1,2,3,\ldots.$$

\proclaim{Conjecture B17 {\rm (raised in an early version of [S11j] and confirmed by Guo and Zeng [GZ2])}}
\ For any $n\in\Z^+$ and $x\in\Z$ we have
$$\sum_{k=0}^{n-1}(2k+1)(-1)^k A_k(x)\eq0\ (\mo\ n).$$
If $p$ is an odd prime, then
$$\sum_{k=0}^{p-1}(2k+1)(-1)^kA_k(x)\eq p\l(\f {1-4x}p\r)\ (\mo\ p^2).$$
For any prime $p>3$ we have
$$\sum_{k=0}^{p-1}(2k+1)(-1)^kA_k\eq p\l(\f p3\r)\ (\mo\ p^3)$$
and
$$\sum_{k=0}^{p-1}(2k+1)(-1)^kA_k(-2)\eq p-\f 43p^2q_p(2)\ (\mo\ p^3).$$
\endproclaim

\proclaim{Conjecture B18 {\rm (raised on April 4, 2010 (see [S-7]) and confirmed by K. Hessami Pilehrood and T. Hessami Pilehrood [HP2])}}
{\rm (i)} For any $n\in\Z^+$ we have
$$a_n:=\f1{8n^2\bi{2n}n^2}\sum_{k=0}^{n-1}(205k^2+160k+32)(-1)^{n-1-k}\bi{2k}k^5\in\Z^+.$$

{\rm (ii)} Let $p>5$ be a prime. Then
$$\sum_{k=0}^{p-1}(205k^2+160k+32)(-1)^k\bi{2k}k^5\eq 32p^2+64p^3H_{p-1}\ (\mo\ p^7).$$
\endproclaim
\Remark. The conjecture was motivated by the identity
$$\sum_{k=1}^\infty\f{(-1)^k(205k^2-160k+32)}{k^5\bi{2k}k^5}=-2\zeta(3)$$
obtained by T. Amdeberhan and D. Zeilberger [AZ].
The author also conjectured that
$$\sum_{k=0}^{(p-1)/2}(205k^2+160k+32)(-1)^k\bi{2k}k^5\eq 32p^2+\f{896}3p^5B_{p-3}\ (\mo\ p^6)$$
for any prime $p>3$; this is still open.

\proclaim{Conjecture B19} {\rm (i) (raised in [S11e] and confirmed by K. Hessami Pilehrood and T. Hessami Pilehrood [HP2])} We have
$$\sum_{k=1}^\infty\f{(15k-4)(-27)^{k-1}}{k^3\bi{2k}k^2\bi{3k}k}=\sum_{k=1}^\infty\f{(\f k3)}{k^2}.$$

{\rm (ii) (raised in [S11e] and confirmed by K. Hessami Pilehrood and T. Hessami Pilehrood [HP3])} For each prime $p>5$ we have
$$\sum_{k=1}^{p-1}\f{21k-8}{k^3\bi{2k}k^3}+\f{p-1}{p^3}
\eq\f{H_{p-1}}{p^2}(15p-6)+\f{12}5p^2B_{p-5}\ (\mo\ p^3).$$
\endproclaim

\proclaim{Conjecture B20 {\rm (raised in [S11j, S11k] and confirmed by Hao Pan [P])}}
Let $m>1$ and $n$ be positive integers and let $x$ be any integer. Then
$$\f1n\sum_{k=0}^{n-1}(2k+1)D_k(x)^m,\ \f1n\sum_{k=0}^{n-1}(2k+1)A_k(x)^m,\ \f1n\sum_{k=0}^{n-1}(2k+1)(-1)^kA_k(x)^m$$
are all integral.
\endproclaim
\Remark. Let $n>0$ and $x$ be integers. The author [S11j] proved that $\f1n\sum_{k=0}^{n-1}(2k+1)D_k(x)\in\Z$ and $\f1n\sum_{k=0}^{n-1}(2k+1)A_k(x)\in\Z$.
Guo and Zeng [GZ2] confirmed the author's conjecture that $\f1n\sum_{k=0}^{n-1}(2k+1)(-1)^kA_k(x)\in\Z$. In [S11j] the author showed that
$n^2\mid \sum_{k=0}^{n-1}(2k+1)D_k(x)^2$.

\enddocument